\documentclass[11pt,a4paper]{article}
\usepackage[latin1]{inputenc}
\usepackage{amsmath}
\usepackage{graphicx}
\usepackage{multirow}
\usepackage{amsfonts}
\usepackage{variations}
\usepackage{dsfont}
\usepackage{amssymb}
\usepackage[bottom]{footmisc}

\author{BERNARD Damien}
\title{\textbf{Small first zeros of $L$-functions}}
\date{}
\def\N{\mathbb{N}}

\def\R{\mathbb{R}}

\def\E{\mathbb{E}}
\def\V{\mathbb{V}}
\def\Z{\mathbb{Z}}
\newcommand{\un}{\mathds{1}}

\def\c{\mathfrak{c}}

\def\m{\mathfrak{m}}
\def\h{\mathfrak{h}}
\numberwithin{equation}{section}

\newtheorem{proposition}{Proposition}

\newtheorem{corollary}{Corollary}
\newtheorem{theorem}{Theorem}
\newcounter{remark}
\newenvironment{remark}{\addtocounter{remark}{1} \textbf{Remark \theremark}}{}
\newtheorem{lemma}{Lemma}

\begin{document}
\maketitle

\begin{abstract} From a family of $L$-functions with unitary symmetry, Hughes and Rudnick obtained results on the height of its lowest zero. We extend their results to other families of $L$-functions according to the type of symmetry coming from statistics for low-lying zeros. 
\end{abstract}

\tableofcontents
 
\def\thefootnote{\fnsymbol{footnote}}

\footnotetext{ D. BERNARD, Université Blaise Pascal, Laboratoire de Mathématiques, Campus des Cézeaux, BP 80026, 63171 Aubière cedex, France  \\
E-mail address: damien.bernard@math.univ-bpclermont.fr
\vspace{1mm}

Mathematics Subject Classification : 11M41, 11M50, 11F67 .

\vspace{1mm}
Key words : $L$-function, density theorem, random matrix theory, smallest zero, differential equation with temporal shifts.}

\def\thefootnote{\arabic{footnote}}

\pagebreak

\section{Introduction} 

  \subsection{Preview of results} 
The existence of a deep link between non-trivial zeros of natural families of $L$-functions and eigenvalues of random matrices has been speculated since Montgomery's work (\cite{Mo}) in the seventies. So, we are able to assign a classical compact group of matrices to many  classical families of $L$-functions. We can refer to  \cite{ILS}, \cite{FI}, \cite{HR}, \cite{RR} or \cite{Mi}. Using the one-level density, Hughes and Rudnick obtained informations about the lowest zero of Dirichlet $L$-functions (\cite{HR}, section 8). Our aim is to generalise these results. 

\paragraph*{Extreme low-lying zeros of a natural family of $L$-functions}
 Let $\mathcal{F}(Q)$ be a finite set of $L$-functions with analytic conductor $Q$. We build the associate family $\mathcal{F} = \bigcup_{Q \geq 1}\mathcal{F}(Q)$. We assume Riemann hypothesis for any function in $\mathcal{F}$. We also assume the density theorem for $\mathcal{F}$ with test functions $\Phi$ satisfying supp $\widehat{\Phi}\subset \left[-\nu;\nu \right]$ with $\nu <\nu_{\max}(\mathcal{F})$. We write $W^*[\mathcal{F}]$ for the one-level density for non-trivial low-lying zeros of functions in $\mathcal{F}$. It turns out that the one-level densities that have been identified up to now have always been of the shape $W[G]$ whose Fourier transform is given by
 \begin{equation}  \label{eqn:expressiontransformeefourierdensite} \widehat{W[G]}(y) = \delta_0(y)+\frac{\delta}{2}\eta(y)+\varepsilon
 \end{equation} 
 where $\delta_0$ is the Dirac function, $\eta$ is defined on $\R$ by 
 $$\eta(y)=\left\{\begin{array}{lcl}   1 &\mbox{ if }& |y|<1 \\ \frac{1}{2} &\mbox{ if }& |y|=1 \\ 0 &\mbox{ if }& |y|>1 \end{array} \right. $$ 
 and $\delta$ and $\varepsilon$ are given in table \ref{tab:tableaudeltaepsilon}. 
\begin{center}
\begin{table}[h]
\begin{center} 
\begin{tabular}{|c|c|c|c|c|c|} 
\hline
$G$ & $U$ & $O$ &$Sp$ & $SO^+$ & $SO^-$ \\
\hline 
$\delta$ & $0$ & $0$ & $-1$ & $1$ & $-1$ \\
\hline
$\varepsilon$ & $0$  & $1/2$ & $0$ & $0$ & $1$\\
\hline
\end{tabular}
\end{center}
\caption{\label{tab:tableaudeltaepsilon} Value of $(\delta,\varepsilon)$. }
 \end{table}
 \end{center}
Finally, for $G=Sp$, $SO^+$ or $SO^-$, let $n \geq 1$ be the only integer such that $n-1 < \nu_{\max}(\mathcal{F}) \leq n$ and consider the equation in $\lambda$ given by
\begin{small}
\begin{eqnarray} \label{eqn:equationaresoudreenoncedanslintro}
 \frac{\delta}{\lambda}\cos\theta_\lambda - \sum_{k=0}^{n-1}U_k(\lambda)\sin\left(\theta_\lambda -k\delta \frac{\pi}{2} \right)\left[\frac{\delta \alpha_R(k)}{2}-1+\varepsilon \beta_R(k)  \right]+\frac{2\varepsilon}{\lambda}\sum_{k=0}^{n-1}U_k(\lambda)\cos\left(\theta_\lambda -k\delta \frac{\pi}{2} \right) = 0 . 
 \end{eqnarray}
 \end{small}In this equation,  let $U_k$ be the k-th Chebyshev polynomial  of the second kind, $(\delta,\varepsilon)$ is given in the table \ref{tab:tableaudeltaepsilon} and parameters $\theta_\lambda$, $\alpha_R$ and $\beta_R$ are defined in lemma \ref{lem:defzlambda} (page \pageref{lem:defzlambda}) and proposition \ref{pro:prodonnantlequationgenerale} (page \pageref{pro:prodonnantlequationgenerale}). We prove the following theorem on the smallest non negative imaginary part $\tilde{\gamma}_{f,1}$ of a non-trivial normalised zero of $L(f,.)$ in $\mathcal{F}$.
 \begin{theorem} \label{th:thm1} We have
 \begin{multline*}  \limsup_{Q\rightarrow +\infty} \min_{L(f,.) \in \mathcal{F}(Q)} \tilde{\gamma}_{f,1}  \leq  \frac{1}{2\nu_{\max}(\mathcal{F})}\times \\
  \left\{ \begin{array}{lcl}               1 & \mbox{if} & W^*[\mathcal{F}]=W[U] \\
4V^{-1}\left( 1+\frac{2}{\nu_{\max}(\mathcal{F})} \right) & \mbox{if} &  W^*[\mathcal{F}]=W[O] \\
4V^{-1}\left( 1+(\delta+2\varepsilon)\frac{2}{\nu_{\max}(\mathcal{F})} \right)& \mbox{if} &  W^*[\mathcal{F}]=W[Sp], W[SO^+] \mbox{ or }W[SO^-] \mbox{ and } \nu_{\max}(\mathcal{F}) \leq 1 \\
\frac{\nu_{\max}(\mathcal{F})}{\pi} \underset{R\rightarrow \nu_{\max}(\mathcal{F})/2^-}{\lim} \lambda_R & \mbox{if} &  W^*[\mathcal{F}]=W[Sp], W[SO^+] \mbox{ or }W[SO^-] \mbox{ and } \nu_{\max}(\mathcal{F}) > 1
\end{array}\right.
\end{multline*}
where $\lambda_R$ is the smallest positive root of equation (\ref{eqn:equationaresoudreenoncedanslintro}) which is not a root of $U_nU_{n-1}$ (with $n-1 < \nu_{\max}(\mathcal{F}) \leq n$) and where $V$ is defined by
$$V : \left\{ \begin{array}{ccl} [0,\frac{1}{4}[\cup]\frac{1}{4},x_1[ & \longrightarrow & \R \\ x & \longmapsto & \frac{\tan (2\pi x) }{2\pi x}\end{array} \right. $$
with $x_1 =\inf\left\{ x>0, \, \frac{\tan(2\pi x)}{2\pi x}=1 \right\} \approx 0,71 $. 
 \end{theorem}
\begin{remark} We plot the upper bound in theorem \ref{th:thm1}, denoted $M_{W^*[\mathcal{F}]}(\nu_{\max}(\mathcal{F}))$, for each $W^*[\mathcal{F}]$. From top to bottom, we have $W^*[\mathcal{F}]=W[Sp]$,  $W^*[\mathcal{F}]=W[U]$, $W^*[\mathcal{F}]=W[SO^+]$, $W^*[\mathcal{F}]=W[O]$ and $W^*[\mathcal{F}]=W[SO^-]$. 
 \begin{center}
 \begin{tabular}{c}
  \includegraphics[scale=0.57]{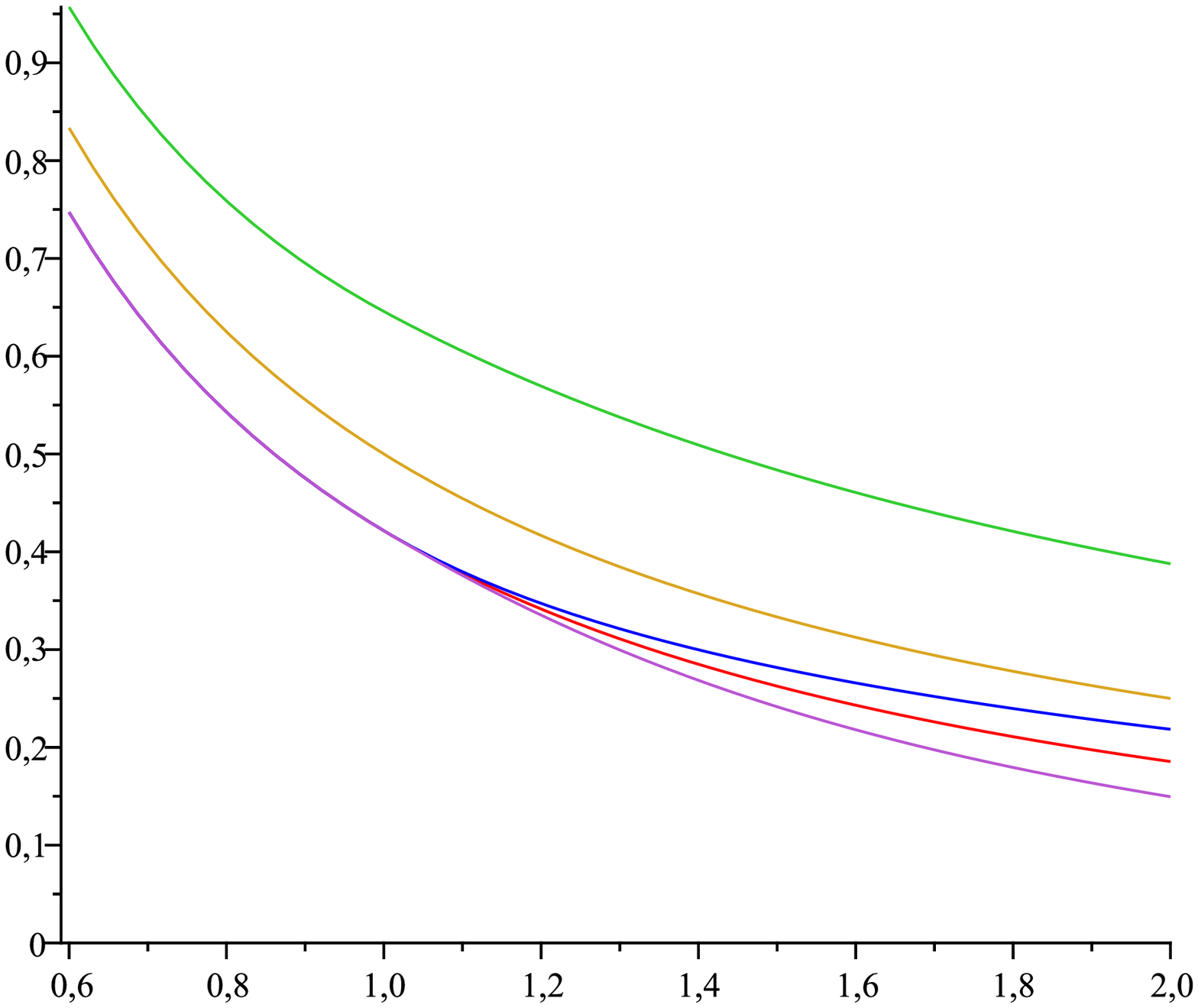}    \\
Representative curves of $\nu_{\max}(\mathcal{F}) \longmapsto M_{W^*[\mathcal{F}]}(\nu_{\max}(\mathcal{F})) $ 
\end{tabular}
\end{center}
\end{remark}
\begin{remark} In the orthogonal case, when $\nu_{\max}(\mathcal{F})$ goes to infinity, we have
\begin{center}
$   \frac{2}{\nu_{\max}(\mathcal{F})}V^{-1}\left( 1+ \frac{2}{\nu_{\max}(\mathcal{F})} \right)= \frac{\sqrt{6}}{\pi \nu_{\max}(\mathcal{F})^{3/2}}\left(1-\frac{6}{5\nu_{\max}(\mathcal{F})}+o\left(\frac{1}{\nu_{\max}(\mathcal{F})^{3/2}} \right) \right) .
$
\end{center}
\end{remark}
 \begin{center}
 \begin{tabular}{c}
\includegraphics[scale=0.3]{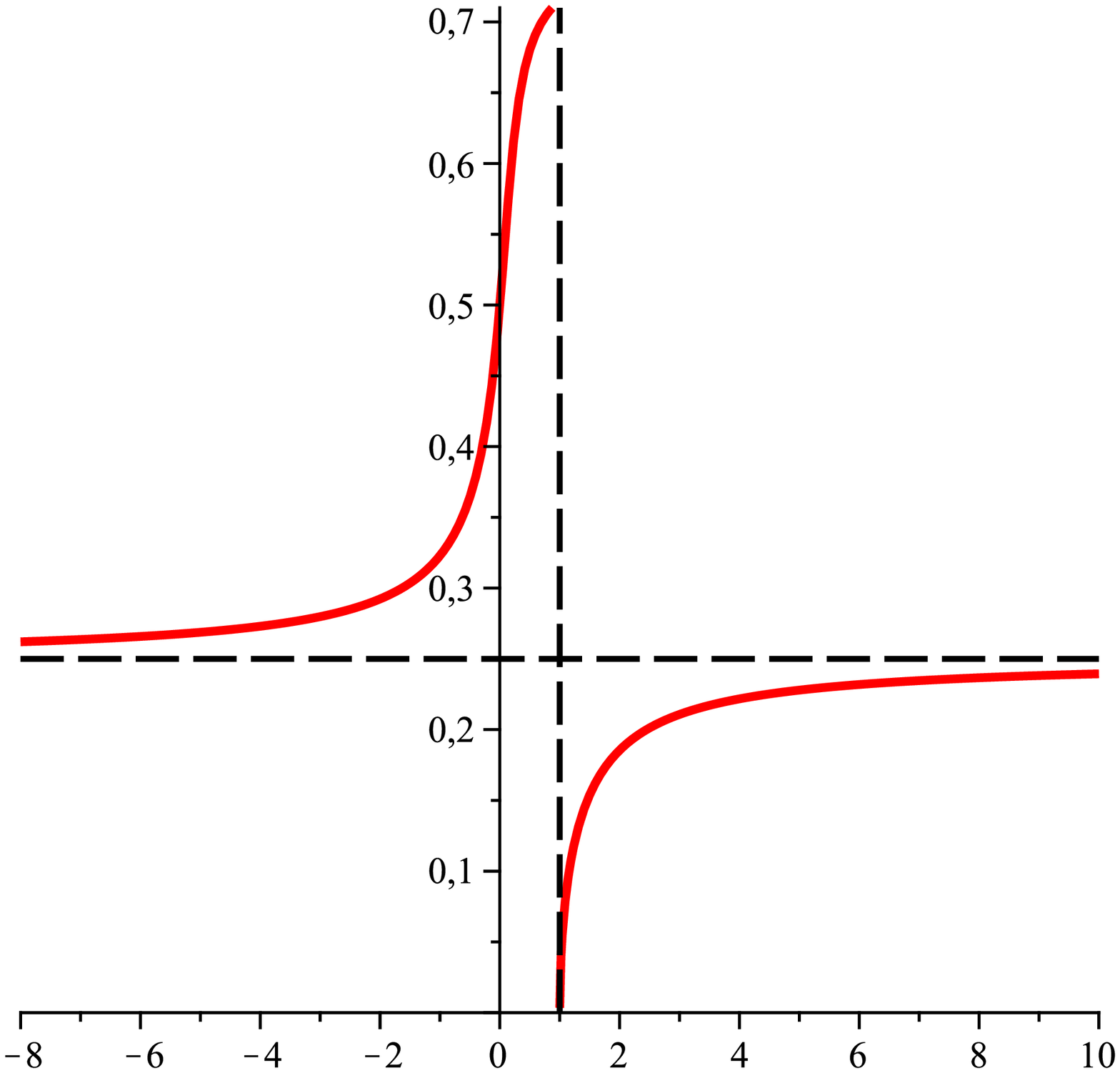} \\
 Representative curve of $V^{-1} $
\end{tabular}
\end{center}
\begin{remark} In the case $1 < \nu_{\max}(\mathcal{F}) \leq 2 $, we can simplify equation (\ref{eqn:equationaresoudreenoncedanslintro}) as follows
\begin{small}
\begin{eqnarray*}& &(\delta+2\varepsilon)\frac{1-4\lambda^2}{\lambda}\left(\sin\lambda(1-R)-2\delta\lambda \cos\lambda R\right)\\ & & -\left[(\delta+2\varepsilon)(1-R)-1+4\varepsilon \right]\left[\cos\lambda (1-R)-2\delta \lambda \sin\lambda R-2\lambda \tan \Theta_R (\sin\lambda(1-R)-2\delta\lambda \cos\lambda R)  \right]=0
\end{eqnarray*}
\end{small}where
\begin{small}
$$\Theta_R = \frac{1}{2}\left(R-\frac{1}{2}\right)+\frac{\pi}{2}\left(1+\frac{\delta}{2} \right) . $$
\end{small}Moreover, this relation always vanishes in $\lambda=1/2$ which is a root of $U_2$.
\end{remark}
\paragraph*{}
In order to give several instances, let \begin{small} $H_k^*(q)$ \end{small} be the set of primitive holomorphic cusp forms of prime level $q$ end even weight $k\geq 2$. Let $r$ be a positive integer. We define
\begin{small}
$$\mathcal{H}_r(q)=\left\{L(Sym^r f,s), \hspace{1mm}f\in H_k^*(q)\right\}$$ 
\end{small}and if $r$ is odd, 
  \begin{small}$$\mathcal{H}_r^{\pm}(q)=\left\{L(Sym^r f,s), \hspace{1mm}f\in H_k^*(q) \mbox{ and }\varepsilon(Sym^r f)=\pm1 \right\} .$$\end{small}Then, we have several families of symmetric power $L$-functions 
\begin{small}
$$\mathcal{H}_r=\underset{\scriptsize{q \mbox{ prime}}}{\bigcup}\mathcal{H}_r(q)\hspace{2mm} \mbox{ and }\hspace{2mm} \mathcal{H}_r^\pm =\underset{\scriptsize{q \mbox{ prime}}}{\bigcup}\mathcal{H}_r^\pm(q). $$
\end{small}These families have been studied in \cite{ILS} (theorem 1.1, for the case $r=1$) and \cite{RR} (theorem B). We may sum up some properties of these natural families of $L$-functions in the following table where
$$\nu_{1,max}(1,k,\theta_0) =2 \mbox{ and }  \nu_{1,max}(r,k,\theta_0) = \left( 1- \dfrac{1}{2(k-2\theta_0)}\right) \dfrac{2}{r^2} \mbox{ if } r\geq 2 \mbox{ with } \theta_0 = \frac{7}{64} ,$$
and for $\epsilon=\pm1$, 
$$\nu^{\epsilon}_{1,max}(1,k,\theta_0) =2 \mbox{ and }  \nu^{\epsilon}_{1,max}(r,k,\theta_0)= \inf \left\{ \nu_{1,max}(r,k,\theta_0),\frac{3}{r(r+1)} \right\}\mbox{ if } r\geq 2. $$
\begin{center}
\begin{table}[h]
\begin{small}
$$
\begin{array}{|c|c|c|c|c|}
\hline
\mathcal{F}     \rule[-4mm]{0cm}{10mm} & \mathcal{H}_r \hspace{2mm}(r\mbox{ even})  & \mathcal{H}_r \hspace{2mm}(r\mbox{ odd}) & \mathcal{H}_r^+  \hspace{2mm}(r\mbox{ odd}) & \mathcal{H}_r^-  \hspace{2mm}(r\mbox{ odd}) \\
\hline
\nu_{\max}(\mathcal{F}) \rule[-4mm]{0cm}{10mm}  &  \nu_{1,\max}(r,k,\theta_0) & \nu_{1,\max}(r,k,\theta_0) & \nu_{1,\max}^\epsilon(r,k,\theta_0) &  \nu_{1,\max}^\epsilon(r,k,\theta_0)\\
\hline
\rho_{\max}(\mathcal{F} ) \rule[-4mm]{0cm}{10mm}  & \frac{1}{r^2} & \frac{1}{r^2} & \frac{1}{2r(r+2)}  & \frac{1}{2r(r+2)}  \\
\hline
W[\mathcal{F}] \rule[-4mm]{0cm}{10mm}  & W[Sp] & W[O] & W[SO^+] & W[SO^-] \\
\hline
W^*[\mathcal{F}]  \rule[-4mm]{0cm}{10mm} & W[Sp] & W[O] & W[SO^+] & W[Sp] \\
\hline
\end{array}
$$
\end{small}
\caption{\label{tab:tableadenumax} Several natural families of $L$-functions }
 \end{table}
 \end{center}
Thus, thanks to theorem \ref{th:thm1}, we have
\begin{eqnarray*}
\limsup_{N \rightarrow +\infty} \min_{f\in H_{k}^{*}(N)}\tilde{\gamma}_{f,1} & \leq  & V^{-1}(2) < 0,19  \\
\limsup_{N \rightarrow +\infty} \min_{f\in H_{k}^{+}(N)}\tilde{\gamma}_{f,1} & \leq  &  0,22 \\
\limsup_{N \rightarrow +\infty} \min_{f\in H_{k}^{-}(N)}\tilde{\gamma}_{f,1} & \leq  & 0,39
\end{eqnarray*}
and
$$\limsup_{\underset{ q \mbox{ \begin{tiny}prime\end{tiny}}}{ q\rightarrow +\infty}} \min_{f \in H_{k}^{*}(q)}  \tilde{\gamma}_{Sym^rf,1} \leq \frac{2}{\nu_{1,max}(r,k,\theta_0)}V^{-1}\left( 1+(-1)^{r+1}\frac{2}{\nu_{1,max}(r,k,\theta_0)} \right)  $$
and if  $r$ is odd, $\epsilon=\pm 1$, we have 
$$\limsup_{\underset{ q \mbox{ \begin{tiny}prime\end{tiny}}}{ q\rightarrow +\infty}} \min_{ \begin{array}{c}  \scriptstyle{ f \in H_{k}^{*}(q)}  \\ \scriptstyle{\varepsilon (Sym^r f)=\epsilon} \end{array}}  \tilde{\gamma}_{Sym^rf,1} \leq \frac{2}{\nu^{\epsilon}_{1,max}(r,k,\theta_0)}V^{-1}\left( 1+\epsilon \frac{2}{\nu^{\epsilon}_{1,max}(r,k,\theta_0)} \right) .$$
Actually, since $\mathcal{H}^-(N) \neq \emptyset$ for large $N$, we get (see lemma \ref{lem:lemmezerotrivial} page \pageref{lem:lemmezerotrivial})
$\underset{\scriptsize{N \rightarrow +\infty}}{\lim} \underset{\scriptsize{f\in H_{k}^{*}(N)}}{\min}\tilde{\gamma}_{f,1}=0 .$

\paragraph*{Proportion of $L$-functions which have a small first zero}
Combining Bienaymé-Chebyshev inequality and statistics for low-lying zeros of symmetric power $L$-functions, we can obtain a positive proportion of $L$-functions in our family which have a small first zero. Hughes and Rudnick exposed this phenomenon in the case of Dirichlet $L$-functions which have a unitary symmetry. Symmetric power $L$-functions allows us to deal with all currently known symmetry group. $H_k^*(q)$ denotes the set of primitive holomorphic cusp forms of prime level $q$ and even weight $k\geq 2$ and let $\omega_q(f)$ be the harmonic weight associated to $f$ in $H_k^*(q)$. Let also $\varepsilon (Sym^r f)$ denotes the sign of the functional equation associated to the $L$-function $L(Sym^rf,s)$. 
\begin{theorem} \label{th:thm3} Let $r$ be a fixed positive integer. Assuming the Riemann hypothesis for all symmetric power $L$-functions of order $r$.
\begin{small}
$$\mbox{If }\beta \geq  \sqrt{\frac{\pi r^4}{4}\frac{ 6\pi^3r^4-24(-1)^{r}\pi r^2+9\pi-\pi^3+2\sqrt{6}\sqrt{-\pi^4+6\pi^4r^4+7\pi^2+12-24(-1)^{r}\pi^2r^2}   }{6\pi^4r^4+48(-1)^{r+1}\pi^2r^2+96-3\pi^2-\pi^4}} $$
\end{small}then
$$   \liminf_{\underset{ q \mbox{ \begin{tiny}prime\end{tiny}}}{ q\rightarrow +\infty}}  \sum_{\underset{\tilde{\gamma}_{Sym^rf,1}  \leq \beta}{f \in H_{k}^{*}(q)} }\omega_q(f)  \geq  1 - \frac{\pi^2}{6}\frac{16(\pi^2+3)\beta^4+8r^4(9-\pi^2)\beta^2+(3+\pi^2)r^8}{[r^6\pi^2 -4\beta^2r^2\pi^2+16\beta^2(-1)^r  ]^2} .  $$ 
Moreover, if $r$ is odd and $\sigma=\pm 1$. If $\frac{\beta}{r(r+2)\sqrt{\pi}} \geq $
\begin{small}
$$\sqrt{\frac{\pi[ 24\pi^2r^2(r+2)^2+48\sigma r(r+2)+9-\pi^2]+2\sqrt{6}\sqrt{-\pi^4+24\pi^4r^2(r+2)^2+7\pi^2+12+48\sigma\pi^2r(r+2)}   }{24\pi^4r^2(r+2)^2+96\sigma\pi^2r(r+2)+96-3\pi^2-\pi^4}} $$
then 
$$ \liminf_{\underset{ q \mbox{ \begin{tiny}prime\end{tiny}}}{ q\rightarrow +\infty}}  2\sum_{  \substack{   f \in H_{k}^{*}(q)  \\ \varepsilon (Sym^r f)=\sigma  \\ \tilde{\gamma}_{Sym^rf,1}  \leq \beta}   }\omega_q(f)  \geq 1-\frac{\pi^2}{24}\frac{(\pi^2+3)\beta^4+2r^2(r+2)^2(9-\pi^2)\beta^2+(\pi^2+3)r^4(r+2)^4}{[2\sigma\beta^2+\pi^2\beta^2r(r+2)-\pi^2r^3(r+2)^3]^2} .
$$
\end{small}
\end{theorem}
In order to give some examples for small values of $r$, we plot the graph of the function with parameter $\beta$ which is associated to the lower bound in the first part of theorem \ref{th:thm3}. On these representative curves, the critical value of $\beta$ corresponds to the minimal value of $\beta$ which appears in theorem \ref{th:thm3}. 
\begin{center}
\begin{tabular}{cccc}
\includegraphics[scale=0.182]{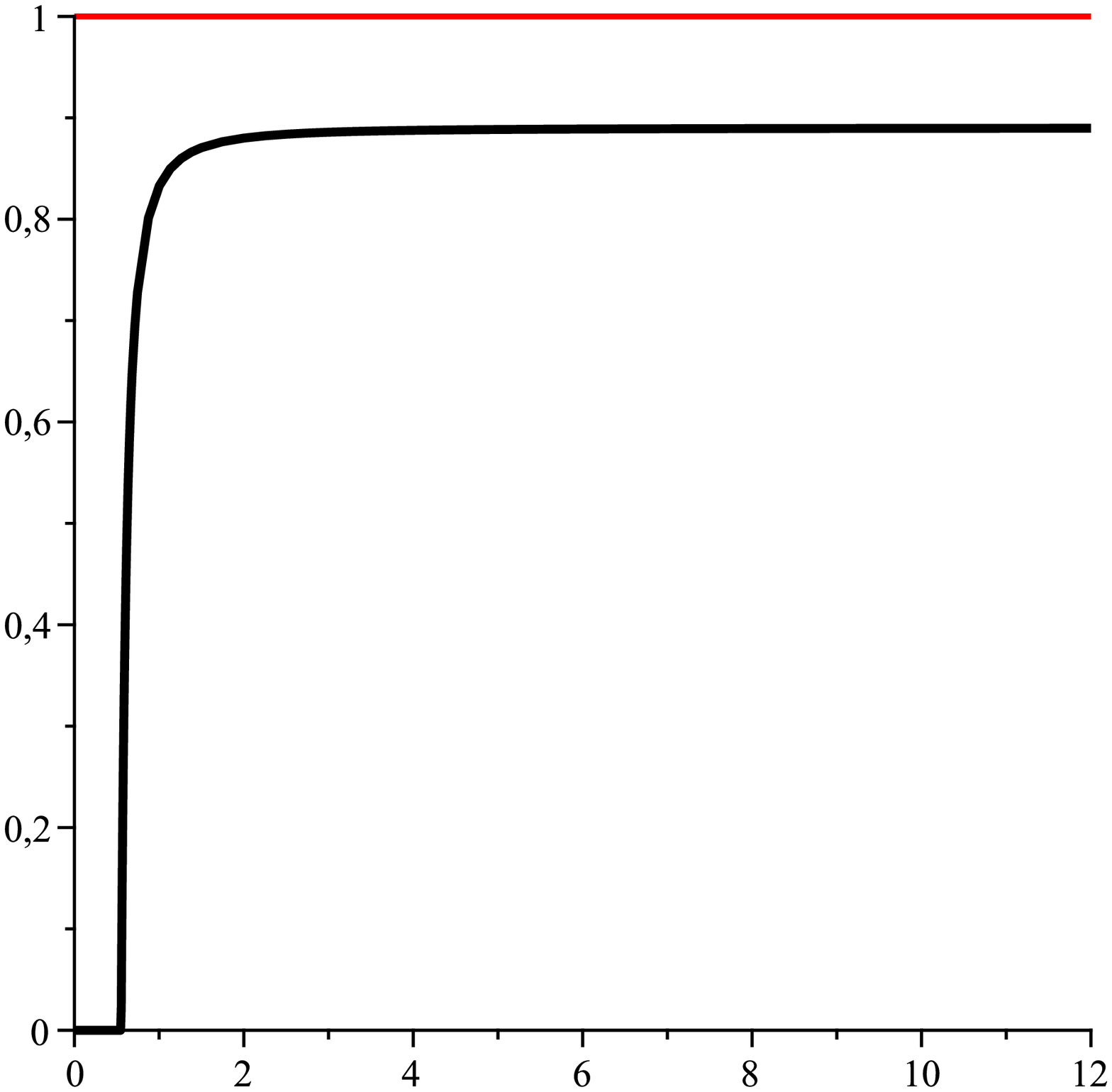} & \includegraphics[scale=0.182]{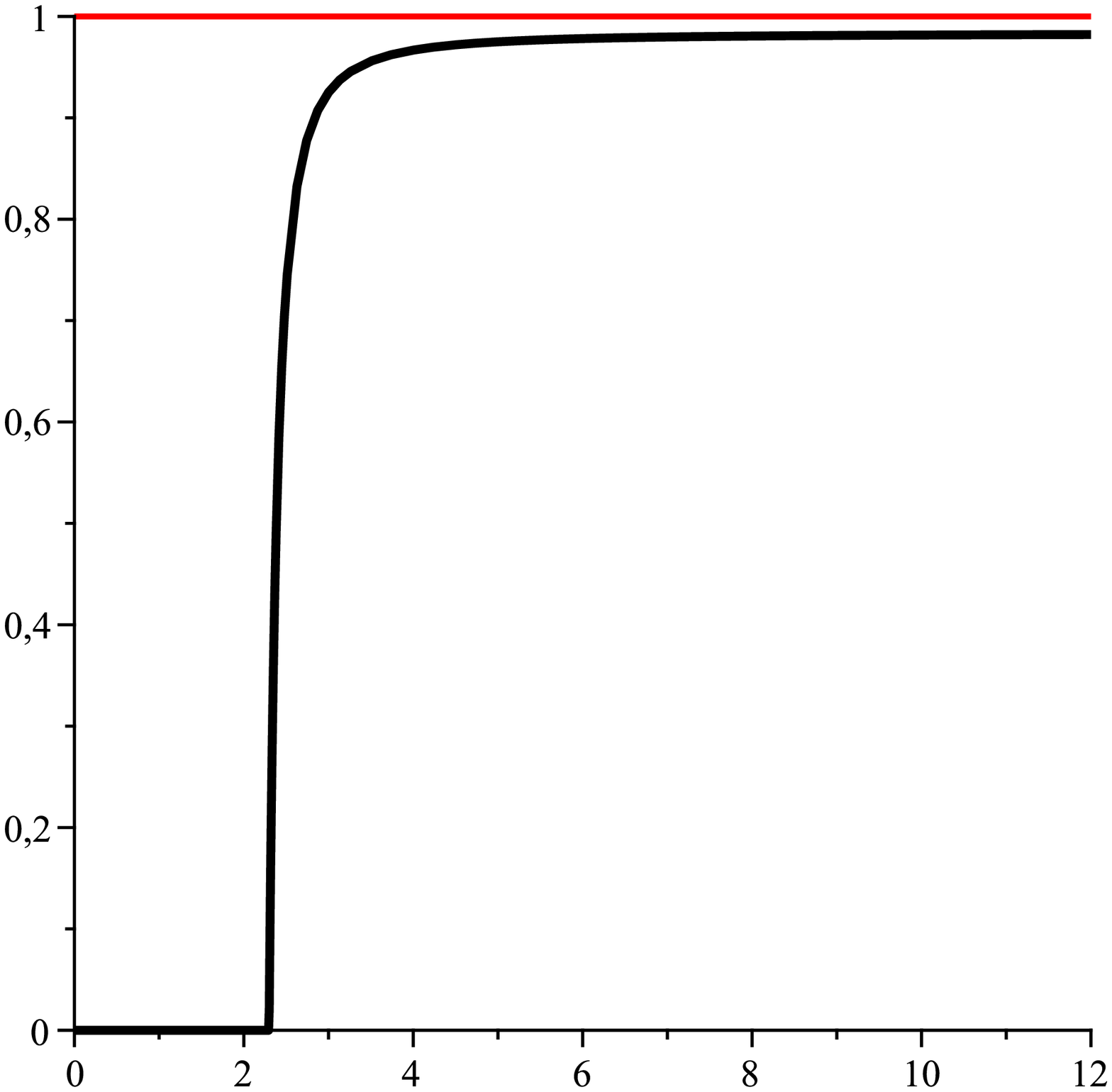}& \includegraphics[scale=0.182]{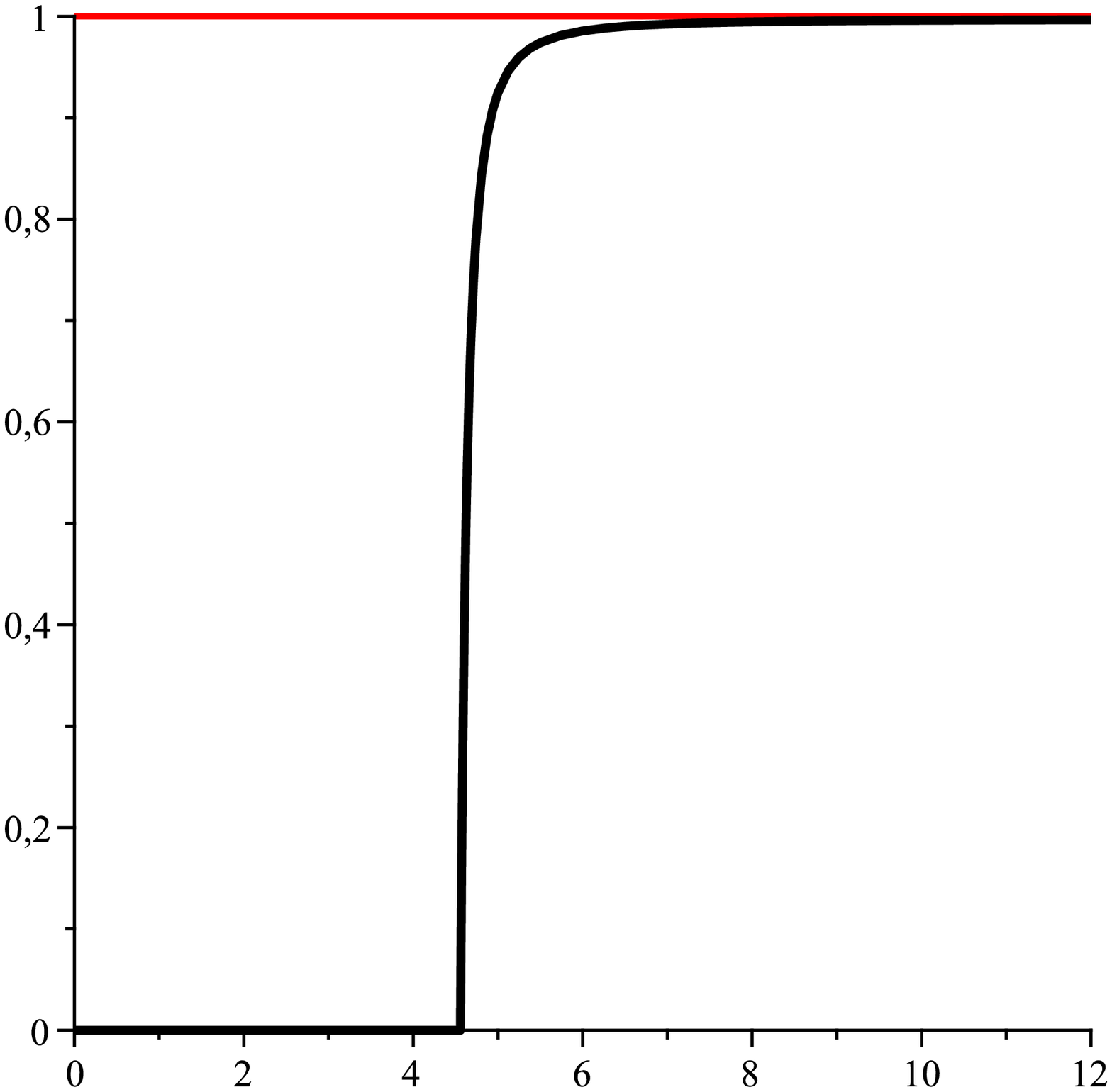} & \includegraphics[scale=0.182]{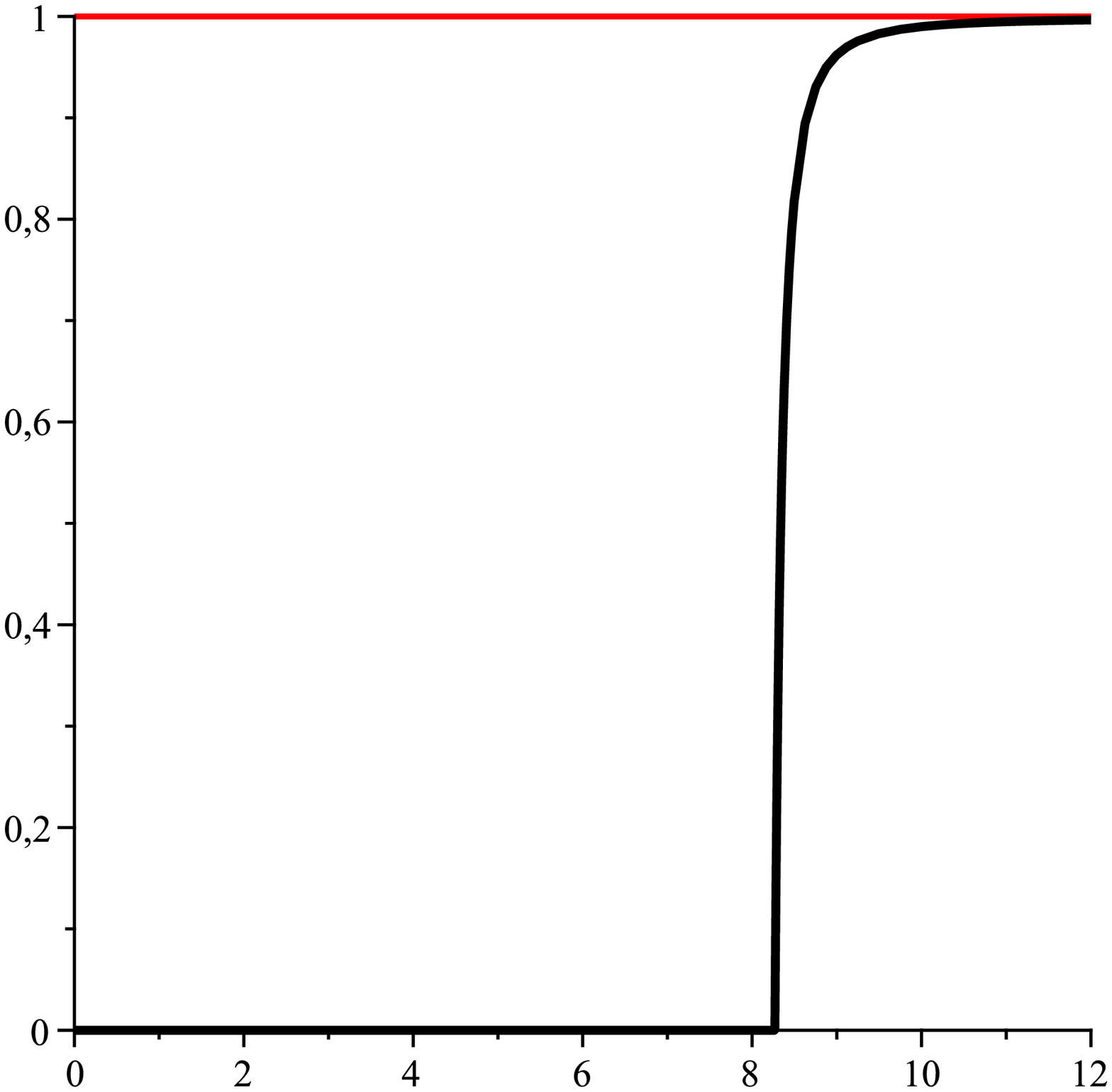} \\
$r=1$ & $r=2$ & $r=3$ & $r=4$  \\
\end{tabular}
\end{center}
Similarly, for a fixed small odd value of $r$, we plot the lower bounds of theorem \ref{th:thm3} in different cases: when there is no restriction on the sign $\varepsilon (Sym^r f)$ of the functional equation, when $\varepsilon (Sym^r f)=+1$ and, finally, when $\varepsilon (Sym^r f)=-1$. 
\newpage
\noindent First, if $r=1$:   
\begin{center}
\begin{tabular}{ccc}
\includegraphics[scale=0.182]{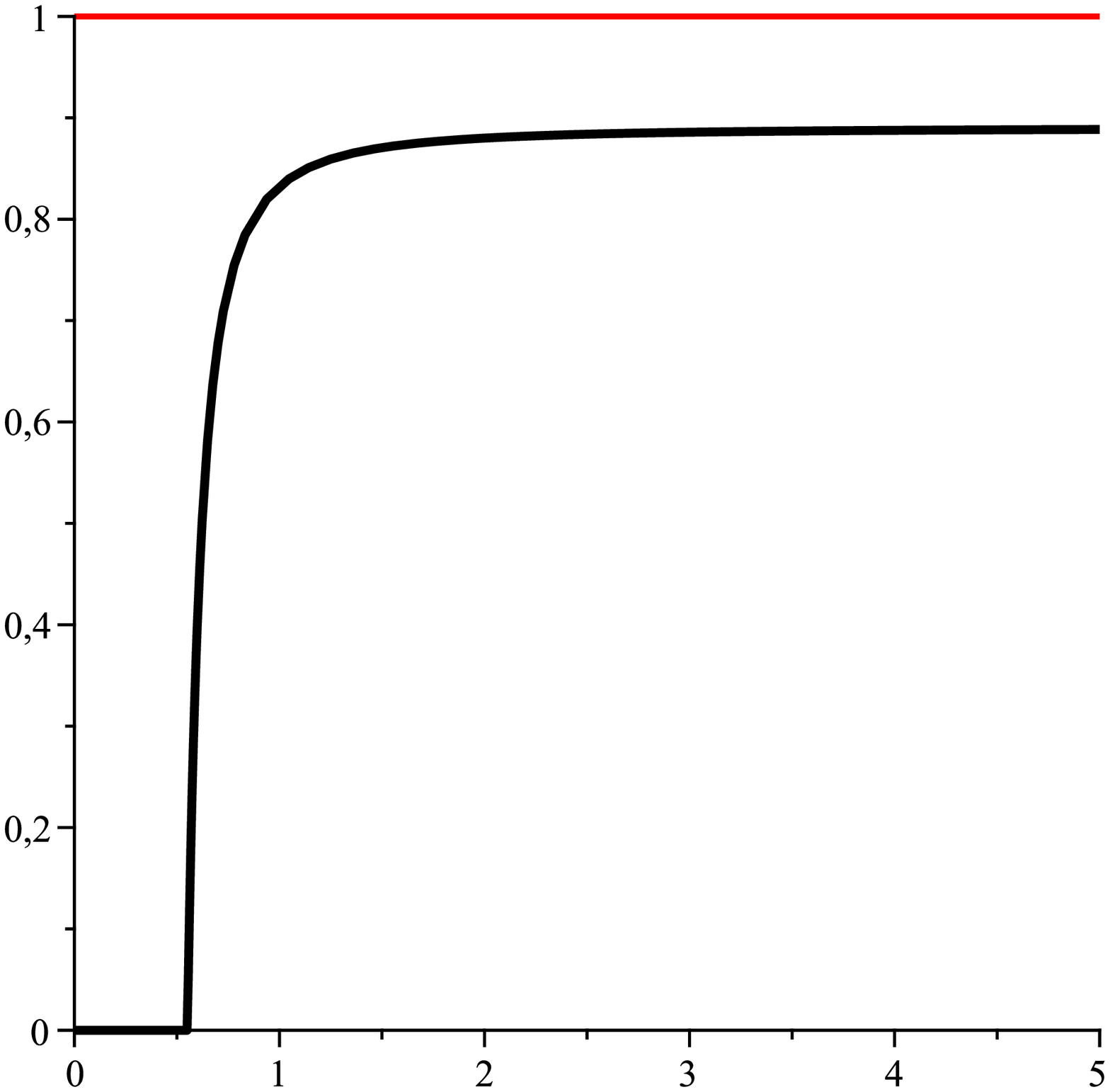} & \includegraphics[scale=0.182]{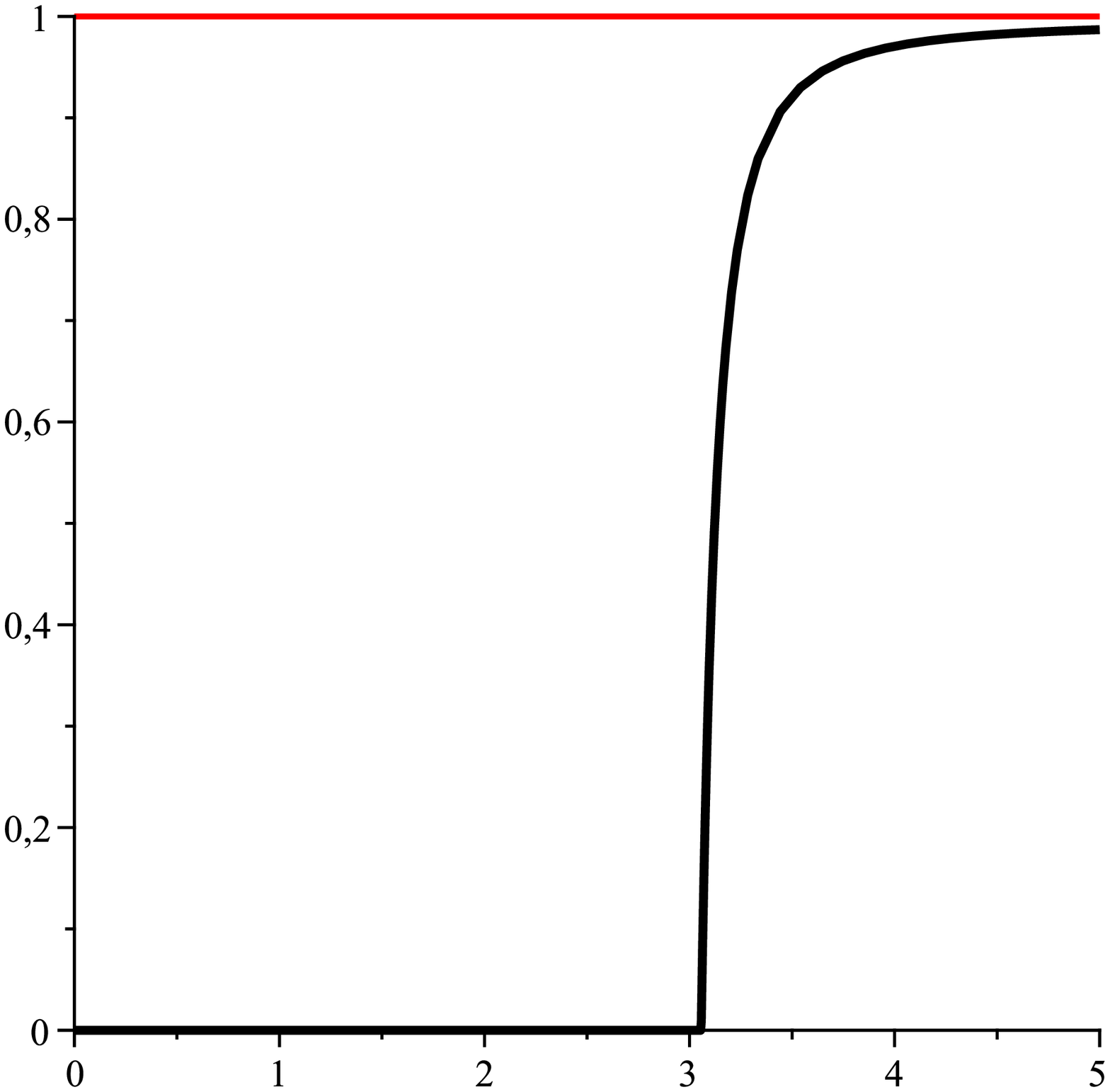}& \includegraphics[scale=0.182]{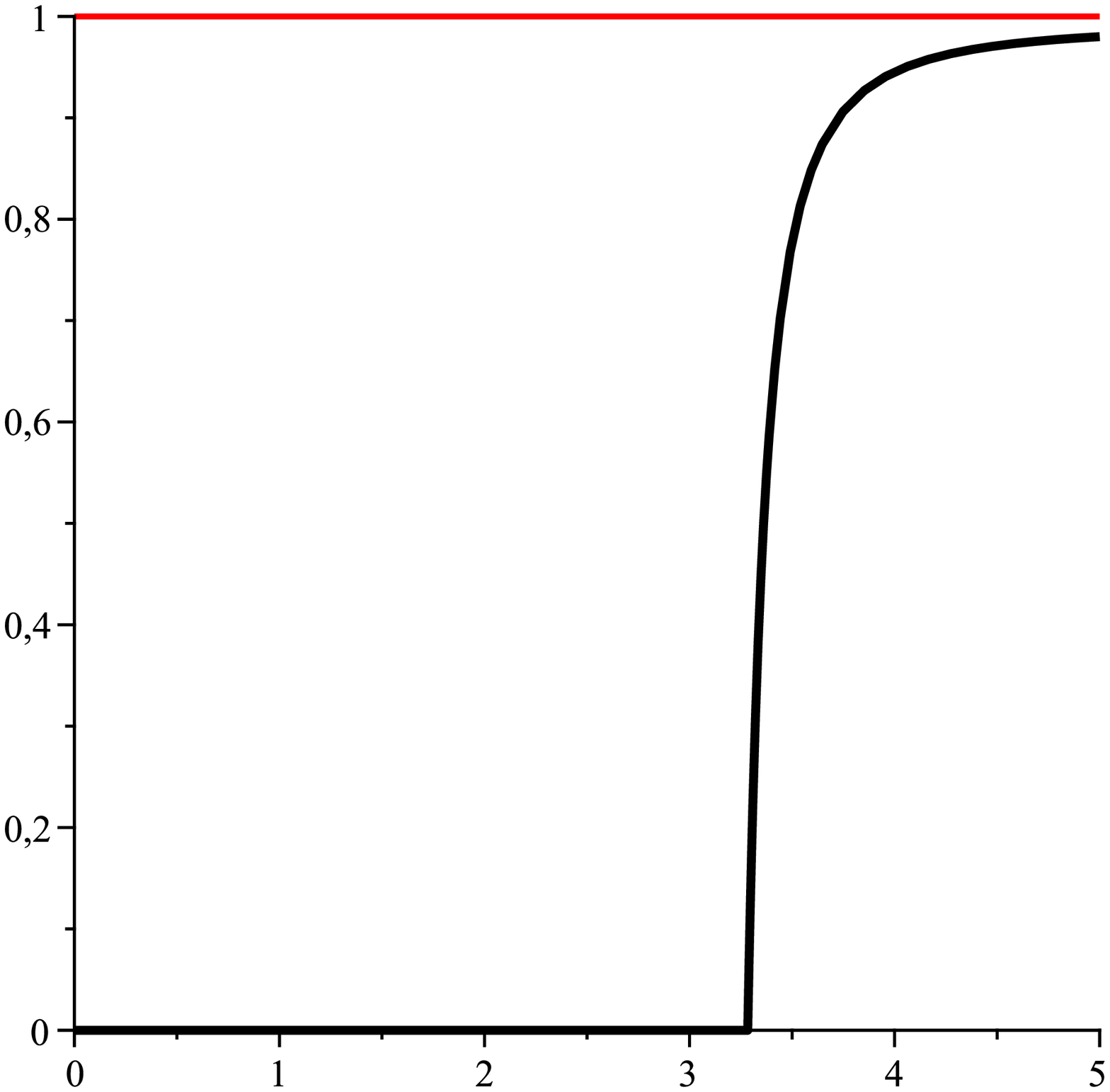} \\
 $\varepsilon (Sym^r f)=\pm 1$ & $\varepsilon (Sym^r f)=+1$ & $\varepsilon (Sym^r f)=-1$  \\
\end{tabular}
\end{center} 
Second, if $r=3$:
\begin{center}
\begin{tabular}{ccc}
\includegraphics[scale=0.182]{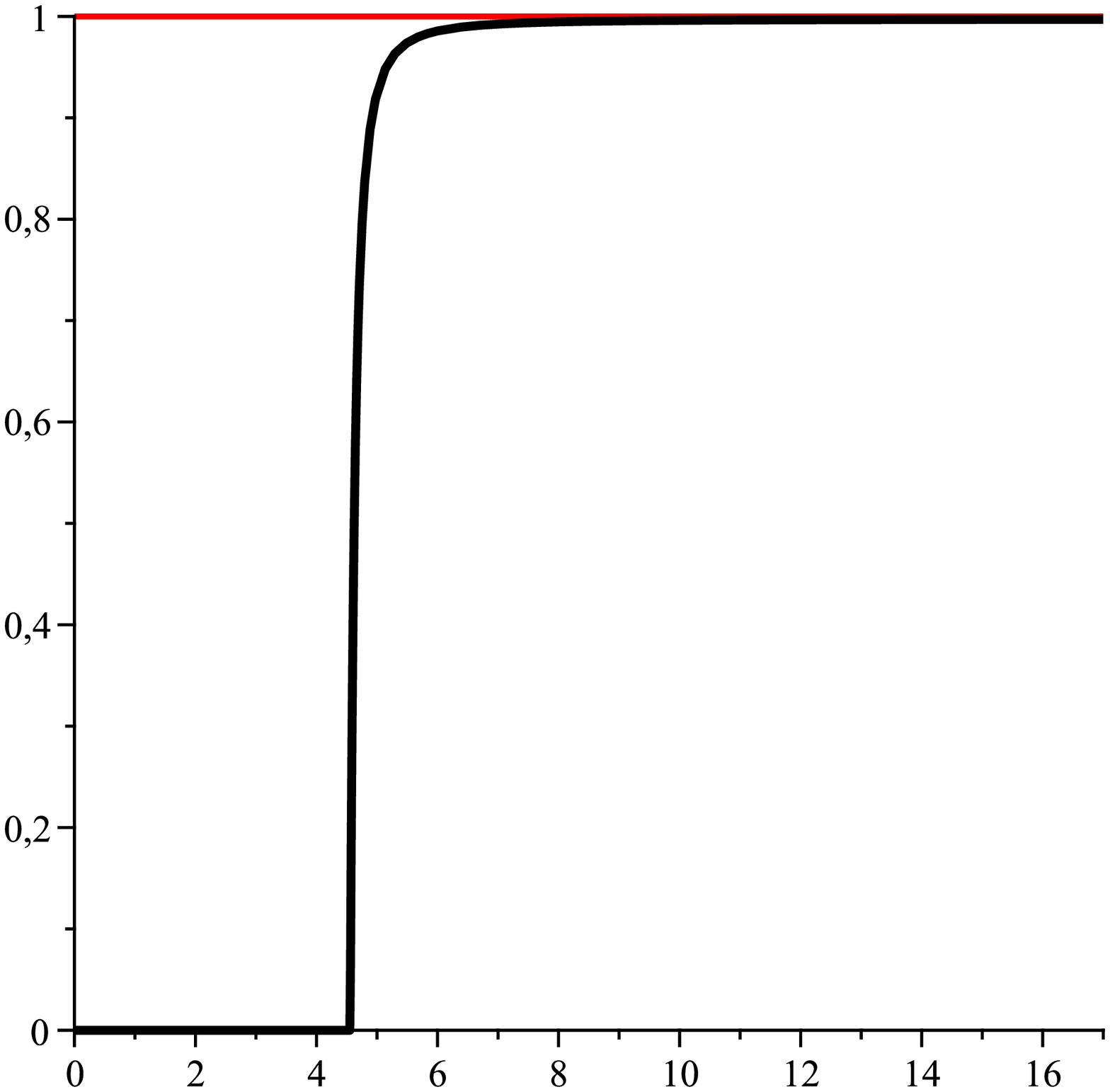} & \includegraphics[scale=0.182]{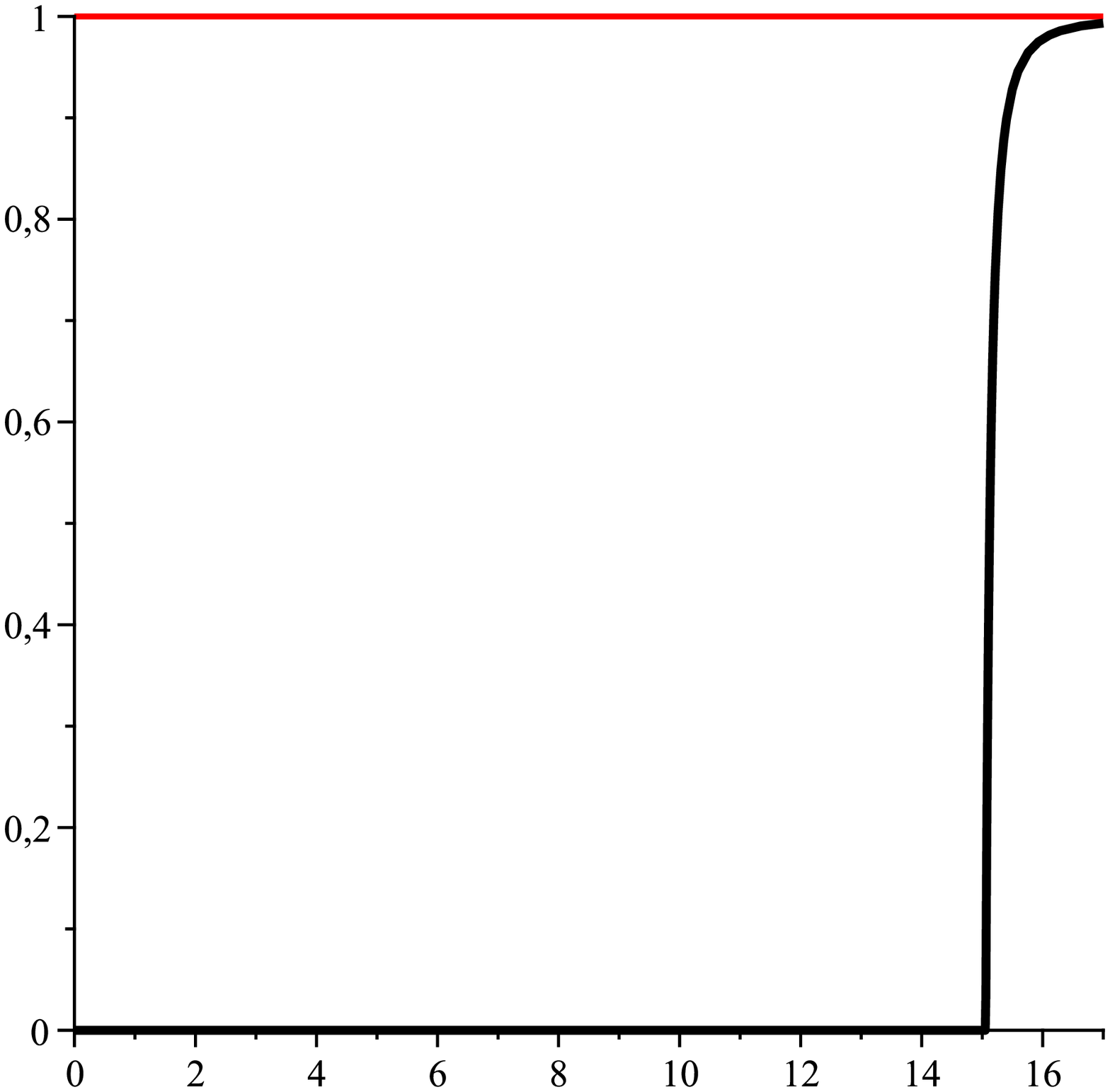}& \includegraphics[scale=0.182]{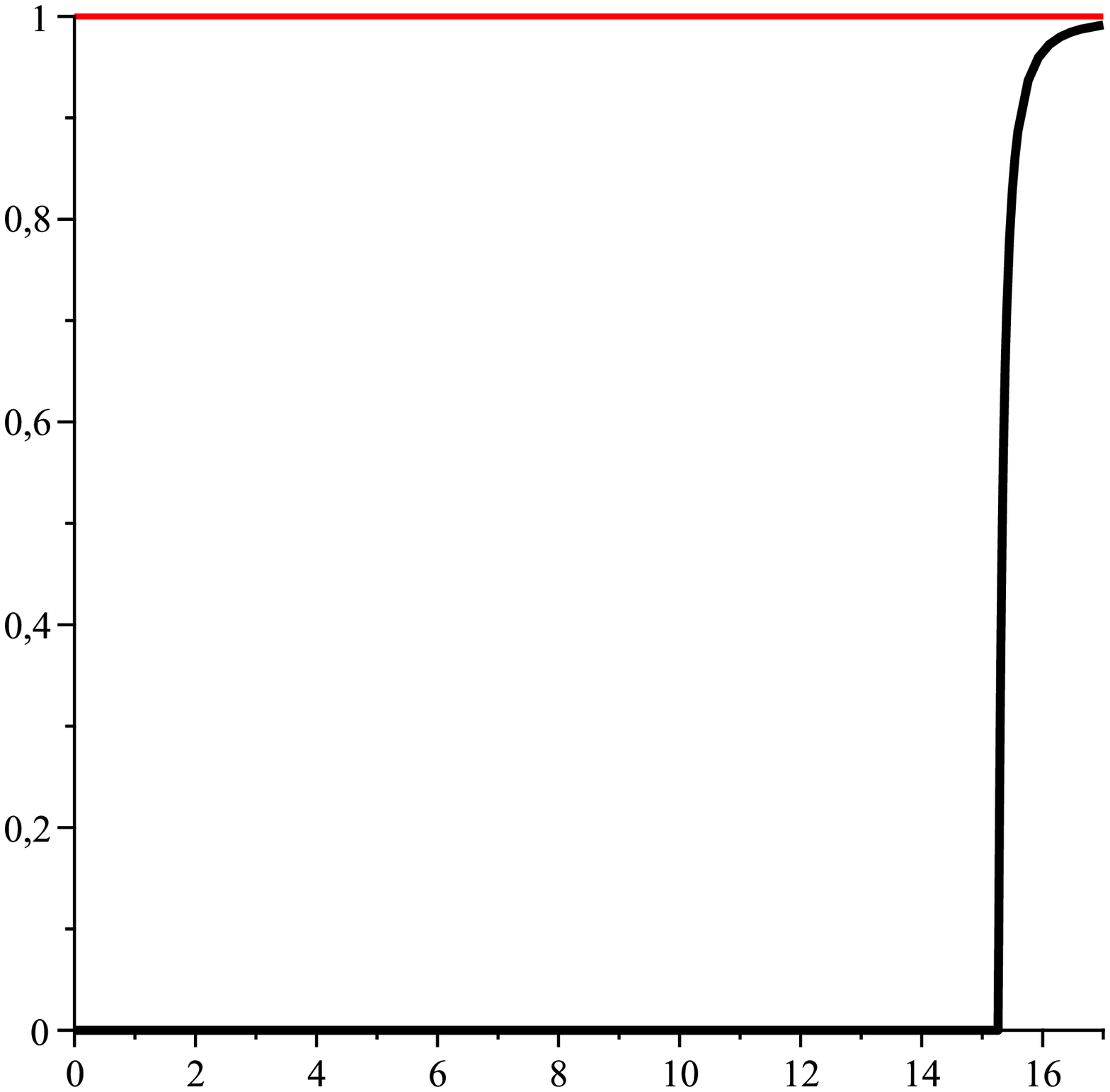} \\
 $\varepsilon (Sym^r f)=\pm 1$ & $\varepsilon (Sym^r f)=+1$ & $\varepsilon (Sym^r f)=-1$  \\
\end{tabular}
\end{center}

  \subsection{Notations}
  
The following notations will be used throughout this paper.
\begin{itemize} 
\item $\lfloor x\rfloor$ denotes the floor of the real number $x$.
\item $\sum^{*}_{n\geq 0}$ means the sum is running over odd non-negative integers.
\item $\mathcal{C}_c^\infty(\R)$  denotes the set of infinitely differentiable functions which are compactly supported.
\item For $1\leq p <+\infty$, $L^p(\R)$ refers to the set of functions $f:\R\rightarrow \R$ such that $\int_\R |f(x)|^p dx<+\infty$. In this case, we put $||f||_p=\left(\int_\R |f(x)|^p dx  \right)^{1/p}$. If $f$ and $g$ are in $L^2(\R)$, let $(f,g)_{L^2}=\int_R f(x)g(x)dx$. 
\item If $\Phi$ is in $L^1(\R)$, $\widehat{\Phi}(u)=\int_{\R}\Phi(x)e^{-2i\pi xu}dx$ is called the Fourier transform of $\Phi$. When it is allowed, we can apply the inverse transform formula $\Phi(x)=\int_{\R}\widehat{\Phi}(u)e^{2i \pi xu}du $.
\item If $f$ and $g$ are in $L^2(\R)$, $f*g(u)=\int_\R f(t)g(u-t)dt$ is the convolution product of $f$ and $g$.
\item $\mathcal{S}_\nu (\R)$ denotes the set of even Schwartz functions whose Fourier transforms are compactly supported in $[-R,R]$ with $0<R<\nu$.
\end{itemize}

\paragraph*{Acknowledgements} I would like to thank Frederic Bayart for his availability and all their advices on Sobolev spaces.

\subsection{Statistics for low-lying zeros}  
\paragraph*{}
Let $\mathcal{F}$ be a natural family of $L$-functions\footnote{Our definition of $L$-function is the one of \cite{IK} chapter 5.} all of whose satisfies the Riemann hypothesis. Consider a $L$-function $L(f,s)$ in $\mathcal{F}$ with analytic conductor $c_f$. Let $\mathcal{F}(Q)= \left\{L(f,s) \in  \mathcal{F} \hspace{1mm} , \hspace{1mm} c_f = Q \right\}$. Each zero $\rho_f$ of $L(f,s)$ which is on the critical line $\Re(s)=1/2$  can be written  $\rho_f=\frac{1}{2}+i\gamma_f$, and we denote $\tilde{\rho}_f=\frac{1}{2}+i\tilde{\gamma}_f$ with $\tilde{\gamma}_f=\gamma_f\frac{\ln c_f}{2\pi}$ the normalized zero.
\paragraph{}
In order to study the distribution of low-lying zeros of $L(f,s)$, for any test function $\Phi$ in $\mathcal{S}_\nu (\R)$, we define the low-zeros sum 
$$ D[\Phi](f) = \sum_{\tilde{\gamma_f}} \Phi ( \tilde{\gamma}_f ) $$
where the sum is running over the imaginary parts of normalised zeros counted with multiplicity. $\mathcal{F}(Q)$ can be seen as a measurable space where measurable sets are all its subsets and which is equipped with the counting probability measure $\mu_{\mathcal{F}(Q)}$. $D[\Phi]$ is a measurable function on $\mathcal{F}(Q)$. However, we are unable to determine the asymptotic behaviour of $D[\Phi](f)$ for a single $L$-function. As a consequence, we must take into account a family of $L$-functions in order to obtain a significant result.  That is why, we define the one-level density as the expectation of $D[\Phi]$:
\begin{small}
$$   \E_{\mathcal{F}(Q)} (D[\Phi]) = \frac{1}{\left|  \mathcal{F}(Q)  \right| }\sum_{L(f,s)\in \mathcal{F}(Q)} D[\Phi](f)    $$ 
\end{small}We also define the variance of the one-level density:
\begin{small}
$$ \V_{\mathcal{F}(Q)} (D[\Phi]) = \E([D[\Phi]-\E_{\mathcal{F}(Q)} (D[\Phi])  ]^2)= \frac{1}{\left|  \mathcal{F}(Q)  \right| }\sum_{L(f,s)\in \mathcal{F}(Q)} \left[  D[\Phi](f)-\frac{1}{\left|  \mathcal{F}(Q)  \right| }\sum_{L(f,s)\in \mathcal{F}(Q)} D[\Phi](f)\right]^2$$
\end{small}Our purpose is to find the asymptotic behaviour of $\E_{\mathcal{F}(Q)} (D[\Phi])$ and $\V_{\mathcal{F}(Q)} (D[\Phi])$ when $Q$ goes to infinity. The density conjecture predicts
$$ \lim_{Q \rightarrow +\infty} \frac{1}{\left|  \mathcal{F}(Q)  \right| }\sum_{L(f,s) \in \mathcal{F}(Q)} D[\Phi](f) = \int_{\R}\Phi(t)W[\mathcal{F}](t) dt$$ 
where $ W[\mathcal{F}]$ is a density function characterised by $\mathcal{F}$. In order to estimate the one-level density, we convert sums over zeros to sums over primes. Unfortunately, we are able to evaluate these sums over primes only if the support of the test function is small. That's why, density theorems are proved only for test functions in $\mathcal{S}_{\nu_{\max}(\mathcal{F})} (\R)$ with $\nu_{\max}(\mathcal{F})$ fixed. Currently, the maximal value for $\nu_{\max}(\mathcal{F})$ is $2$ whereas the density conjecture does not predict any restriction on the support of test functions. 
\subparagraph*{}
Sometimes, we will prefer using the harmonic measure for technical conveniences, rather than the Dirac one. For instance, in the case of symmetric power $L$-functions. Precisely, if $A$ is a subset of $\mathcal{H}_r(q)$, let
$$\mu_{{\mathcal{H}_r}(q)}^h(A)=\sum_{L(Sym^r f,s)\in A}w_q(f)$$
and 
\begin{small}
$$\E^h_{{\mathcal{H}_r}(q)} (D[\Phi]) = \sum_{f\in H_k^*(q)} \omega_q(f) D[\Phi](Sym^r f) \mbox{ where }\omega_q(f)= \frac{\Gamma(k-1)}{(4\pi)^{k-1}\langle f,f \rangle_q}  $$
$$\mbox{and }\hspace{2mm}\V^h_{{\mathcal{H}_r}(q)} (D[\Phi])= \sum_{f\in H_k^*(q)} \omega_q(f) \left[ D[\Phi](Sym^r f) - \E^h_{{\mathcal{H}_r}(q)} (D[\Phi])\right]^2.$$
\end{small}Then, if $r$ is odd and if $A$ is a subset of $\mathcal{H}_r^\pm(q)$, let
$$\mu_{{\mathcal{H}_r^\pm}(q)}^h(A)=2\sum_{L(Sym^r f,s)\in A}w_q(f)$$
and 
\begin{small}
$$ \E^h_{{\mathcal{H}_r^\pm}(q)} (D[\Phi]) = 2 
\sum_{\begin{array}{c}  \scriptstyle{ f \in H_{k}^{*}(q)}  \\ \scriptstyle{\varepsilon (Sym^r f)=\pm 1} \end{array}} \omega_q(f) D[\Phi](Sym^r f)$$
$$\mbox{and }\hspace{2mm}\V^h_{{\mathcal{H}_r^\pm}(q)} (D[\Phi])= 2\sum_{\begin{array}{c}  \scriptstyle{ f \in H_{k}^{*}(q)}  \\ \scriptstyle{\varepsilon (Sym^r f)=\pm 1} \end{array}} \omega_q(f) \left[ D[\Phi](Sym^r f) - \E^h_{{\mathcal{H}_r^\pm}(q)} (D[\Phi])\right]^2. $$
\end{small}Actually,  harmonic measures are  asymptotic probability measures since we only have 
\begin{eqnarray}   \label{eqn:mesureharmonique}   \lim_{\underset{ q \mbox{ \begin{tiny}prime\end{tiny}}}{q \rightarrow +\infty}}   \sum_{f\in H_{k}^{*}(q)} \omega_q(f)=1        \mbox{ and }\lim_{\underset{ q \mbox{ \begin{tiny}prime\end{tiny}}}{q \rightarrow +\infty}} \sum_{ \begin{array}{c}  \scriptstyle{ f \in H_{k}^{*}(q)}  \\ \scriptstyle{\varepsilon (Sym^r f)=\pm 1} \end{array}       } \omega_q(f) = \frac{1}{2} . \end{eqnarray}
The first relation comes from Petersson trace formula (\cite{RR}, proposition 2.2) and the second one is subject to an assumption (Hypothesis $Nice(r,f)$ of \cite{RR})  we are assuming in order to get density theorems for these families.

\subsection{What is the smallest zero of a $L$-function ?}
If $L(f,s)$ is a self-dual $L$-function, the sign of the functional equation is equal to $\pm1$. Moreover, due to the following observation, we need to define $\tilde{\gamma}_{f,1}$ which appears in theorems \ref{th:thm1} and \ref{th:thm3} and to explain the consequences on density theorems.
\begin{lemma}[\cite{IK}, proposition 5.1] \label{lem:lemmezerotrivial}
Let $L(f,s)$ be a self-dual $L$-function with $\varepsilon(f)=-1$. Then $L\left(f,\frac{1}{2}\right)=0$.
\end{lemma}
If $\mathcal{F}$ is a natural family, all of whose $L$-functions $L(f,s)$ are self-duals and satisfy $\varepsilon(f)=-1$. Due to the previous lemma, we can denote non-trivial  zeros of $L(f,s)$ by
 $$   \{\rho_{f,0} \} \cup \{   \rho_{f,i}, i \in \Z^*   \}    $$
 where $\rho_{f,0}=1/2$ and $\rho_{f,i}=1-\rho_{f,-i}$ if $i\neq 0$. Moreover, we have:
 $$  ...\leq Im(\rho_{f,-1}) \leq Im(\rho_{f,0}) = 0  \leq   Im(\rho_{f,1}) \leq Im(\rho_{f,2}) \leq ... $$ 
In the other cases, we use the same notations without $\rho_{f,0}$.
\subparagraph*{}
We are going to expose its consequences on statistics of low-lying zeros. Let $\mathcal{F}$ a natural family of $L$-functions, we define
$$D^*[\Phi](f) = \sum_{i\in \Z^*} \Phi ( \tilde{\gamma}_{f,i}) . $$
The density function $W^*[\mathcal{F}]$  is defined to satisfy
\begin{eqnarray} \label{eqn:defdensitepointe}
 \lim_{Q \rightarrow +\infty} \E_{\mathcal{F}(Q)}( D^*[\Phi]) = \int_{\R}\Phi(t)W^*[\mathcal{F}](t) dt .
\end{eqnarray}
If  all $L$-functions $L(f,s)$ in $\mathcal{F}$ are self-duals and satisfy $\varepsilon(f)=-1$, we have
\begin{eqnarray} \label{eqn:relationentreesperance} \E_{\mathcal{F}(Q)} (D[\Phi])=\Phi(0)+\E_{\mathcal{F}(Q)} (D^*[\Phi]) \hspace{2mm} \mbox{ and }\hspace{2mm} \V_{\mathcal{F}(Q)} (D[\Phi])=\V_{\mathcal{F}(Q)} (D^*[\Phi]). \end{eqnarray}
This phenomenon occurs for the family $\mathcal{H}^-$. It has been shown in \cite{ILS} (equation (1.18)) that the symmetry group associated to  $\mathcal{H}^-$ is $SO^-$. In other words, we have
$$\lim_{\underset{ N \mbox{ \begin{tiny}squarefree\end{tiny}}}{ N\rightarrow +\infty}} \E_{\mathcal{H}^-(N)} (D[\Phi])=\int_{\R}\Phi(t)W[SO^-](t) dt\hspace{2mm} \mbox{ ie }\hspace{2mm} W[\mathcal{H}^-]=W[SO^-]. $$
As a result, we deduce
\begin{eqnarray}
\lim_{\underset{ N \mbox{ \begin{tiny}squarefree\end{tiny}}}{ N\rightarrow +\infty}} \E_{\mathcal{H}^-(N)} (D^*[\Phi])=\int_{\R}\Phi(t)W[Sp](t) dt\hspace{2mm} \mbox{ ie }\hspace{2mm} W^*[\mathcal{H}^-]=W[Sp] .
\end{eqnarray}
Similarly, the family $\mathcal{H}_r^-$ ($r$ odd) has also $SO^-$ as symmetry group (\cite{RR}, theorem $A$). Nevertheless, for the harmonic measure,  relation $(\ref{eqn:relationentreesperance})$ becomes  
\begin{small} $$\E_{\mathcal{H}_r^-(q)}^h (D[\Phi])=2\Phi(0)\mu_{\mathcal{H}_r(q)}^h\left(\mathcal{H}_r^-(q) \right)+\E_{\mathcal{H}_r^-(q)}^h (D^*[\Phi]) $$ 
\end{small}\mbox{and }
\begin{small}
$$
 \V_{\mathcal{H}_r^-(q)}^h (D[\Phi])= \V_{\mathcal{H}_r^-(q)}^h (D^*[\Phi])+2\Phi (0)  \left[\Phi (0)\mu_{\mathcal{H}_r(q)}^h\left(\mathcal{H}_r^-(q) \right)+ \E_{\mathcal{H}_r^-(q)}^h (D^*[\Phi]) \right]\left[1-2\mu_{\mathcal{H}_r(q)}^h\left(\mathcal{H}_r^-(q) \right) \right]^2 . $$ \end{small}Thanks to relation $(\ref{eqn:mesureharmonique})$, we have
\begin{small}
$$\lim_{\underset{ q \mbox{ \begin{tiny}prime\end{tiny}}}{q \rightarrow +\infty}} \E_{\mathcal{H}_r^-(q)}^h (D[\Phi])=\Phi(0)+\lim_{\underset{ q \mbox{ \begin{tiny}prime\end{tiny}}}{q \rightarrow +\infty}}\E_{\mathcal{H}_r^-(q)}^h (D^*[\Phi]) \hspace{1mm} \mbox{ and }\hspace{1mm} \lim_{\underset{ q \mbox{ \begin{tiny}prime\end{tiny}}}{q \rightarrow +\infty}}\V_{\mathcal{H}_r^-(q)}^h (D[\Phi])=\lim_{\underset{ q \mbox{ \begin{tiny}prime\end{tiny}}}{q \rightarrow +\infty}}\V_{\mathcal{H}_r^-(q)}^h (D^*[\Phi]).$$
\end{small}As a result, with theorems $A$ and $D$ of \cite{RR}, we get
\begin{small} 
\begin{eqnarray}
\lim_{\underset{ q \mbox{ \begin{tiny}prime\end{tiny}}}{ q\rightarrow +\infty}} \E_{\mathcal{H}_r^-(q)}^h (D^*[\Phi])=\int_{\R}\Phi(t)W[Sp](t) dt \hspace{2mm} \mbox{ and } \hspace{2mm} \lim_{\underset{ q \mbox{ \begin{tiny}prime\end{tiny}}}{ q\rightarrow +\infty}} \V_{\mathcal{H}_r^-(q)}^h (D^*[\Phi])=2\int_\R |u|\widehat{\Phi}(u)^2du. 
\end{eqnarray}
\end{small}

\section{Proportion of $L$-functions with a small smallest zero} 
In this part, we prove theorem \ref{th:thm3}. The starting point is the following proposition. We do not give a proof of this result since it is essentially the same than the proof theorem 8.3 from \cite{HR}. 
\begin{proposition} \label{pro:proportioncasgeneral}Let $g$ be in $\mathcal{S}_R(\R)$  and $\Phi(x)=(x^2-\beta^2)g^2(x)$.  Let
 $$B(g)=\sqrt{\frac{\int_{\R}  x^2g^2(x)W^*[\mathcal{F}](x)dx}{\int_{\R}g^2(x)W^*[\mathcal{F}](x) dx }}.$$  We assume
\begin{small} 
\begin{eqnarray*}
\lim_{\underset{ q \mbox{ \begin{tiny}prime\end{tiny}}}{ q\rightarrow +\infty}} \E_{\mathcal{F}(q)}^h (D^*[\Phi])= \int_{\R}\Phi(t)W^*[\mathcal{F}](t) dt  \hspace{2mm} \mbox{ and } \hspace{2mm} \lim_{\underset{ q \mbox{ \begin{tiny}prime\end{tiny}}}{ q\rightarrow +\infty}} \V_{\mathcal{F}(q)}^h (D^*[\Phi])=\V_\mathcal{F}(\Phi) .
\end{eqnarray*}
\end{small}Then, if  $\beta > B(g)$, we have
\begin{eqnarray} \label{eqn:proportioncasgeneral}
\liminf_{\underset{ q \mbox{ \begin{tiny}prime\end{tiny}}}{ q\rightarrow +\infty}}\mu_{\mathcal{F}(q)}^h \left(\left\{ L(f,s)\in \mathcal{F}; \, \tilde{\gamma}_{f,1} < \beta \right\} \right) \geq 1 - \frac{\V_\mathcal{F}(\Phi) }{\left[ \int_{\R}\Phi(t)W^*[\mathcal{F}](t) dt   \right]^2} =: Borne_\mathcal{F}(\beta^2).  
\end{eqnarray}
\end{proposition} 
\begin{remark}
The same result holds for all natural family of $L$-functions with the counting probability measure instead of the harmonic measure.
\end{remark} 
\paragraph*{}
This proposition gives a result only if the right member term is positive. We denote $\beta_{\min} ( g)$  the smallest value of $\beta > B(g)$ such that this property is satisfied. Thanks to $\beta_{\min} ( g)>B(g)$, we may detect a zero. We have:
$$  \liminf_{\underset{ q \mbox{ \begin{tiny}prime\end{tiny}}}{ q\rightarrow +\infty}}  \min_{L(f,s)\in \mathcal{F}(q)} \tilde{\gamma}_{f,1} \leq \beta_{\min} ( g) $$   
Actually, we prove a better upper bound in theorem \ref{th:thm1}.
\paragraph*{}
Ricotta and  Royer proved in \cite{RR} (theorems A, B and D), that if $\Phi$ is in $\mathcal{S}_{\rho_{\max}(\mathcal{F})}(\R)$, then 
\begin{eqnarray*}
\lim_{\underset{ q \mbox{ \begin{tiny}prime\end{tiny}}}{ q\rightarrow +\infty}} \E_{\mathcal{F}(q)}^h (D[\Phi])= \int_{\R}\Phi(t)W[\mathcal{F}](t) dt  \hspace{2mm} \mbox{ and } \hspace{2mm} \lim_{\underset{ q \mbox{ \begin{tiny}prime\end{tiny}}}{ q\rightarrow +\infty}} \V_{\mathcal{F}(q)}^h (D[\Phi])=2 \int_\R |u|\widehat{\Phi}(u)^2du
\end{eqnarray*}
where $\rho_{\max}(\mathcal{F})$ and $W[\mathcal{F}]$ are given in table \ref{tab:tableadenumax}. Since $\rho_{\max}(\mathcal{F})<1$, we have
$$ \int_{\R}\Phi(t)W^*[\mathcal{F}](t) dt = \widehat{\Phi}(0)+\frac{\sigma_{\mathcal{F}}}{2}\Phi(0) $$
with:
$$
\begin{array}{|c|c|c|c|c|}
\hline
\mathcal{F}     \rule[-4mm]{0cm}{10mm} & \mathcal{H}_r \hspace{2mm}(r\mbox{ even})  & \mathcal{H}_r \hspace{2mm}(r\mbox{ odd}) & \mathcal{H}_r^+  \hspace{2mm}(r\mbox{ odd}) & \mathcal{H}_r^-  \hspace{2mm}(r\mbox{ odd}) \\
\hline
\sigma_{\mathcal{F}}\rule[-4mm]{0cm}{10mm} & (-1)^{r+1}=-1 & (-1)^{r+1}=1 & 1 & -1 \\
\hline
\end{array}
$$
In order to obtain explicit lower bound in the previous proposition, we specialise relation (\ref{eqn:proportioncasgeneral}) in a fixed test function. 
\begin{lemma} \label{lem:quelquescalculs} If $0< R<1/2$,  $\widehat{g}_0(u)=\cos\left(\frac{\pi u}{2R}  \right)\un_{[-R,R]}(u)$ and $\Phi(x)=(x^2-\beta^2)g_0^2(x)$, then
\begin{small}
\begin{eqnarray*}
Borne_\mathcal{F}(\beta^2) = 1 - \frac{2\pi^2R^2}{3}\frac{256(3+\pi^2)R^4\beta^4+32(9-\pi^2)R^2\beta^2+\pi^2+3}{\left( 128\sigma_{\mathcal{F}} R^3\beta^2+16\pi^2R^2\beta^2-\pi^2\right)^2}\hspace{1mm}\mbox{ and } \hspace{1mm}  B(g_0)  =  \frac{1}{4R\sqrt{1+\sigma_{\mathcal{F}}\frac{8R}{\pi^2} }}.
\end{eqnarray*}
\end{small}
\end{lemma}  
Proof: We have
\begin{small}
\begin{eqnarray*}
 & \widehat{\Phi}(u) = \left[\frac{2R-|u|}{2}\left( \frac{1}{16R^2}-\beta^2 \right)\cos \left(  \frac{\pi u}{2R}\right)-\frac{1}{\pi}\left( \frac{1}{16R}+\beta^2 R \right)\sin \left(\frac{\pi |u|}{2R}  \right)  \right].\un_{[-2R,2R]}(u)&  \\
 & \widehat{\Phi}(0)  =  \frac{1}{16R}-\beta^2R  \hspace{1cm} \mbox{ and }\hspace{1cm} \Phi(0) =  -\frac{16\beta^2R^2}{\pi^2} . &
 \end{eqnarray*}
\end{small}Thanks to several integration by parts, we obtain
\begin{small}
\begin{eqnarray*}
\int_{\R}|u|\widehat{\Phi}(u)^2du = \frac{768R^4\beta^4+3+288\beta^2R^2+\pi^2-32\beta^2R^2\pi^2+256R^4\beta^4\pi^2}{768\pi^2}.
\end{eqnarray*}
\end{small}
\begin{flushright}
$\Box$
\end{flushright}
\begin{lemma} \label{lem:variationdelaborne} Let  $Borne_\mathcal{F}^+(\beta^2)=\max\left\{ 0;Borne_\mathcal{F}(\beta^2) \right\}$. The following table sums up variations of  $\beta \longmapsto Borne_\mathcal{F}^+(\beta^2)$ where 
\begin{small}
$$ \beta_{\min}(g_0)= \frac{1}{4R} \sqrt{\frac{\pi^2 (3\pi^2+24\sigma_\mathcal{F} R-2(\pi^2-9)R^2) +4\pi R\sqrt{9\pi^4+72\sigma_\mathcal{F} R\pi^2-6(\pi^4-7\pi^2-12)R^2 } }{3\pi^4+48\sigma_\mathcal{F} R\pi^2+(192-6\pi^2-2\pi^4)R^2} }  .$$
\end{small}
\begin{center}
\begin{footnotesize}
\begin{variations}
\beta & B(g_0) &      & \beta_{\min}(g_0) &       & \pI \\
\filet
\m{Borne_\mathcal{F}^+(\beta^2)} & 0 & \longrightarrow & 0 & \c & \h{1-\frac{2\pi^2R^2}{3}\frac{\pi^2+3}{(\pi^2+8R\sigma_\mathcal{F})^2}} \\
\end{variations}
\end{footnotesize}
\end{center}
\end{lemma}
Proof: First, we study the sign of $Borne_\mathcal{F}(X)$. It is the same as the quadratic polynomial's
\begin{small}
\begin{eqnarray*}
Y(X)=3\pi^4-6\pi^2R^2-2\pi^4R^2+X\left( -96R^2\pi^4-768\sigma_\mathcal{F}\pi^2R^3-576R^4\pi^2+64R^4\pi^4  \right)  \\   \hspace{40mm} + X^2\left(768R^4\pi^4+12288\sigma_\mathcal{F} R^5\pi^2+49152R^6-1536R^6\pi^2-512R^6\pi^4  \right).
\end{eqnarray*}
\end{small}Let $L(R)$ be its leading coefficient and $\Delta(R)$ its discriminant. We get
\begin{small}
\begin{eqnarray*}
 \Delta(R)   & = & \underbrace{-98304(\pi^4-7\pi^2-12)\pi^2}_{<0}R^6 (R-R_1)(R-R_2)
\end{eqnarray*}
\end{small}with
\begin{small} \begin{eqnarray*} 
 R_1 & = & \min \left\{ \pi^2 \frac{12\sigma_\mathcal{F} +\sqrt{6(\pi^2-3)(\pi^2-4)}}{2(\pi^4-7\pi^2-12)};\pi^2 \frac{12\sigma_\mathcal{F} -\sqrt{6(\pi^2-3)(\pi^2-4)}}{2(\pi^4-7\pi^2-12)} \right\} \\
 R_2 & = & \max \left\{ \pi^2 \frac{12\sigma_\mathcal{F} +\sqrt{6(\pi^2-3)(\pi^2-4)}}{2(\pi^4-7\pi^2-12)};\pi^2 \frac{12\sigma_\mathcal{F} -\sqrt{6(\pi^2-3)(\pi^2-4)}}{2(\pi^4-7\pi^2-12)} \right\}.
\end{eqnarray*}\end{small}
\begin{quote}
\textit{\textbf{Numerical values} If $\sigma_\mathcal{F}=-1$, then $R_1\approx -8,330$ and $R_2 \approx 1,074$.\\ \hspace{3cm} If $\sigma_\mathcal{F}=1$, then $R_1\approx -1,074$ and $R_2 \approx 8,330$.}
\end{quote}
In addition, 
\begin{small}\begin{eqnarray*}
L(R)   & = & \underbrace{256(192-6\pi^2-2\pi^4)}_{<0}R^4(R-R_3)(R-R_4)
\end{eqnarray*}\end{small}with
\begin{small}\begin{eqnarray*}
R_3 & = & \min\left\{ \pi^2\frac{-24\sigma_\mathcal{F}+\sqrt{6\pi^2(\pi^2+3)}}{192-6\pi^2-2\pi^4}, \pi^2\frac{-24\sigma_\mathcal{F} -\sqrt{6\pi^2(\pi^2+3)}}{192-6\pi^2-2\pi^4} \right\} \\
R_4 & = & \max\left\{ \pi^2\frac{-24\sigma_\mathcal{F} +\sqrt{6\pi^2(\pi^2+3)}}{192-6\pi^2-2\pi^4}, \pi^2\frac{-24\sigma_\mathcal{F} -\sqrt{6\pi^2(\pi^2+3)}}{192-6\pi^2-2\pi^4} \right\}.
\end{eqnarray*}\end{small}
\begin{quote}
\textit{\textbf{Numerical values} If $\sigma_\mathcal{F}=-1$, then $R_3 \approx-8,210$ and $R_4\approx 0,573$. \\If $\sigma_\mathcal{F}=1$, then $R_3\approx -0,573$ and $R_4\approx 8,210$.}
\end{quote} 
Therefore, we always have
\begin{small}$$R_1 < R_3 < 0 <\frac{1}{2} < R_4 < R_2 .$$\end{small}Let $W(R)$ be degree $1$ monomial coefficient in $Y$, then
\begin{small}\begin{eqnarray*} W(R)   & = & 32\pi^2R^2(-3\pi^2-24\sigma_\mathcal{F} R+2(\pi^2-9)R^2).
\end{eqnarray*}\end{small}Since $0<R<1/2$, we have $\Delta(R)\geq 0$ and $L(R)\geq 0$. The polynomial $Y$ has two real roots,
\begin{small}$$  X_1(R)=\frac{-W(R)-\sqrt{\Delta(R)}}{2L(R)} \mbox{ and } X_2(R)=\frac{-W(R)+\sqrt{\Delta(R)}}{2L(R)} . $$ \end{small}Now, we want compare to $X_1(R)$, $X_2(R)$ and $B(g_0)^2$. Since $L(R)$ is positive , $X_1(R)\leq X_2(R)$. Moreover,
\begin{small}\begin{eqnarray*}
& &[B(g_0)^2-X_1(R)][B(g_0)^2-X_2(R)]\\
& = &\underbrace{\frac{\pi^2(\pi^4+3\pi^2-96)(8\pi^2R^2+24R^2+24\sigma_\mathcal{F} \pi^2R+3\pi^4)}{8R^2[\pi^2+8R\sigma_\mathcal{F}]^2[-3\pi^4-48\sigma_\mathcal{F} R\pi^2+2(\pi^4+3\pi^2-96)R^2]^2}}_{>0}(R-R_3)(R-R_4)  . 
\end{eqnarray*}\end{small}As a result, we have $X_1(R) \leq B(g_0)^2 \leq X_2(R) $ and $\beta_{\min}(g_0)^2=X_2(R)$.\\
Then, we prove that $Borne_\mathcal{F}(X)$ is increasing on $[X_2(R);+\infty[$. We compute
\begin{small}$$ \frac{ d Borne_\mathcal{F} }{dX}(X)=-\frac{256R^4B(g_0)^6}{3\pi^4(X-B(g_0)^2)^3}\left[16R^2(2R\sigma_\mathcal{F} (\pi^2-9)-3\pi^2)X-2R\sigma_\mathcal{F}(\pi^2+3)-3\pi^2  \right] . $$\end{small}If $0<R<1/2$, then $2R\sigma_\mathcal{F}(\pi^2-9)-3\pi^2 < 0$. As a consequence:
\begin{small}\begin{eqnarray*}
\frac{d Borne_\mathcal{F}}{dX}(X) \geq 0   & \Longleftrightarrow & X \geq \frac{1}{16R^2}\frac{2R\sigma_\mathcal{F} (\pi^2+3)+3\pi^2}{2R\sigma_\mathcal{F}(\pi^2-9)-3\pi^2}
\end{eqnarray*}\end{small}Since the right member of this inequality is always negative (when $0<R<1/2$), then $Borne_\mathcal{F}^+$ is an increasing function.
\begin{flushright}
$\Box$
\end{flushright}
\paragraph{Proof of theorem \ref{th:thm3}}
We specialise the inequality in proposition \ref{pro:proportioncasgeneral} in the test function $g_0$ with $0<R<\rho_{\max}(\mathcal{F})$. Using lemma \ref{lem:quelquescalculs}, if $\beta >B(g_0)$, we have
$$ \liminf_{\underset{ q \mbox{ \begin{tiny}prime\end{tiny}}}{ q\rightarrow +\infty}}\mu_{\mathcal{F}(q)}^h \left(\left\{ L(f,s)\in \mathcal{F}; \, \tilde{\gamma}_{f,1} < \beta \right\} \right) \geq Borne_\mathcal{F}^+(\beta^2) . $$
Thanks to lemma \ref{lem:variationdelaborne}, if $\beta > \beta_{\min}(g_0)$, we get
$$ \liminf_{\underset{ q \mbox{ \begin{tiny}prime\end{tiny}}}{ q\rightarrow +\infty}}\mu_{\mathcal{F}(q)}^h \left(\left\{ L(f,s)\in \mathcal{F}; \, \tilde{\gamma}_{f,1} < \beta \right\} \right) \geq Borne_\mathcal{F}(\beta^2). $$
To obtain theorem \ref{th:thm3} , we evaluate the right member in $R=\rho_{\max}(\mathcal{F})/2$.

\section{The smallest zero} \label{sec:pboptimisation}
In this part, we prove theorem \ref{th:thm1}. Our starting point is the following proposition. We do not write the proof that is essentially given in \cite{HR} (theorem 8.1).
\begin{proposition} \label{pro:alaplaceth1} The infimum being taken over functions $g$ in $\mathcal{S}_{\nu_{\max}(\mathcal{F})/2}(\R)$, we have
\begin{eqnarray} \limsup_{Q\rightarrow +\infty} \min_{L(f,s) \in \mathcal{F}(Q)} \tilde{\gamma}_{f,1}  \leq   \inf\left\{ \sqrt{ \frac{\int_{\R}  x^2g^2(x)W^*[\mathcal{F}](x)dx}{\int_{\R}g^2(x)W^*[\mathcal{F}](x) dx }}      \right\} .   \label{eqn:majorant}  \end{eqnarray}
\end{proposition}
We give a sketch of the proof of theorem \ref{th:thm1}.
\paragraph{Sketch of the proof}
We have to determine an explicit expression of the right member term in (\ref{eqn:majorant}). Let 
$$\m:=\inf_{g\in \mathcal{S}_{\nu}(\R)\backslash \{0\}} \frac{\int_{\R}  x^2g^2(x)W[G](x)dx}{\int_{\R}g^2(x)W[G](x) dx } $$
with $\nu =\nu_{\max}(\mathcal{F})/2$ and $W[G]=W^*[\mathcal{F}]$.
\begin{itemize}
\item In subsection \ref{ss:subdensitettop} (lemma \ref{lem:newexpression}), we prove the existence of a function $\tilde{B}$ and a set $H_{\nu}^{\infty}$ such that $\m=\inf_{\widehat{g} \in H_{\nu}^{\infty}\backslash \{0\}} \tilde{B}(\widehat{g})$. Precisely, we define $\tilde{B}$ by \begin{eqnarray} \label{eqn:defBdeh}
\tilde{B}(h) =  \frac{1}{4\pi^2}\frac{\int_{\R}h'(u)^2du -\frac{\delta}{2}\int_{-1}^{1}h'*h'(u)du}{\int_{\R}h(u)^2du+\frac{\delta}{2}\int_{-1}^{1}h*h(u)du + \varepsilon \left(\int_\R h(u)du  \right)^2}
\end{eqnarray}
where $(\delta,\varepsilon)$ is given in table \ref{tab:tableaudeltaepsilon}. 
\item In subsection \ref{ss:subdensitettop} (lemma \ref{lem:densityarg}), by a density argument, if 
$$ H_{R} = \left\{h : \R\rightarrow \R \mbox{ even, continuous, }\mathcal{C}^1 \mbox{ on }]-R,R[,\mbox{ with supp }h\subset [-R,R] \right\}$$ and   $$\tilde{\m}_R = \inf_{h \in H_{R} \backslash \{0\}} \tilde{B}(h),$$ we prove $\m=\lim_{R\rightarrow \nu^-}\tilde{\m}_R$.
\item In subsection \ref{ss:subdensitettop} (lemma \ref{lem:topoarg}), by topological arguments, there exists $\h_R$ in the Sobolev space $H_{0}^1$ such that $ \tilde{B}(\h_R)=\tilde{\m}_R$.
\item In subsection \ref{ss:sseV} (lemma \ref{lem:equationvolterra}), thanks to Fourier theory, $\h_R$ satisfies a Volterra equation with temporal shifts
\begin{eqnarray} \h_R(u)=\varphi(u)+\frac{\delta}{2}\int_u^R\h_R(t+1)-\h_R(t-1)dt  \label{eqn:eqnvoltpremiereappa} \end{eqnarray}
where $\varphi$ is an explicit function which is defined with the unknown parameter $\tilde{\m}_R$. Then, $\h_R$ is in $H_R$.
\item Then, we solve the previous Volterra equation in subsection \ref{ss:optimaltestfunction}. As a result, we obtain an explicit expression of $\h_R$ which also depends on the unknown parameter $\tilde{\m}_R$.
\item We determine $\sqrt{\tilde{\m}_R}$ in subsection \ref{ss:exactvaleur} by solving the equation $\tilde{\m}_R= \tilde{B}(\h_R)$. 
\item To conclude, we use $\sqrt{\m}=\lim_{R\rightarrow \nu^-}\sqrt{\tilde{\m}_R}$.
\end{itemize}

\subsection{Density and topological arguments} \label{ss:subdensitettop}

This section is devoted to the proof of the following lemma which sums up the three first steps of the preceding sketch of the proof. Let  $H_0^{1}$ be the Sobolev space defined by
\begin{small}
\begin{multline*}
H_0^1 = \left\{ u \in L^2(]-R,R[) \text{ with }u(-R)=u(R)=0 ,\right. \\ \left. \exists v\in L^2(\R) \mbox{ such that }\forall \phi \in \mathcal{C}_c^\infty(]-R,R[), \hspace{1mm}\int_\R u\phi'=-\int_\R v\phi  \right\}.
\end{multline*}
\end{small}If $u$ is in $H^1_0$ then $v$ is called the weak derivative of $u$ and is denoted by $u'$. The space $H^1_0$ is equipped with the inner product
$$ \langle u,v\rangle_{H^1_0}=\langle u',v'\rangle_{L^2}= \int_{-R}^Ru'v'.$$
$H^1_0$ is a reflexive separable  Hilbert space (\cite{B}, paragraphe VIII.3).
\begin{lemma} \label{lem:lemlimetmin} We have $$\m= \lim_{R\rightarrow \nu^-} \tilde{\m}_R.$$
In addition, for each $R>0$, there exists $\h_R$ in $H_0^1$ such that $\tilde{m}_R =\tilde{B}(\h_R)$.
\end{lemma}
Before proving this result, we need to prove some technical lemmas.

\paragraph*{A new expression}

\begin{lemma}\label{lem:newexpression} Let 
$ H_{\nu}^{\infty} = \left\{h \in \mathcal{C}_c^\infty(\R), \mbox{ even and with  supp }h\subset [-R,R] \mbox{ such that }0<R<\nu \right\}$.\\ Then,
$$\m=  \inf_{h \in H_{\nu}^{\infty}\backslash \{0\}} \tilde{B}(h) . $$
\end{lemma}
Proof: Thanks to Plancherel theorem,  inversion formula, Parseval formula and relation $(\ref{eqn:expressiontransformeefourierdensite})$, we have:
\begin{small}\begin{eqnarray*}
\int_{\R}g^2(x)W[G](x) dx & = &\int_{\R}\widehat{g^2}(y)\widehat{W[G]}(y)dy =\widehat{g^2}(0)+\frac{\delta}{2}\int_{-1}^{1}\widehat{g^2}(y)dy+\varepsilon\int_{\R}\widehat{g^2}(y)dy \\
  & = & \int_\R g^2(x)dx  +\frac{\delta}{2}\int_{-1}^{1}\widehat{g}(y)*\widehat{g}(y)dy+\varepsilon g^2(0)\\
  & = & \int_\R \widehat{g}^2(y)dy+\frac{\delta}{2}\int_{-1}^{1}\widehat{g}(y)*\widehat{g}(y)dy+ \varepsilon \left( \int_\R \widehat{g}(y)dy\right)^2
\end{eqnarray*}\end{small}Similarly, since $\widehat{xg(x)}=\frac{-1}{2i\pi}\widehat{g}'$, we prove
\begin{small}$$\int_{\R}x^2 g^2(x)W[G](x) dx=\frac{1}{4\pi^2}\int_{\R}\widehat{g}'(u)^2du -\frac{\delta}{2}\int_{-1}^{1}\widehat{g}'*\widehat{g}'(u)du . $$\end{small}In addition, since the set of Schwartz functions is invariant by Fourier transformation, we have:
\begin{small}$$g\in \mathcal{S}_\nu(\R) \Leftrightarrow \widehat{g} \in H_{\nu}^{\infty}$$\end{small}Therefore, \begin{small}$$\m = \inf_{\widehat{g} \in H_{\nu}^{\infty}\backslash \{0\}}\tilde{B}(\widehat{g}).$$\end{small}
\begin{flushright}
$\Box$
\end{flushright}

\paragraph*{A density argument}

\begin{lemma} \label{lem:densityarg} We have
$$\inf_{h \in H_{\nu}^{\infty}\backslash \{0\}} \tilde{B}(h)=\lim_{R\rightarrow \nu^-} \tilde{m}_R . $$
\end{lemma}
Proof:  We define mollifiers $(\rho_n)$ by
$$\rho_n:\left\{\begin{array}{l}\R\longrightarrow \R \\ x \longmapsto \frac{n}{\int_\R \rho(t)dt}\rho(nx) \end{array} \right.\hspace{2mm} \mbox{ where }\hspace{2mm} \rho(x) = \left\{ \begin{array}{cl} e^\frac{1}{|x|^2-1} & \mbox{if }|x|<1 \\ 0  & \mbox{if } |x|\geq 1 \end{array} \right. .$$
The function $\rho_n$ is  non-negative, smooth  with supp $\rho_n \subset [-1/n;1/n]$ and such that $\int_\R \rho_n(u)du=1$. We recall two properties (see e.g. \cite{B}, Theorem 4.22 and Theorem 4.15).
\begin{description}
\item[P1] Let $1\leq p <\infty$, if $g\in L^p(\R)$ then $\rho_n*g$ tends to $g$ in $L^p(\R)$.
\item[P2] Let $1\leq p \leq \infty$, if $f\in  L^1(\R)$ and $g\in L^p(\R)$ then $f*g\in L^p(\R)$ and $||f*g||_p\leq ||f||_1||g||_p$.
\end{description}
Let $0<R<\nu$  and $\eta >0$. \\
There exists $h$ in $H_R$ such that $ \tilde{m}_R  \leq \tilde{B}(h) \leq  \tilde{m}_R +\eta  $.
Let $h_n=\rho_n*h$. For large $n$, $h_n$ is in $H_{\nu}^{\infty} $. We also have $h_n'=\rho_n*h'$. Thanks to property \textbf{P1}, we have
$$\lim_{n\rightarrow +\infty}||h_n||_{2}^{2}=||h||_{2}^{2}, \hspace{4mm}  \lim_{n\rightarrow +\infty}||h_n'||_{2}^{2}=||h'||_{2}^{2}\hspace{4mm} \mbox{ and }\hspace{4mm}\lim_{n\rightarrow  +\infty}\int_{\R}h_n(u)du=\int_{\R}h(u)du .$$
Thanks to property \textbf{P2}, we may write
\begin{eqnarray*}& & \left|\int_{-1}^{1}h_n*h_n(u)-h*h(u)du\right| \leq   \int_\R  \left|h_n*h_n(u)-h*h(u)\right|du  \\
   & \leq &  ||(h_n-h)*h_n||_1+||h*(h_n-h)||_1 \leq  (||h_n||_1+||h||_1)||h-h_n||_1 .
\end{eqnarray*} 
Therefore,
$$  \lim_{n\rightarrow +\infty}\int_{-1}^{1}h_n*h_n(u)du =\int_{-1}^{1}h*h(u)du \mbox{ and similarly  }\lim_{n\rightarrow +\infty}\int_{-1}^{1}h_n'*h_n'(u)du =\int_{-1}^{1}h'*h'(u)du.$$
Then, $\lim_{n\rightarrow +\infty}\tilde{B}(h_n)= \tilde{B}(h)$. There exists $h_n$ in $H_{\nu}^{\infty} $ such that $|\tilde{B}(h)-\tilde{B}(h_n)| \leq \eta $.
As a result, for all $\eta>0$, there exists $h_n$ in $H_{\nu}^{\infty} $ such that $|\tilde{m}_R -\tilde{B}(h_n)| \leq 2\eta $. Then, for all $R<\nu$, we have $  \inf_{h\in H_{\nu}^{\infty}\backslash \{0\}} \tilde{B}(h) \leq \tilde{m}_R  $. As a consequence, we have
$$  \inf_{h\in H_{\nu}^{\infty}\backslash \{0\}} \tilde{B}(h) \leq \lim_{R\rightarrow \nu^-}\tilde{m}_R . $$
On the other side, let $h$ be in $H_{\nu}^{\infty}$. There exists $R_0<\nu$ such that supp $h \subset [-R_0;R_0]$. Since $h$ is in $H_{R_0}$ and $ R \mapsto \tilde{m}_R$ is decreasing, we get
$$  \lim_{R\rightarrow \nu^-} \tilde{m}_R \leq \tilde{m}_{R_0} \leq \tilde{B}(h) . $$
Therefore, $$\lim_{R\rightarrow \nu^-} \tilde{m}_R  \leq \inf_{h\in H_{\nu}^{\infty}\backslash \{0\}} \tilde{B}(h) .$$
\begin{flushright}
$\Box$
\end{flushright}

\paragraph*{Some compact operators}
Relation (\ref{eqn:defBdeh}) allows us to extend  $\tilde{B}$ to $H_0^1 \backslash \{0\}$.
\begin{lemma} \label{lem:topoarg} If $R>0$, there exists $\h_R$ in $H_{0}^1\backslash \{0\}$ such that $\tilde{m}_R=  \tilde{B}(\h_R) $.
\end{lemma}
Proof: Let $K$ be the operator of  $L^2(]-R,R[)$ defined by
$$ K[h](u) = \frac{\delta}{2}\int_{u-1}^{u+1}h(t)dt +\varepsilon \int_\R h(t)dt.$$
We may write $\tilde{B}$ on the shape 
$$ \tilde{B}(h) =\frac{\langle h',h' \rangle_{L^2}+\langle Kh',h' \rangle_{L^2}}{\langle h,h \rangle_{L^2}+\langle Kh,h \rangle_{L^2}}.$$
Since $K$ is a Hilbert-Schmidt operator, it is a compact operator of  $L^2(]-R,R[)$. Denote by $I$ the identity function of $L^2(]-R,R[)$, the transformations that  have been done in lemma \ref{lem:newexpression} show that, for all $h$ in $L^2(]-R,R[)$, 
$\langle (I+K)[h],h\rangle_{L^2} \geq 0$ with equality if and only if $h=0$. 

Since smooth function compactly supported in $]-R,R[$ are dense in $H_0^1$, the infimum  of $\tilde{B}$ over $H_0^1$ is also equal to $\tilde{m}_R$.

 Let $(g_n)$ be a sequence of non zero functions in $H_R$ such that
$\lim_{n\rightarrow +\infty}\tilde{B}(g_n)=\tilde{m}_R$. We consider $h_n=g_n / \langle g_n',g_n' \rangle_{L^2}$.  Since for all real number $t\neq 0$ and all $h$ in $H_R$, $\tilde{B}$ satisfies  $\tilde{B}(t.h)=\tilde{B}(h)$, we have
$$\lim_{n\rightarrow +\infty}\tilde{B}(h_n)=\tilde{m}_R \hspace{2mm} \mbox{ and }\hspace{2mm} \|h_n'\|_{L^2}= \|h_n\|_{H_0^1}=1.$$

Sequences $(h_n)$ and $(h_n')$ are bounded in the Hilbert space $L^2(]-R,R[)$. Since the unit ball of $L^2(]-R,R[)$ is compact for the weak topology, up to consider sub-sequences, there exists   $h$ and $k$ in $L^2(]-R,R[)$ such that $(h_n)$ (respectively $(h_n')$) tends to $h$ (respectively $k$) weakly. Moreover, for all function $\phi$ in $\mathcal{C}_c^\infty(]-R,R[)$, we may write
$$ \int g\phi' = \lim_{n\rightarrow +\infty}\int g_n\phi' =-\lim_{n\rightarrow +\infty}\int g_n'\phi =- \int k\phi. $$
Then $g$ belongs to $H^1$, $g'=k$ and $g_n$ tends to $g$ weakly in $H^1$. In addition, there exists a compact embedding of $H^1$ into the set of continuous function on $[-R,R]$ equipped with the norm of uniform convergence (\cite{B}, theorem 8.2 et 8.8), we may conclude (\cite{B}, 6.1 remark 2)  $h_n$ converges uniformly to $h$ (precisely to its continuous representative which will be always identified with $h$). As a result, $h$ is an even functions belonging to $H_0^1$ and $\lim_{n \rightarrow +\infty} \|h_n\|_{L^2}=\|h\|_{L^2}$.

Since $(h_n)$ and $(h_n')$ are weakly convergent in $L^2(]-R,R[)$ and since $K$ is a compact operator,  the sequence  $(Kh_n)$ and  $(Kh_n')$ converge respectively to $Kh$ and $Kh'$  strongly in $L^2(]-R,R[)$. Then, we get (\cite{B}, proposition 3.5 (iv))
$$\lim_{n\rightarrow +\infty}\langle Kh_n,h_n\rangle_{L^2} = \langle Kh,h\rangle_{L^2} \mbox{ and }\lim_{n\rightarrow +\infty}\langle Kh_n',h_n'\rangle_{L^2} = \langle Kh',h'\rangle_{L^2}.$$
To sum up, we have 
\begin{align*}
\lim_{n\rightarrow +\infty}\langle (I+K)h_n',h_n'\rangle_{L^2}   = 1+  \langle Kh',h'\rangle_{L^2}, \hspace{4mm}
\lim_{n\rightarrow +\infty}\langle (I+K)h_n,h_n\rangle_{L^2}  =  \langle (I+K)h,h\rangle_{L^2} 
\end{align*} and, since $(h_n')$ tends weakly to $h'$ in $L^2$, we may deduce $\|h'\|_{L^2} \leq \liminf \|h_n'\|_{L^2}=1$. 

Furthermore, $1+  \langle Kh',h'\rangle_{L^2} $ is non zero. Indeed, if it was zero, we would have
$$0=1+  \langle Kh',h'\rangle_{L^2} \geq \|h'\|_{L^2}+  \langle Kh',h'\rangle_{L^2} = \langle (I+K)h',h'\rangle_{L^2} \geq 0.$$
So, $\|h'\|_{L^2}=1$ and $\langle (I+K)h',h'\rangle_{L^2}=0$, thus $\|h'\|_{L^2}=1$ and $h'$ should be zero ! 

As a result $\langle (I+K)h,h\rangle_{L^2}\neq 0$ because $\tilde{m}_R$ is finite and $1+  \langle Kh',h'\rangle_{L^2} $  is non zero.

Then
$$\frac{\|h'\|_{L^2}+ \langle Kh',h'\rangle_{L^2}}{\langle (I+K)h,h\rangle_{L^2}} = \tilde{B}(h) \geq \tilde{m}_R = \lim_{n\rightarrow +\infty} \tilde{B}(h_n)=\frac{1+ \langle Kh',h'\rangle_{L^2}}{\langle (I+K)h,h\rangle_{L^2}}, $$
and $\|h'\|_{L^2} \geq 1$. Since we already have $\|h'\|_{L^2} \leq 1$, it comes  $\|h'\|_{L^2} = 1$. To conclude, the function $h$ is non zero and satisfies $\tilde{B}(h) = \tilde{m}_R$.\begin{flushright}
$\Box$
\end{flushright}
\begin{remark} $\tilde{B}$ is not continuous when $B_{||\cdot ||_{H^1}}(0,1) \rightarrow \R $ is equipped with the weak topology. Indeed, since the unit ball is weakly compact, $\tilde{B}$ should be bounded. However, by considering $h_n(u)=\cos\left( \frac{(2n+1)\pi u}{2R}\right)\un_{[-R,R]}(u)$, we remark that $\tilde{B}$ cannot be bounded. 
\end{remark}
\paragraph*{}
The proof of lemma \ref{lem:lemlimetmin} comes from lemmas \ref{lem:newexpression}, \ref{lem:densityarg} and \ref{lem:topoarg}.

\subsection{Fourier analysis} \label{ss:ssFourier}

For technical conveniences, let $B(h)=16R^2 \tilde{B}(h)$ and $\m_R = \inf_{h\in H_R}B(h)$. Then $$\m_R= B(\h_R) = 16R^2\tilde{\m}_R .$$
Let $$\Omega_R=\left\{(c_n)\in \R^\N,\mbox{ such that }c_{2n}=0\mbox{ and }x\mapsto \sum_{n \geq 0}c_n \cos\left( \frac{\pi n x}{2R}\right) \mbox{ is continuous and in}H_0^1   \right\} . $$
In this section, we prove some technical lemmas.
\begin{lemma} \label{lem:transformationBavecSFpart1} 
Let $$\Psi:\left\{ \begin{array}{lcl }H_0^1\cap \{continuous , even \} & \longrightarrow & \Omega_R \\ h & \longmapsto & \left( \frac{1-(-1)^n}{2R}\int_{0}^{R}h(t)\cos\left(\frac{\pi n t}{2R} \right)dt \right)_{n\geq 0} \end{array}  \right. . $$ 
Then, $\Psi$ is a bijective function.
\end{lemma}
Proof: For each $h$ in $H_0^1$, we associate an even $4R$-periodic function $\tilde{h}$ which is defined by:
$$ \mbox{for all } x \in [0,R], \hspace{2mm}\tilde{h}(x)= h(x) \mbox{ and } \tilde{h}(R+x)=-h(R-x). $$ 
We also have $h=\tilde{h}\cdot \un_{[-R,R]}$.
\begin{center}
\begin{tabular}{cc}
\includegraphics[scale=0.2]{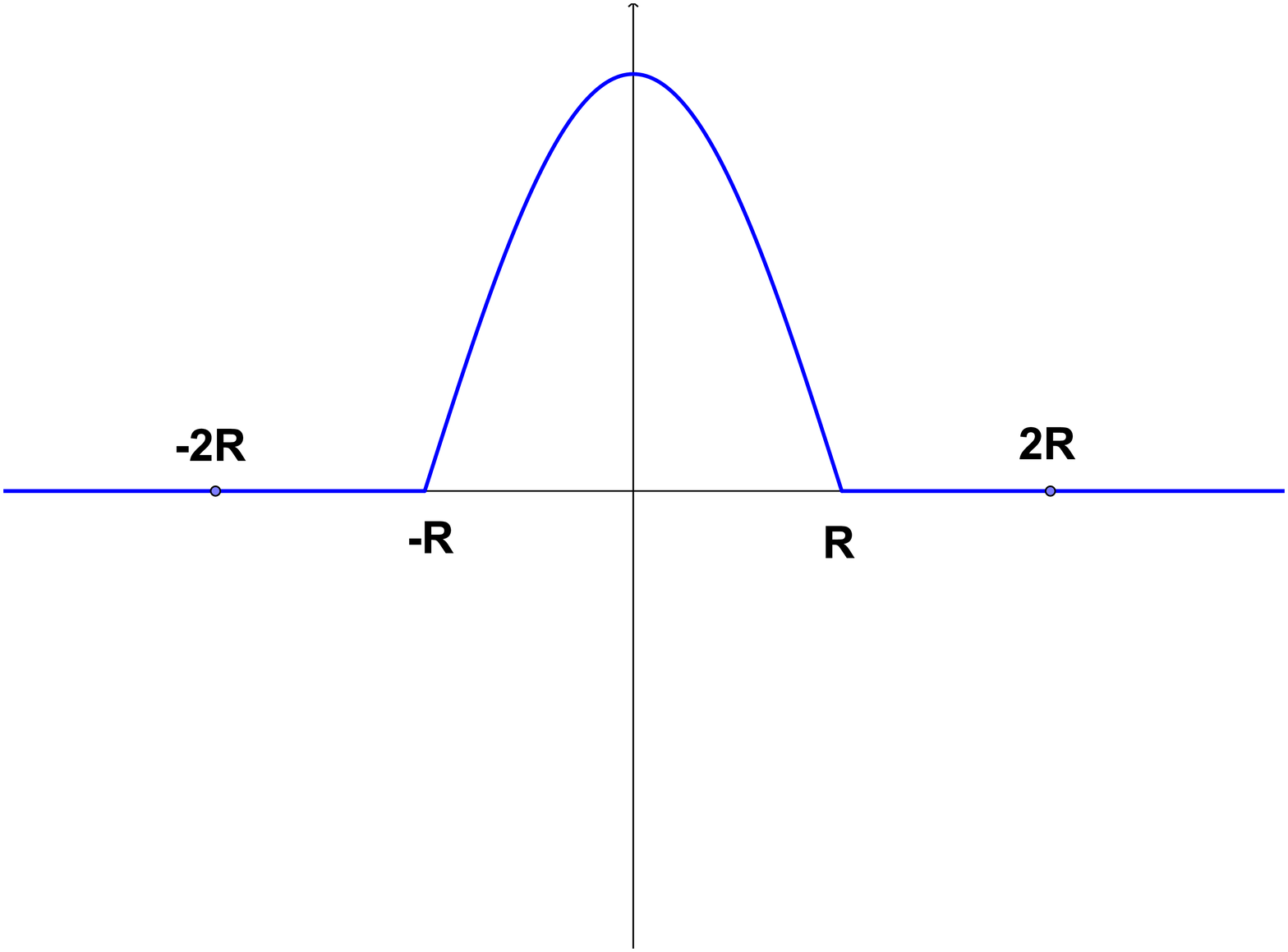} & \includegraphics[scale=0.2]{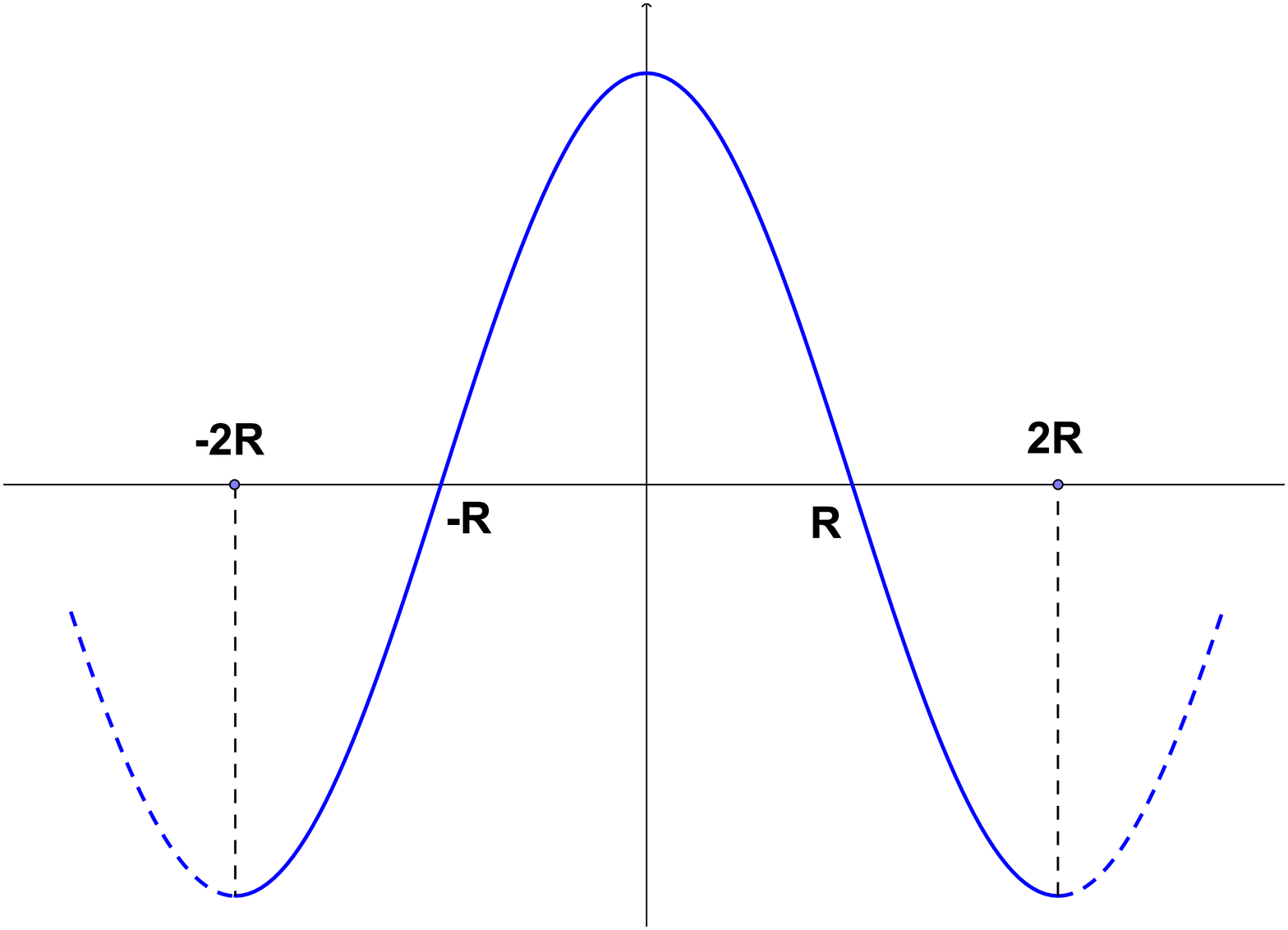} \\
representative curve of $h$ & representative curve of $\tilde{h}$
\end{tabular}
\end{center}
As a result, the mean value of $\tilde{h}$ is equal to zero. Actually, if $c_n$ denotes the $n$-th Fourier coefficient of $\tilde{h}$, we have
$$c_n=\frac{1}{4R}\int_{-2R}^{2R}\tilde{h}(t)e^{-\frac{i\pi n t}{2R}}dt = \frac{1-(-1)^n}{2R}\int_{0}^{R}h(t)\cos\left(\frac{\pi n t}{2R} \right)dt. $$ 
Since the Fourier series of a function in $H_0^1$ is normally convergent, we get
$$h(u)=\tilde{h}(u)\cdot \un_{[-R,R]}(u) = 2\left[\sum_{n\geq 0} c_n  \cos\left( \frac{\pi n x}{2R}\right)\right]\un_{[-R,R]}(u).$$
Then $\Psi$ is well defined and the last relation gives us an explicit expression of $\Psi^{-1}$. 
\begin{flushright}
$\Box$
\end{flushright}
As a consequence, since $B$ is defined on the set of continuous and even functions in $H_0^1$, we may define $B$ on $\Omega_R$ by $B(c):=B(\Psi^{-1}(c))$. 
\begin{lemma} \label{lem:transformationBavecSF} Let $c=(c_n)$ be in $\Omega_R$. If $R\leq 1/2$ then
$$B(c) = \frac{\sum_{n\geq 0}^{*}n^2c_n^2 }{\sum_{n\geq 0}^{*}c_n^2+(\delta+2\varepsilon)\frac{8R}{\pi^2}\left( \sum_{n\geq 0}^{*} \frac{ (-1)^{ \frac{n-1}{2}} }{n}   c_n \right)^2 } $$
and, if $R>1/2$, we have
$$B(c) = \frac{\sum_{n\geq 0}^{*} n^2c_n^2-\frac{\delta}{2R}\sum\sum_{m,n\geq 0}mn\mu_{m,n}c_mc_n}{\sum_{n\geq 0}^{*}c_n^2+\frac{\delta}{2R}\sum\sum_{m,n\geq 0}\lambda_{m,n}c_mc_n+\frac{16R\varepsilon}{\pi^2}\left( \sum_{n\geq 0}^* \frac{(-1)^{\frac{n-1}{2}}c_n}{n} \right)^2}$$ 
where, for $n$ and $m$  odd numbers:
$$ \lambda_{m,n}=\left\{\begin{array}{lc} \frac{8R^2mn(-1)^{\frac{m+n}{2}}}{\pi^2(m^2-n^2)}\left[ \frac{1}{n^2}\cos \frac{\pi n}{2R}-\frac{1}{m^2}\cos \frac{\pi m}{2R} \right]-\frac{8R^2(-1)^{\frac{m+n}{2}}}{mn\pi^2} & \mbox{si }  m\neq n \\
                                        \frac{2R(2R-1)}{n\pi }\sin \frac{\pi n}{2R}-\frac{8R^2}{\pi^2n^2}\cos \frac{n\pi}{2R}+\frac{8R^2}{\pi^2n^2} & \mbox{si }  m=n
                     \end{array} \right.   $$
$$\mbox{and }\hspace{1cm}\mu_{m,n}=\left\{\begin{array}{lc} 
 \frac{8R^2(-1)^{\frac{m+n}{2}}}{\pi^2(m^2-n^2)}\left[ \cos \frac{\pi m}{2R}-\cos \frac{\pi n}{2R} \right]& \mbox{si }  m\neq n \\
 -\frac{2R(2R-1)}{n\pi }\sin \frac{\pi n}{2R}  & \mbox{si }  m=n
 \end{array} \right.   $$                     
\end{lemma} 
Proof: We are giving an explicit expression of $B$ on $\Omega_R$.  Let $h$ be an even and continuous function in $H_0^1$ such that $c=(c_n)=\Psi (h)$. Thanks to Parseval formula, we have:
\begin{small}\begin{eqnarray*}
\frac{1}{2R}\int_{\R}h(u)^2du  & = & \frac{1}{4R}\int_{-2R}^{2R}\tilde{h}(u)^2du  =  \sum_{n\in \Z} |c_n|^2 =2\sum_{n\geq 0}c_n^2\label{eqn:perse1} \\
\frac{1}{2R}\int_{\R}h'(u)^2du &  = & \frac{1}{4R}\int_{-2R}^{2R}\tilde{h}'(u)^2du =  \frac{\pi^2}{4R^2}\sum_{n\in \Z}n^2|c_n|^2 = \frac{\pi^2}{2R^2}\sum_{n\geq 0}n^2c_n^2 \label{eqn:perse2}
\end{eqnarray*} \end{small}In addition, we may write
\begin{small}\begin{eqnarray*}
\int_{\R}h(u)du = \int_{-R}^{R}\tilde{h}(u)du = \frac{4R}{\pi}\sum_{n \in \Z}\frac{(-1)^n}{(2n+1)}c_{2n+1}=\frac{8R}{\pi}\sum_{n\geq 0}\frac{(-1)^n}{(2n+1)}c_{2n+1}. \label{eqn:perse3}
\end{eqnarray*}\end{small}Remember $B(h)=16R^2 \tilde{B}(h)$ where  $\tilde{B}(h)$ is given in (\ref{eqn:defBdeh}). We obtain an explicit expression of $B(c)$ using
\begin{small}$$\int_{-1}^{1}h*h(x)dx = \left\{ \begin{array}{ll} \left( \int_\R h(x)dx \right)^2 & \mbox{ if }R\leq \frac{1}{2} \\ 4\sum_{m\geq 0}\sum_{n\geq 0}c_m c_n \lambda_{m,n} & \mbox{ if }R  > \frac{1}{2} \end{array}\right.$$\end{small}and
\begin{small}$$\int_{-1}^{1}h'*h'(x)dx = \left\{ \begin{array}{ll} 0 & \mbox{ if }R\leq \frac{1}{2}\\ \frac{\pi^2}{R^2}\sum_{m\geq 0}\sum_{n\geq 0}mnc_m c_n \mu_{m,n} & \mbox{ if }R  > \frac{1}{2} \end{array}\right. . $$ \end{small}
\begin{flushright}
$\Box$
\end{flushright}
Hughes and Rudnick proved that $\m=\frac{1}{16 \nu^2}$ if $G=U$. We extend their result.
\begin{corollary} \label{cor:taillem}If $G=O$ or if $G=SO^\pm$ and $0<R \leq 1/2$  then
$$ \frac{1}{1+R}<\m_R \leq \frac{1}{1+\frac{8}{\pi^2}R} < 1 .$$
If  $G=SO^\pm$ with $R>1/2$, then
$   \m_R < 1$. If $G=Sp$  and $0<R<1/2$, then
$$ 1 < \m_R < \frac{1}{1-\frac{8}{\pi^2}R} < 2 . $$
\end{corollary}
Proof: Thanks to the Cauchy-Schwarz inequality,  if $G=O$ or if $G=SO^\pm$ with $0<R \leq 1/2$, we have
\begin{small}$$B(c) = \frac{\sum_{n\geq 0}^{*}n^2c_n^2 }{\sum_{n\geq 0}^{*}c_n^2+(\delta+2\varepsilon)\frac{8R}{\pi^2}\left( \sum_{n\geq 0}^{*} \frac{ (-1)^{ \frac{n-1}{2}} }{n}   c_n \right)^2 } \geq \frac{\sum_{n\geq 0}^{*}n^2c_n^2 }{(1+R)\sum_{n\geq 0}^{*}c_n^2 }\geq \frac{1}{1+R} $$\end{small}and one of both inequalities is strict.
In order to obtain an upper bound for $\m_R$, we specialise $B$ in $e=(e_n)$ which is defined by $e_n=0$ for all $n$ except $e_1=1$. The upper bound comes easily if $G=O$, if $G=SO^\pm$ or if $G=Sp$ and $0<R \leq 1/2$. If $G=SO^\pm$ and $R>1/2$, we have
\begin{small}$$B(e)=\frac{1-\frac{\delta}{2R}\mu_{1,1}}{1+\frac{\delta}{2R}\lambda_{1,1}+\frac{16R\varepsilon }{\pi^2}}=  \frac{1+\frac{\delta(2R-1)}{\pi}\sin \frac{\pi}{2R} }{1+\frac{\delta(2R-1)}{\pi}\sin \frac{\pi}{2R}+\frac{4R\delta }{\pi^2}\left(1-\cos \frac{\pi}{2R} \right)+\frac{16R \varepsilon}{\pi^2}} < 1 . $$\end{small}Finally, if $G=Sp$ and $R<1/2$, for all $c$ in $\Omega_R$, we have
\begin{small}$$B(c) = \frac{\sum_{n\geq 0}^{*}n^2c_n^2 }{\sum_{n\geq 0}^{*}c_n^2-\frac{8R}{\pi^2}\left( \sum_{n\geq 0}^{*} \frac{ (-1)^{ \frac{n-1}{2}} }{n}   c_n \right)^2 } \geq \frac{\sum_{n\geq 0}^{*}n^2c_n^2 }{\sum_{n\geq 0}^{*}c_n^2 } \geq 1 $$\end{small}and  one of both inequalities is strict.
Therefore, thanks to lemma \ref{lem:topoarg}, we have $\m_R  > 1$. 
\begin{flushright}
$\Box$
\end{flushright}
Corollary \ref{cor:taillem} implies, except if $G=Sp$ and $R>1/2$, that $\m_R$ is not the square of an integer. We make the following hypothesis which will be justified page \pageref{lem:pourref}.\begin{center}
\textbf{``$Sp$ hypothesis'':} If $G=Sp$ and $R>1/2$ then $\m_R$ is not the square of an integer.
\end{center}
\begin{remark} In the symplectic case with $R>1/2$, by evaluating $B$ in $e$, we get:
$$\m_R \leq B(e)=  \frac{1-\frac{2R-1}{\pi}\sin \frac{\pi}{2R} }{1-\frac{2R-1}{\pi}\sin \frac{\pi}{2R}-\frac{4R }{\pi^2}\left(1-\cos \frac{\pi}{2R} \right)}= \frac{12 R}{\pi^2} + O(1) $$
Therefore, since $\m_R = 16R^2 \tilde{\m}_R$, corollary \ref{cor:taillem} and proposition \ref{pro:alaplaceth1} show:
\begin{eqnarray*} \limsup_{Q\rightarrow +\infty} \min_{L(f,s) \in \mathcal{F}(Q)} \tilde{\gamma}_{f,1}  \leq \lim_{R\rightarrow \nu_{\max}(\mathcal{F})/2} \sqrt{\tilde{\m}_R} \ll \left\{ \begin{array}{ll} \frac{1}{\nu_{\max}(\mathcal{F})^{1/2}} & \mbox{if } W^*[\mathcal{F}]=W[Sp] \\  \frac{1}{\nu_{\max}(\mathcal{F})} & \mbox{if } W^*[\mathcal{F}]=W[SO^+] \\  \frac{1}{\nu_{\max}(\mathcal{F})^{3/2}} & \mbox{if } W^*[\mathcal{F}]=W[O] \\  \frac{1}{\nu_{\max}(\mathcal{F})^2} & \mbox{if } W^*[\mathcal{F}]=W[SO^-]  \end{array}  \right.  \end{eqnarray*}
As a result, density conjecture agrees with random matrix prediction.
\end{remark}

\subsection{A Volterra equation with temporal shifts} \label{ss:sseV}
Since the unitary case has been solved, we assume $G$ is one of the following compact groups $O$, $SO^-$, $SO^+$ or $Sp$.
In this section, we prove the following lemma.
\begin{lemma}  \label{lem:equationvolterra} The optimal test function $\h_R$ satisfies, for all $0\leq u \leq R$, 
\begin{eqnarray} \h_R(u)=\varphi(u) +\frac{\delta}{2}\int_{u}^{R}\h_R(t+1)-\h_R(t-1)dt \label{eqn:reftardive} \end{eqnarray}
where $(\delta,\varepsilon)$ is given in the table \ref{tab:tableaudeltaepsilon} and with 
$$\varphi(u) = \frac{\pi k_{\delta,\varepsilon}}{2 \m_R \cos\frac{\sqrt{\m_R}\pi }{2}}\left[\cos\frac{\sqrt{\m_R}\pi u}{2R} - \cos\frac{\sqrt{\m_R}\pi }{2}\right] \cdot \un_{[-R,R]}(u) $$
where 
\begin{eqnarray}  k_{\delta,\varepsilon} = \frac{2\m_R}{\pi}\left[\frac{\delta}{2}\int_{R-1}^{R} \h_R (x)dx + \varepsilon\int_{-R}^{R} \h_R(x)dx \right] . \label{eqn:equationfinalepremiereapparition}
 \end{eqnarray}
\end{lemma}
\begin{remark} We identify $\h_R$ to its continuous representative. As a result,  equation (\ref{eqn:reftardive}) proves that $\h_R$ is of class $\mathcal{C}^1$ on $]-R,R[$. In other words,  the function $\h_R$ belongs to $H_R$.
\end{remark}

Before proving this result, we need to prove some technical lemmas. However, we may immediately deduce the following corollary.
\begin{corollary} \label{cor:remasimplificatrice} If $G=O$ or if $G=SO^\pm$ or $Sp$ with $R \leq 1/2$, then 
$$\sqrt{\m} = \frac{1}{\nu}V^{-1}\left( 1+ (\delta+2\varepsilon)\frac{1}{\nu} \right) $$ 
where $V$ has been defined in section ``Preview of results''.
\end{corollary}
Proof: If $G=O$, then $\delta=0$ and lemma \ref{lem:equationvolterra} gives us immediately $\h_R=\varphi$. Similarly, if $G=SO^+$, $SO^-$ or $Sp$, $R\leq 1/2$ and if $0< t< R$  then $t+1 \notin [-R,R]$ and $t-1 \notin [-R,R]$. Therefore, in this case, we also have $\h_R=\varphi$. We may sum up these remarks by
\begin{small}
$$\h_R(u)=\varphi(u)=\frac{(\delta+2\varepsilon)\int_{-R}^{R}\h_R(x)dx}{2\cos \left( \frac{\pi \sqrt{\m_R}}{2}\right)} \left[\cos\frac{\sqrt{\m_R}\pi u}{2R} - \cos\frac{\sqrt{\m_R}\pi }{2}\right] \cdot \un_{[-R,R]}(u) . $$
\end{small}As a result, we have 
\begin{small}
$$\int_{-R}^{R}\h_R(u)du \cdot \left[1-(\delta+2\varepsilon)R\left(\frac{2}{\pi \sqrt{\m_R}} \tan \left(\frac{\pi \sqrt{\m_R}}{2} \right)-1 \right)  \right] =  0 . $$
\end{small}Since $\int_{-R}^{R}\h_R(u)du \neq 0$ (otherwise $\h_R=0$), we may conclude
\begin{small}
$$\frac{2}{\pi \sqrt{\m_R}} \tan \left(\frac{\pi \sqrt{\m_R}}{2} \right) = 1 +\frac{\delta+2\varepsilon}{R} .$$
\end{small}The result comes easily from lemmas \ref{lem:newexpression} and \ref{lem:densityarg}.
\begin{flushright}
$\Box$
\end{flushright}
\begin{remark}  We will use this phenomenon in the next section in order to determine $\h_R$ in full generality.
\end{remark}
\paragraph*{}
Let
$$S_n(t)= \sin \left( \frac{\pi n t}{2R}\right)\un_{[-R,R]}(t) \hspace{3mm} \mbox{and } \hspace{3mm} C_n(t)= \cos \left( \frac{\pi n t}{2R}\right)\un_{[-R,R]}(t) . $$ 
\begin{lemma} \label{lem:lemmecalculsommemuetlambda} $n$ denotes an odd positive integer. Let $h$ be in $H_R$ with $c=(c_m)=\Psi(h)$, then:
\begin{eqnarray*}
 & & \sum_{m\geq 0}\left[mn\mu_{m,n} + \m_R \lambda_{m,n}\right]c_m \\
 & & \hspace{19mm}= -\frac{2R}{n\pi}(n^2-\m_R)(S_n*h)(1)+\frac{8R^2(-1)^{\frac{n-1}{2}}}{n\pi^2}\m_R\sum_{m\geq 0}\frac{(-1)^{\frac{m-1}{2}}}{m}\left[1-\cos \frac{m\pi}{2R}\right]c_m 
\end{eqnarray*}
\end{lemma}
Proof: We have:
\begin{small}
\begin{eqnarray*}
& & \sum_{m\geq 0}\left[mn\mu_{m,n} + \m_R \lambda_{m,n}\right]c_m \\
 & &= -(n^2-\m_R)\left( \frac{2R(2R-1)}{n\pi}\sin\left(\frac{\pi n}{2R}\right)c_n+\frac{8R^2}{n\pi^2}\sum_{m\neq n}\frac{m(-1)^{\frac{m+n}{2}}}{m^2-n^2}\cos\left(\frac{\pi n}{2R}\right)c_m \right)\\
 & & -\frac{8R^2\m_R}{n^2\pi^2}\cos\left(\frac{\pi n}{2R}\right)c_n-\m_R\sum_{m\geq 0}\frac{8R^2(-1)^{\frac{m+n}{2}}}{nm\pi^2}c_m+\frac{8R^2}{\pi^2}\sum_{m\neq n}\frac{n(-1)^{\frac{m+n}{2}}}{m(m^2-n^2)}(m^2-\m_R)\cos\left(\frac{\pi m}{2R}\right)c_m 
\end{eqnarray*}
\end{small}Since \begin{small}$$\frac{n}{m(m^2-n^2)}= \frac{1}{n}\left(\frac{m}{m^2-n^2}-\frac{1}{m}\right) ,$$\end{small}we may conclude
\begin{small}
\begin{eqnarray*}
& & \sum_{m\geq 0}\left[mn\mu_{m,n} + \m_R \lambda_{m,n}\right]c_m \\
& & = -(n^2-\m_R)\left( \frac{2R(2R-1)}{n\pi}\sin\left(\frac{\pi n}{2R}\right)c_n+\frac{8R^2}{n\pi^2}{\sum_{m\neq n}}^*\frac{(-1)^{\frac{m+n}{2}}m}{m^2-n^2}\left[\cos\left(\frac{\pi n}{2R}\right)-\cos\left(\frac{\pi m}{2R}\right)\right] c_m \right)\\
& & \hspace{65mm} +\frac{8R^2\m_R(-1)^{\frac{n-1}{2}}}{n\pi^2}{\sum_{m\geq 0}}^{*}\frac{(-1)^{\frac{m-1}{2}}}{m}\left[1-\cos\left(\frac{\pi m}{2R}\right)\right]c_m .
\end{eqnarray*}
\end{small}In addition, we have
\begin{small}
\begin{eqnarray*}
& & (S_n*h)(1)= \int_{1-R}^{R}h(u)\sin \left(\frac{\pi n}{2R}(1-u) \right) du \\
& &= {\sum_{m\geq 0}}^* c_m \int_{1-R}^{R}\sin \left(\frac{\pi (m-n) u}{2R}+\frac{\pi n}{2R} \right)- \sin \left(\frac{\pi (m+n)u}{2R}-\frac{\pi n}{2R} \right) du \\
& & =(2R-1)\sin \left( \frac{\pi n}{2R}\right)c_n+\frac{4R}{\pi}{\sum_{m\neq n}}^*\frac{(-1)^{\frac{m+n}{2}}m}{m^2-n^2}\left[\cos\left(\frac{\pi n}{2R}\right)-\cos\left(\frac{\pi m}{2R}\right)\right] c_m .
\end{eqnarray*}
\end{small}
\begin{flushright}
$\Box$
\end{flushright}

\begin{lemma}\label{lem:pour Sphypo} Let $\c=(\c_n)=\Psi(\h_R)$. Then, for all odd positive integer $n$, we have
$$\c_n =  \frac{ (-1)^\frac{n-1}{2}k_{\delta,\varepsilon}}{n[n^2-\m_R]}-\frac{\delta}{n\pi}(S_n*\h_R)(1) =\frac{ (-1)^\frac{n-1}{2}k_{\delta,\varepsilon}}{n[n^2-\m_R]}+\frac{2R\delta}{n^2\pi^2}(C_n*\h_R ')(1) . $$
\end{lemma}
Proof: We define a norm on $\Omega_R$ by
\begin{small}$$||c||=\sum_{n=0}^{+\infty}n |c_n| +\left( \sum_{n=0}^{+\infty}n^2c_n^2 \right)^{1/2}$$\end{small}and we consider five differentiable functions $H$, $T$, $S$, $A$ and $G$ on $(\Omega_R, ||\cdot ||)$ defined by:
\begin{small}\begin{eqnarray*}
& & H(c)=\sum_{n\geq 0}n^2 c_n^2, \hspace{11mm} T(c)=\sum_{n\geq 0}c_n^2,\hspace{11mm} S(c) = {\sum_{n\geq 0}}^*\frac{(-1)^{\frac{n-1}{2}}}{n}c_n ,\\
& & A(c)=\sum_{m\geq 0}\sum_{n\geq 0}mn\mu_{m,n}c_mc_n \hspace{3mm}\mbox{ and }\hspace{14mm}G(c)=\sum_{m\geq 0}\sum_{n\geq 0}\lambda_{m,n}c_mc_n
. \end{eqnarray*}\end{small}With lemma \ref{lem:transformationBavecSF}, we may write 
\begin{small}$$B = \frac{H-\frac{\delta}{2R}A}{T+\frac{\delta}{2R}G+\frac{16R\varepsilon}{\pi^2}S^2} \mbox{ if }R>\frac{1}{2} \hspace{4mm}\mbox{or} \hspace{4mm} B= \frac{H}{T+(\delta+2\varepsilon) \frac{8R}{\pi^2}S^2}\mbox{ if }R\leq \frac{1}{2}. $$\end{small}Since $B(\c)=\m_R$, $\c$ is a critical point. Therefore, if $R>1/2$, we have $dB_{|_\c}=0$ hence:
\begin{small}$$ \forall n \mbox{ odd}, \hspace{3mm} (n^2-\m_R) \c_n= \frac{\delta}{2R}\sum_{m\geq 0}\left[mn\mu_{m,n} + \m_R\lambda_{m,n}\right]\c_m +\varepsilon \frac{16RS(\c)\m_R}{\pi^2}\frac{(-1)^{\frac{n-1}{2}}}{n} $$\end{small}Thanks to lemma \ref{lem:lemmecalculsommemuetlambda}, ``$Sp$ hypothesis'' and corollary \ref{cor:taillem}, we may conclude
\begin{small}\begin{eqnarray}  \label{eqn:expressioncn} \c_n =  \frac{ (-1)^\frac{n-1}{2}k_{\delta,\varepsilon}}{n[n^2-\m_R]}-\frac{\delta}{n\pi}(S_n*\h_R)(1) =\frac{ (-1)^\frac{n-1}{2}k_{\delta,\varepsilon}}{n[n^2-\m_R]}+\frac{2R\delta}{n^2\pi^2}(C_n*\h_R ')(1)\end{eqnarray}\end{small}where 
\begin{small}\begin{eqnarray}  k_{\delta,\varepsilon}&  = & \frac{4R\m_R}{\pi^2}\left[\delta{\sum_{m\geq 0}}^*\frac{(-1)^{\frac{m-1}{2}}}{m}\left( 1-\cos\frac{m\pi}{2R}\right)\c_m+4\varepsilon{\sum_{m\geq 1}}^*\frac{(-1)^{\frac{m-1}{2}}}{m}\c_m  \right]\nonumber \\
 & = & \frac{2\m_R}{\pi}\left[\frac{\delta}{2}\int_{R-1}^{R} \h_R(x)dx + \varepsilon\int_{-R}^{R} \h_R(x)dx \right] . \label{eqn:constantemagique}
 \end{eqnarray}\end{small}Similarly, $\c$ satisfies (\ref{eqn:expressioncn}) and (\ref{eqn:constantemagique}) in the case $R\leq 1/2$ due to the fact that, if $R\leq 1/2$, we have $(S_n*\h_R)(1)=0$ and $\int_{R-1}^{R} \h_R(x)dx=\int_{-R}^{R} \h_R(x)dx$. 
\begin{flushright}
$\Box$
\end{flushright}
\paragraph*{Proof of lemma \ref{lem:equationvolterra}}
Assume $0\leq u \leq R$. Since $\h_R(u)=2\sum_{n\geq 0}\c_n \cos\left( \frac{\pi n u }{2R} \right)$, lemma \ref{lem:pour Sphypo} gives
\begin{small}$$ \h_R(u) = 2k_{\delta,\varepsilon}{\sum_{n\geq 0}}^*\frac{(-1)^{\frac{n-1}{2} }\cos\left( \frac{\pi n u }{2R} \right)}{n(n^2-\m_R)}+\frac{4R\delta}{\pi^2}{\sum_{n\geq 0}}^*\frac{(C_n*\h_R')(1)}{n^2}\cos\left( \frac{\pi n u }{2R} \right) . $$\end{small}First, using relation 1.444.6 of \cite{GR}, we have
\begin{small}
\begin{eqnarray*}
 {\sum_{n\geq 0}}^*\frac{(C_n*\h_R')(1)}{n^2}\cos\left( \frac{\pi n u }{2R} \right)  & = & \frac{\pi}{8}\int_{-R}^{R}\left( \pi -\frac{\pi |u+t|}{2R} -\frac{\pi |u-t|}{2R} \right)\h_R'(1-t)dt \\ & = &  \frac{\pi^2}{8R}\int_{u}^{R}\h_R(t+1)-\h_R(t-1)dt .
\end{eqnarray*}
\end{small}Then, using relation 1.445.6 of \cite{GR}, we compute
\begin{small}$$2k_{\delta,\varepsilon}{\sum_{n\geq 0}}^*\frac{(-1)^{\frac{n-1}{2} }\cos\left( \frac{\pi n u }{2R} \right)}{n(n^2-\m_R)}  = \frac{\pi k_{\delta,\varepsilon}}{2 \m_R \cos\frac{\sqrt{\m_R}\pi }{2}}\left[\cos\frac{\sqrt{\m_R}\pi u}{2R} - \cos\frac{\sqrt{\m_R}\pi }{2}\right] = \varphi(u) . $$
\end{small}

\subsection{Optimal test Function} \label{ss:optimaltestfunction}
In this section, we solve the Volterra equation with temporal shifts which appears in lemma \ref{lem:equationvolterra}. Precisely, we prove that, except for a countable or finite number of value of $R$, the previous Volterra equation admits one and only one solution in $H_R$. We are giving an explicit expression of $\h_R$.

\subsubsection*{An appropriate partition} For technical convenience, since $R$ will tend to $\nu$ (thanks to lemma \ref{lem:densityarg}), we may assume that $2R$ is not an integer. Similarly, since $n \geq 1$ is the only integer such that $n-1 < 2\nu \leq n$, we may assume $\frac{n-1}{2}< R < \frac{n}{2}$, then $n =\lfloor 2R\rfloor +1$. \\
Let $D$ be the derivative operator and $T$, $T^{-1}$ the shift operators defined by:
\begin{small}$$D[f]=f' ,\hspace*{3mm} T[f](u)= f(u+1)\hspace{2mm}\mbox{ and  }\hspace{3mm} T^{-1}[f](u)= f(u-1). $$\end{small}We build a partition of $[-R,R]$ which is invariant by $T[id]$, $T^{-1}[id]$ and by symmetry. We define
$$ \left\{ \begin{array}{lcl} a_{n-2i}=R-i & \mbox{if} & 0\leq i \leq \frac{n-1}{2}  \\ a_{n-2i-1}= \lfloor 2R \rfloor -R -i & \mbox{if} & 0\leq i \leq \frac{n-2}{2} \end{array} \right.  $$
and $I_0=]-a_1,a_1[$, $I_k=]a_k,a_{k+1}[$, $I_{-k}=]-a_{k+1},-a_k[$  for $1\leq k \leq n-1$. We have \begin{small}$$[-R,R]= \bigcup_{k=-\lfloor 2R\rfloor}^{\lfloor 2R\rfloor}\overline{I_k}.$$\end{small}Moreover, if $|k|<n-2$ then $T^{\pm 1}(I_k)=I_{k\pm 2}$ and, if $k\in \{n-2,n-1\}$ then $T(I_k) \cap \mbox{ supp }\h_R = \emptyset$.

\subsubsection*{A differential equation with temporal shifts}
By derivation and integration, we remark that $\h_R$ satisfies the previous Volterra equation is equivalent to $\h_R$ satisfies the following differential equation with temporal shifts: \\for all $u$ in $]-R, R[$,
\begin{eqnarray}
\h_R '(u)  = \varphi'(u) -\frac{\delta}{2}\left[\h_R (u+1)-\h_R (u-1) \right] \label{eqn:eddecalale} .
\end{eqnarray} 
Let 
$$\lambda = \frac{\sqrt{\m_R}\pi }{2R}  \hspace{5mm}\mbox{ and }\hspace{9mm} w= \frac{-\pi^2 k_{\delta,\varepsilon} }{4R\sqrt{\m_R}\cos \frac{\pi \sqrt{\m_R}}{2}}  . $$
We consider polynomial sequences $(T_n)$ and $(U_n)$ which are defined by induction:
$$\left\{ \begin{array}{l} T_0=1 \\ T_1=X \\T_{n+1}=2XT_n-T_{n-1} \end{array} \right. \hspace{2cm} \left\{ \begin{array}{l} U_0=1 \\ U_1=2X \\U_{n+1}=2XU_n-U_{n-1} \end{array} \right. $$
Polynomials $T_n$ (resp. $U_n$) are called Chebychev polynomials of first kind (resp. second kind). They satisfy, for all real number $\theta$ and all non-negative integer $n$,
$$ T_n(\cos \theta) = \cos n\theta\hspace{3mm} \mbox{and } \hspace{2mm} U_n(\cos \theta)=\frac{\sin (n+1)\theta}{\sin \theta} . $$
As a result, we have
$$T_n(X) =\prod_{k=0}^{n-1} \left[ X- \cos\left(\frac{\pi}{2n}+\frac{k\pi}{n}  \right)\right], \hspace{2mm} U_n(X) = \prod_{k=1}^{n}\left[ X- \theta_n(k)\right]  \mbox{ where } \theta_n(k)=\cos \left(\frac{k\pi}{n+1}  \right) . $$
\paragraph*{}
In this subsection, we prove the following proposition:
\begin{proposition} \label{pro:solutionEDnoncontinue} Even functions $h$ which are compactly supported on $[-R,R]$, $\mathcal{C}^1$ on each $I_k$ and which satisfy the differential equation $(\ref{eqn:eddecalale})$ on $\cup I_k$ are those that satisfy:
\begin{itemize} \item  If $0\leq k \leq n-1$, then:
\begin{small}
\begin{multline*}
h|_{I_{n-(2k+1)}}(u) = \sum_{j=1}^{\lfloor \frac{n+1}{2}\rfloor}r_j(I_{n-1})U_k\left(\theta_n(j)\right)  \sin \left( \left[u-\frac{n-2k-1}{2}\right]\theta_{n}(j)-\frac{\pi}{2}\left[ j+\delta \frac{n-2k-1}{2} \right]\right)\\
+ r_n^\lambda(k) \sin \left(  \lambda \left[ u-\frac{n-2k-1}{2} \right]+\theta_n^\lambda(k) \right)
\end{multline*}
\end{small}
\item 
If $0\leq k \leq n-2$, then:
\end{itemize}
\begin{small}
\begin{multline*}
h|_{I_{n-2(k+1)}}(u) = \sum_{j=1}^{\lfloor \frac{n}{2}\rfloor}r_j(I_{n-2})U_k\left(\theta_{n-1}(j)\right)  \sin \left( \left[u-\frac{n-2(k+1)}{2}\right]\theta_{n-1}(j)-\frac{\pi}{2}\left[ j+\delta \frac{n-2(k+1)}{2} \right]\right)\\
+ r_{n-1}^\lambda(k)\sin \left(  \lambda \left[ u-\frac{n-2k-2}{2} \right]+\theta_{n-1}^\lambda(k) \right)
\end{multline*}
\end{small}where $r_j(I_{n-1})$ with $1\leq j \leq \left\lfloor \frac{n+1}{2} \right\rfloor$ and  $r_j(I_{n-2})$ with $1\leq j \leq \left\lfloor \frac{n}{2} \right\rfloor$ are arbitrary real numbers and where $r_n^\lambda(k)$ and $\theta_n^\lambda(k)$ refer respectively to the modulus and the argument of the complex number 
$$ \frac{i w}{ \lambda +\delta \sin \lambda}  \left[\frac{U_k \left(\lambda \right) }{U_n \left(\lambda\right)} e^{-i\left(\lambda+\frac{\pi \delta}{2}\right)\frac{n+1}{2}} -e^{i\left(\lambda+\frac{\pi \delta}{2}\right) \frac{n-2k-1}{2}}+ \frac{U_{n-k-1} (\lambda)}{U_n \left(\lambda \right)}e^{i\left(\lambda+\frac{\pi \delta}{2}\right)\frac{n+1}{2}} \right] e^{-i\frac{\pi \delta}{2}\frac{n-2k-1}{2}}. $$
\end{proposition}
Before proving this result, we need to prove some technical lemmas.

\paragraph*{Transformation in a linear differential equation}
With corollary \ref{cor:remasimplificatrice}, we may assume $\delta\neq 0$ (hence $\delta^2=1$) and $n\geq 2$. Moreover, throughout this section $h$ refers to an even function which is compactly supported on $[-R,R]$, $\mathcal{C}^1$ on each $I_k$ and which satisfies the differential equation $(\ref{eqn:eddecalale})$ on $\cup I_k$.\\
Let $Q_k$ the operator which has the following  recursive definition:
$$\left\{ \begin{array}{l} Q_1=I \\ Q_2=D+\frac{\delta}{2}(T^{-1}-T) \\ Q_{k+1}=DQ_k+ \frac{1}{4}Q_{k-1}+\left(\frac{\delta}{2}\right)^k\left(T^{-k}+(-1)^kT^{k}\right) \end{array} \right. $$
\begin{lemma} \label{lem:EDAVEC2} $h$ satisfies the following linear differential equations:
\begin{itemize}
\item On $I_{n-1}$, if $1\leq k \leq n$,  $h$ satisfies:
$$\left(\frac{i}{2}\right)^kU_k\left(\frac{1}{i}D\right)[h] = Q_k[\varphi']+\left(\frac{\delta}{2}\right)^k T^{-k}[h]   $$
\item On $I_{n-2}$, if $1\leq k \leq n-1$,  $h$ satisfies:
$$\left(\frac{i}{2}\right)^kU_k\left(\frac{1}{i}D\right)[h] = Q_k[\varphi']+\left(\frac{\delta}{2}\right)^k T^{-k}[h]   $$
\item On $I_0$, if $1\leq k \leq \lfloor \frac{n+1}{2} \rfloor $, $h$ satisfies:
$$2\left(\frac{i}{2}\right)^kT_k\left(\frac{1}{i}D\right)[h]    = Q_k[\varphi']+\left(\frac{\delta}{2}\right)^k \left( T^{-k}[h]+(-1)^kT^{k}[h ] \right)   $$
\end{itemize}
\end{lemma}
Proof: We prove the first relation by induction for $k=1,...,n$.  
If $k=1$,the differential equation (\ref{eqn:eddecalale}) gives:
\begin{eqnarray}
\mbox{if }u\in I_{n-1}, \hspace{2mm} h'(u)=\varphi'(u)+\frac{\delta}{2}h(u-1). \label{eqn:rec1}
\end{eqnarray}
The result comes  easily from the definition of $Q_1$ and $U_1$.
If $k=2$, we derive relation (\ref{eqn:rec1}) and we apply the differential equation (\ref{eqn:eddecalale}).\\
We assume the result holds for $k$ and $k-1$ with $2\leq k\leq n-1$. Then, with recursive relation on $U_k$, we get:
\begin{small}\begin{eqnarray*}
& & \left(\frac{i}{2}\right)^{k+1}U_{k+1}\left( \frac{1}{i}D\right)[h] = \frac{2}{i}\left(\frac{i}{2}\right)^{k+1}DU_{k}\left( \frac{1}{i}D\right)[h]-\left(\frac{i}{2}\right)^{k+1}U_{k-1}\left( \frac{1}{i}D\right)[h] \\
& = & D\left( \left(\frac{i}{2}\right)^{k}U_{k}\left( \frac{1}{i}D\right)[h] \right)+\frac{1}{4}\left(\frac{i}{2}\right)^{k-1}U_{k-1}\left( \frac{1}{i}D\right)[h] \\
 & = & DQ_k[\varphi']+\frac{1}{4}Q_{k-1}[\varphi']+\left(\frac{\delta}{2}\right)^kT^{-k}[h']+\left(\frac{\delta}{2}\right)^{k+1}T^{-(k-1)}[h]\\
 & & \hspace{32mm}\mbox{since } k \leq n-1 \mbox{ we may apply differential equation } (\ref{eqn:eddecalale}). \\
 & = & DQ_k[\varphi']+\frac{1}{4}Q_{k-1}[\varphi']+\left(\frac{\delta}{2}\right)^kT^{-k}[\varphi']+\left(\frac{\delta}{2}\right)^{k+1}T^{-(k+1)}[h]
\end{eqnarray*}\end{small}Since we apply this relation with $u$ in $I_{n-1}$, we have $T^{k+1}[h](u)=0$. Therefore, we have
\begin{small}$$  \left(\frac{i}{2}\right)^{k+1}U_{k+1}\left( \frac{1}{i}D\right)[h] = Q_{k+1}[\varphi']+\left(\frac{\delta}{2}\right)^{k+1}T^{-(k+1)}[h]\hspace{2mm}\mbox{ on }I_{n-1}.$$ \end{small}We prove the other relations of this lemma in the same way.
\begin{flushright}
$\Box$
\end{flushright}
\begin{lemma} \label{lem:expressionQ} Let $f$ be a function which is $\mathcal{C}^1$ on $]-R,R[$ and with supp $f \subset [-R,R]$. Then, for all $1\leq k \leq n$  we have
$$ Q_k[f|_{I_{n-1}}]=\left( \frac{i}{2} \right)^{k-1}\sum_{j=0}^{k-1}\left( \frac{\delta}{i} \right)^j U_{k-j-1}\left(\frac{1}{i}D\right)T^{-j}[f|_{I_{n-1}}] . $$
Similarly, for all $1\leq k \leq n-1$  we have
$$ Q_k[f|_{I_{n-2}}]=\left( \frac{i}{2} \right)^{k-1}\sum_{j=0}^{k-1}\left( \frac{\delta}{i} \right)^j U_{k-j-1}\left(\frac{1}{i}D\right)T^{-j}[f|_{I_{n-2}}] . $$
\end{lemma} 
Proof: By induction on $k=1,...,n$. The result comes easily from  the definition of $Q_k$ if $k=1$ or $k=2$. We assume the result holds for $k$ and $k-1$ with $2\leq k\leq n-1$. Since $T^{k}[f|_{I_{n-1}}]=0$, we have
\begin{small}
\begin{eqnarray*}
 & & Q_{k+1} = DQ_k+\frac{1}{4}Q_{k-1} + \left( \frac{\delta}{2} \right)^k T^{-k} \\
& &  = \left( \frac{i}{2} \right)^{k}\left( \sum_{j=0}^{k-1}\left( \frac{\delta}{i} \right)^j \frac{2D}{i}U_{k-j-1}\left(\frac{1}{i}D\right)T^{-j}          -\sum_{j=0}^{k-2}\left( \frac{\delta}{i} \right)^j U_{k-j-2}\left(\frac{1}{i}D\right)T^{-j}\right)+ \left( \frac{\delta}{2} \right)^k T^{-k} \\
 & &= \left( \frac{i}{2} \right)^{k}\sum_{j=0}^{k-2}\left( \frac{\delta}{i} \right)^j \left[ \frac{2D}{i}U_{k-j-1}\left(\frac{1}{i}D\right)-U_{k-j-2}\left(\frac{1}{i}D\right)\right]T^{-j} +\left(\frac{i}{2} \right)^{k}\left( \frac{\delta}{i}\right)^{k-1} U_{1}\left(\frac{1}{i}D\right)T^{-(k-1)}
 \\ & & \hspace{130mm} +\left( \frac{\delta}{2}\right)^{k} T^{-k}  \\
 & & = \left(\frac{i}{2} \right)^{k} \sum_{j=0}^{k}\left( \frac{\delta}{i} \right)^j U_{k-j}\left(\frac{1}{i}D\right)T^{-j} . 
\end{eqnarray*}
\end{small}Similarly, we prove the other relation of this lemma.
\begin{flushright}
$\Box$
\end{flushright}
Using lemmas \ref{lem:EDAVEC2} and \ref{lem:expressionQ}, we prove that $h|_{I_{n-1}}$ and $h|_{I_{n-2}}$ satisfy linear non-homogeneous differential equations with constant coefficients. Precisely,
\begin{corollary} \label{cor:eqdlin} We have
$$i^nU_n\left(\frac{1}{i}D\right)[h|_{I_{n-1}}] = \emph{Im }\left[iw\delta^n e^{i\lambda(u-n)}\frac{1-(i\delta e^{i\lambda})^n U_n(\lambda)+(i\delta e^{i\lambda})^{n+1}U_{n-1}(\lambda)}{\lambda+\delta\sin  \lambda}\right] $$
and 
$$i^{n-1}U_{n-1}\left(\frac{1}{i}D\right)[h|_{I_{n-2}}] = \emph{Im }\left[iw\delta^{n-1} e^{i\lambda(u-n+1)}\frac{1-(i\delta e^{i\lambda})^{n-1} U_{n-1}(\lambda)+(i\delta e^{i\lambda})^n U_{n-2}(\lambda)}{\lambda+\delta\sin  \lambda} \right] .$$
\end{corollary}
Proof: We have $T^{-n}[h|_{I_{n-1}}]=0$ and $T^{-(n-1)}[h|_{I_{n-2}}]=0$. Thus, with lemma \ref{lem:EDAVEC2}, we get
\begin{small}$$i^nU_n\left(\frac{1}{i}D\right)[h|_{I_{n-1}}] =2^n Q_n[\varphi'|_{I_{n-1}}]\hspace{3mm} \mbox{ and } \hspace{3mm} i^{n-1}U_{n-1}\left(\frac{1}{i}D\right)[h|_{I_{n-2}}] =2^{n-1}Q_{n-1}[\varphi'|_{I_{n-2}}] . $$\end{small}Let $e_\lambda(u)= we^{i\lambda u}$. We have $\varphi'(u)=\mbox{Im }[e_\lambda(u)]$. Since  $D^2(e_\lambda)=(i\lambda)^2e_\lambda$, thanks to lemma \ref{lem:expressionQ}, we write
\begin{small}$$Q_n(e_\lambda|_{I_{n-1}}) = \left( \frac{i}{2} \right)^{n-1}\sum_{j=0}^{n-1}\left( \frac{\delta}{i} \right)^j U_{n-j-1}(\lambda)T^{-j}(e_\lambda|_{I_{n-1}}) $$\end{small}and 
\begin{small}$$Q_{n-1}(e_\lambda|_{I_{n-2}}) = \left( \frac{i}{2} \right)^{n-2}\sum_{j=0}^{n-2}\left( \frac{\delta}{i} \right)^j U_{n-j-2}(\lambda)T^{-j}(e_\lambda|_{I_{n-2}}) .$$\end{small}Due to the fact that
\begin{eqnarray} \sum_{j=0}^{n-1}U_j(X)z^j = \frac{1-z^nU_n(X)+z^{n+1}U_{n-1}(X)}{1-2zX+z^2}, \label{eqn:seriegeneratricetronquee}\end{eqnarray}
we have
\begin{small}$$ Q_n(e_\lambda|_{I_{n-1}}) =\frac{\delta^{n-1}ie^{-i\lambda}}{2^n}\frac{1-(i\delta e^{i\lambda})^n U_n(\lambda)+(i\delta e^{i\lambda})^{n+1}U_{n-1}(\lambda)}{\lambda\delta+\sin\lambda}e_\lambda(u-n+1) $$\end{small}and
\begin{small}$$ Q_n(e_\lambda|_{I_{n-2}}) =\frac{\delta^{n-2}ie^{-i\lambda}}{2^{n-1}}\frac{1-(i\delta e^{i\lambda})^{n-1} U_{n-1}(\lambda)+(i\delta e^{i\lambda})^{n}U_{n-2}(\lambda)}{\lambda\delta+\sin\lambda}e_\lambda(u-n+2) .$$\end{small}We may conclude since $Q_n(\varphi')=\mbox{Im} (Q_n(e_\lambda))$.
\begin{flushright}
$\Box$
\end{flushright}

\begin{lemma}  $R\mapsto \tilde{\m}_R$ is a strictly decreasing function.
\end{lemma}
Proof: We are assuming that there exist $R_1<R_2$ such that $\tilde{\m}_{R_1}= \tilde{\m}_{R_2}$. Since $R\mapsto \tilde{\m}_R$ is a decreasing function, we may assume there exists an integer $n_0$ such that
\begin{small}
$$\frac{n_0-1}{2} <R_1<R_2 < \frac{n_0}{2}.$$
\end{small}Let $\lambda = \lambda_{R_1}=\lambda_{R_2} =2\pi \sqrt{\tilde{\m}_{R_1}}$. 
There exists $\h_{R_1}$ in $H_{R_1}$ such that $\tilde{B}(\h_{R_1})= \tilde{m}_{R_1}= \tilde{m}_{R_2}$. Thanks to the Volterra equation (\ref{eqn:reftardive}), $\h_{R_1}|_{I_{n_0-1}(R_2)}$ is smooth. Using corollary \ref{cor:eqdlin}, we have:
\begin{small} $$
i^{n_0}U_{n_0}\left(\frac{1}{i}D\right)[\h_{R_1}|_{I_{n_0-1}(R_2)}](u) = \mbox{Im }\left[iw_{R_2}\delta^{n_0} e^{i\lambda(u-n_0)}\frac{1-(i\delta e^{i\lambda})^{n_0} U_{n_0}(\lambda)+(i\delta e^{i\lambda})^{n_0+1}U_{n_0-1}(\lambda)}{\lambda+\delta\sin  \lambda}\right] $$ \end{small}Since supp $\h_{R_1} \subset [-R_1,R_1]$, we have $\h_{R_1}|_{I_{n_0-1}(R_2) \backslash I_{n_0-1}(R_1) } =0$. Thus, for all $u$ in $I_{n_0-1}(R_2)$, we have \begin{small}
$$\mbox{Im }\left[iw_{R_2}\delta^n e^{i\lambda(u-n_0)}\frac{1-(i\delta e^{i\lambda})^{n_0} U_{n_0}(\lambda)+(i\delta e^{i\lambda})^{n_0+1}U_{n_0-1}(\lambda)}{\lambda+\delta\sin  \lambda}\right] =0. $$
\end{small}Then, on $I_{n_0-1}(R_2)$, $\h_{R_1}|_{I_{n_0-1}(R_2)}$ satisfies $ U_{n_0}\left(\frac{1}{i}D\right)[\h_{R_1}|_{I_{n_0-1}(R_2)}](u) =0$ and $\h_{R_1}(u)=0$ on $I_{n_0-1}(R_2) \backslash I_{n_0-1}(R_1) \neq \emptyset$. Using Picard-Lindelöf theorem, we may conclude $\h_{R_1}|_{I_{n_0-1}(R_2)}=0$. Thus, supp $\h_{R_1} \subset [-a_{n_0-1}(R_2),a_{n_0-1}(R_2) ]$ and $\tilde{\m}_{a_{n_0-1}(R_2)} = \tilde{\m}_{R_2}$ with $a_{n_0-1}(R_2)=\lfloor 2R_2\rfloor -R_2$. As a result, for all $R$ such that $a_{n_0-1}(R_2)<R<R_2$, we have $\tilde{m}_R=\tilde{m}_{R_2}$.
By induction, there exists $r_k$ (with $r_1=a_{n-1}(R_2)$) in $]\frac{n_0-k-1}{2},\frac{n_0-k}{2}[$ such that for all $R$ in $[r_k,R_2[$, we have
$\tilde{m}_R=\tilde{m}_{R_2}$. As a result, $R\mapsto \tilde{m}_R$ is a constant function on $]r_{n_0-1},1/2[$ which contradicts corollary \ref{cor:remasimplificatrice}.
\begin{flushright}
$\Box$
\end{flushright}
\begin{corollary} \label{cor:numnonnul} Except for at most $n-1$ values of $R$ in $]\frac{n-1}{2},\frac{n}{2}[$, we have $U_n(\lambda)U_{n-1}(\lambda)\neq 0$.
\end{corollary}
Proof: $U_n(\lambda)U_{n-1}(\lambda)= 0$ if $\lambda$ is one of the $n-1$ positive roots of $U_nU_{n-1}$. Since $R\mapsto \lambda =2\pi \sqrt{\tilde{m}_R}$ is a strictly decreasing function, $U_n(\lambda)U_{n-1}(\lambda)\neq 0$ excepted for at most $n-1$ values of $R$ in $]\frac{n-1}{2},\frac{n}{2}[$.
\begin{flushright}
$\Box$
\end{flushright}
Since $R$ will tends to $\nu$, we may assume $U_n(\lambda)U_{n-1}(\lambda)\neq 0$.
\begin{lemma} \label{lem:SEDComplex1} There exists $2n-1$ complex numbers \begin{small}$z_1(I_{n-1}),...,z_{n}(I_{n-1})$\end{small} and \begin{small}$z_1(I_{n-2}),...,z_{n-1}(I_{n-2})$\end{small}such that
 \begin{eqnarray*}& &h|_{I_{n-1}}(u) = \emph{Im} \left[\sum_{j=1}^{n}z_j(I_{n-1})e^{iu\theta_n(j)}+ z_n^\lambda(0)e^{i\lambda\left(u-\frac{n-1}{2}\right)}  \right] \\
& \mbox{and }&   h|_{I_{n-2}}(u) = \emph{Im} \left[\sum_{j=1}^{n-1}z_j(I_{n-2})e^{iu\theta_{n-1}(j)}+z_{n-1}^\lambda(0)e^{i\lambda\left(u-\frac{n-2}{2}\right)}   \right] \\
\end{eqnarray*}
where
$$ z_n^\lambda(0)=r_n^\lambda(0)e^{i\theta_n^\lambda(0)}= \frac{iw}{(i\delta)^nU_n(\lambda)} e^{-i\lambda \frac{n+1}{2}}\frac{1-(i\delta e^{i\lambda})^n U_n(\lambda)+(i\delta e^{i\lambda})^{n+1}U_{n-1}(\lambda)}{\lambda+\delta\sin  \lambda} . $$
\end{lemma}
Proof: We solve differential equations of corollary \ref{cor:eqdlin}. The general solution of the homogeneous equation is any functions which may be written \begin{small} $$ u\longmapsto \mbox{Im } \left[\sum_{j=1}^{n}z_j(I_{n-1})e^{iu\theta_n(j)} \right] . $$ \end{small}Finally, with corollary \ref{cor:numnonnul}, it is easy to check that $ u\longmapsto\mbox{Im }\left[ z_n^\lambda(0)e^{i\lambda\left(u-\frac{n-1}{2}\right)}  \right]$ is a particular solution of our differential equation. Similarly, we success in obtaining an explicit expression of $h|_{I_{n-2}}$.
\begin{flushright}
$\Box$
\end{flushright}

\paragraph*{Extension of the optimal test function}
In lemma \ref{lem:SEDComplex1}, we have an explicit expression of $h$ only on $I_{n-2}\cup I_{n-1}$. Thanks to the differential equation with temporal shifts (\ref{eqn:eddecalale}), we extend this explicit expression to $[-R,R]$. 

\begin{lemma}\label{lem:SEDComplex2} With $U_{-1}=0$, if $k=0,...,n-1$ then:
$$ h|_{I_{n-(2k+1)}} = \left(\frac{i}{\delta} \right)^{k} T^{k} U_k \left(\frac{1}{i} D \right)[h|_{I_{n-1}}]- \frac{2}{\delta} \sum_{j=0}^{k} \left(\frac{i}{\delta} \right)^{j-1}U_{j-1}\left( \frac{1}{i}D\right)T^j[\varphi']  $$
Similarly, if $k=0,...,n-2$, then:
$$ h|_{I_{n-2(k+1)}} = \left(\frac{i}{\delta} \right)^{k} T^{k} U_k \left(\frac{1}{i} D \right)[h|_{I_{n-2}}]-\frac{2}{\delta} \sum_{j=0}^{k} \left(\frac{i}{\delta} \right)^{j-1}U_{j-1}\left( \frac{1}{i}D\right)T^j[\varphi']  $$
\end{lemma}
Proof: We prove the first relation by induction on $k=0,..,n-1$. If $k=0$, there is nothing to prove. If $k=1$, the result comes easily from relation (\ref{eqn:eddecalale}). Assume the result holds for $k-1$ and $k$ with $1\leq k \leq n-2$. Using relation (\ref{eqn:eddecalale}), we get
\begin{small}$$ h'|_{I_{n-(2k+1)}}(u+1)= \varphi'(u+1)+\frac{\delta}{2}\left( h|_{I_{n-(2k+3)}}(u)- h|_{I_{n-(2k-1)}}(u+2) \right). $$ \end{small}Therefore, we have
\begin{small}\begin{eqnarray*}
 h|_{I_{n-(2k+3)}} =T^2[ h|_{I_{n-(2k-1)}}]+\frac{2}{\delta}\left(TD[h|_{I_{n-(2k+1)}}] -T[\varphi']\right).
\end{eqnarray*}\end{small}The result follows easily from the induction hypothesis.
\begin{flushright}
$\Box$
\end{flushright}

\begin{lemma} \label{lem:expressionparitedeg} If $0\leq k \leq n-1$, then
$$ h|_{I_{n-(2k+1)}}(u) = \emph{Im} \left[(i\delta)^k \sum_{j=1}^{n}z_j(I_{n-1}) U_k(\theta_n(j))e^{i(u+k)\theta_n(j)}+z_n^\lambda (k)e^{i\lambda \left(u-\frac{n-2k-1}{2}\right)} \right]$$
and if $0\leq k \leq n-2$, then
$$h|_{I_{n-2(k+1)}}(u) = \emph{Im} \left[(i\delta)^k \sum_{j=1}^{n-1}z_j(I_{n-2}) U_k(\theta_{n-1}(j))e^{i(u+k)\theta_{n-1}(j)}+z_{n-1}^\lambda (k)e^{i\lambda \left(u-\frac{n-2k-2}{2}\right)} \right]$$
where
\begin{eqnarray} \label{eqn:defznlambdak} z_n^\lambda(k) & = & r_n^\lambda(k)e^{i\theta_n^\lambda(k)} \\
& = &   \frac{iw}{\lambda +\delta \sin \lambda }  \left[  (i\delta)^{k-n}e^{-i\lambda\frac{n+1}{2}}\frac{U_k(\lambda)}{U_n(\lambda)}-e^{i\lambda\frac{n-2k-1}{2}}  +(i\delta)^{k+1}e^{i\lambda\frac{n+1}{2}} \frac{U_{n-k-1}(\lambda)}{U_n(\lambda) } \right] . \nonumber
\end{eqnarray}
\end{lemma}
Proof: With lemmas \ref{lem:SEDComplex1} and \ref{lem:SEDComplex2}, we may write
\begin{small}
\begin{eqnarray*}
h|_{I_{n-(2k+1)}}(u) = \mbox{Im } \left[ \sum_{j=1}^{n}z_j(I_{n-1})\left(\frac{i}{\delta} \right)^k U_k(\theta_n(j))e^{i(u+k)\theta_n(j)}+z_n^\lambda(0)\left(\frac{i}{\delta} \right)^k U_k(\lambda)e^{i\lambda\left(u-\frac{n-2k-1}{2}\right)} \right. \\ \left.- \frac{2w}{\delta} \sum_{j=0}^{k} \left(\frac{i}{\delta} \right)^{j-1}U_{j-1}\left( \frac{1}{i}D\right)e^{i\lambda (u+j)}\right] .
\end{eqnarray*}
\end{small}Thanks to relation  (\ref{eqn:seriegeneratricetronquee}), we get
\begin{small}\begin{eqnarray*}
\frac{2w}{\delta} \sum_{j=0}^{k} \left(\frac{i}{\delta} \right)^{j-1}U_{j-1}\left( \frac{1}{i}D\right)e^{i\lambda (u+j)}&  = & 2w\delta e^{i\lambda(u+1)}\sum_{j=0}^{k-1} \left(\frac{i}{\delta} \right)^{j}U_j(\lambda)e^{i\lambda j}\\
 & = & iwe^{i\lambda u}\frac{1-(i\delta e^{i\lambda})^k U_k(\lambda)+(i\delta e^{i\lambda})^{k+1} U_{k-1}(\lambda) }{\lambda+\delta \sin \lambda} .
\end{eqnarray*}\end{small}Thus  we have:
\begin{small}
\begin{eqnarray*}
& & z_n^\lambda(0)\left(\frac{i}{\delta} \right)^k U_k(\lambda)e^{i\lambda\left(u-\frac{n-2k-1}{2}\right)} - \frac{2w}{\delta} \sum_{j=0}^{k} \left(\frac{i}{\delta} \right)^{j-1}U_{j-1}\left( \frac{1}{i}D\right)e^{i\lambda (u+j)} \\
& & \hspace{8mm} = \frac{iw}{\lambda+\delta \sin \lambda}e^{i\lambda \left(u-\frac{n-2k-1}{2} \right)}\left((i\delta)^{k-n}e^{-i\lambda\frac{n+1}{2}}\frac{U_k(\lambda)}{U_n(\lambda)}+ (i\delta)^{k+1}e^{i\lambda\frac{n+1}{2}}\frac{U_k(\lambda)U_{n-1}(\lambda)}{U_n(\lambda)} -e^{i\lambda\frac{n-2k-1}{2}} \right.\\
& & \hspace{113mm} \left. - (i\delta)^{k+1}e^{i\lambda\frac{n+1}{2}}U_{k-1}(\lambda)\right)
\end{eqnarray*}
\end{small}We may conclude since, if  $n\geq 2$ and $1\leq j \leq n-1$, then
\begin{eqnarray} U_{n-1}U_j-U_nU_{j-1}=U_{n-1-j} .\label{eqn:relationmultiplicativitedesfibo} \end{eqnarray}
\begin{flushright}
$\Box$
\end{flushright}

\paragraph*{Even conditions} $h$ is assumed to be even. We are exploiting this fact in order to obtain some restrictions on complex numbers $z_j$.
\begin{lemma} $h$ is given lemma \ref{lem:expressionparitedeg}. Then $h$ is an even function if and only if we have:\\
For  $1\leq j \leq \left\lfloor \frac{n+1}{2} \right\rfloor$, 
$$ \overline{z_j(I_{n-1})}-z_{n+1-j}(I_{n-1})+(i\delta)^{n-1}U_{n-1}(\theta_n(j))e^{i(n-1)\theta_n(j)}\left( z_j(I_{n-1})-\overline{z_{n+1-j}(I_{n-1})} \right)=0 ,$$
and, for $1\leq j \leq \left\lfloor \frac{n}{2} \right\rfloor$, 
$$\overline{z_j(I_{n-2})}-z_{n-j}(I_{n-2})+(i\delta)^{n-2}U_{n-2}(\theta_{n-1}(j))e^{i(n-2)\theta_{n-1}(j)}\left( z_j(I_{n-2})-\overline{z_{n-j}(I_{n-2})} \right)=0 .$$
\end{lemma}
Proof: First, using relation (\ref{eqn:eddecalale}), we prove that $h$ is even if and only if $h(u)=h(-u)$ for all $u\in I_{n-2}\cup I_{n-1}$. Then, due to the fact that  $\overline{z_n^\lambda(n-1)} =- z_n^\lambda(0)$, if $u\in I_{n-1}$, we have:
\begin{small}\begin{eqnarray*}
& & h|_{I_{n-1}}(u)=h_{I_{-(n-1)}}(-u) \\
& & \Longleftrightarrow \mbox{Im } \left[ \sum_{j=1}^{n} z_j(I_{n-1}) e^{iu\theta_n(j)} \right]= \mbox{Im } \left[ (i\delta)^{n-1}\sum_{j=1}^{n} z_j(I_{n-1}) U_{n-1}(\theta_n(j))e^{i(n-1)\theta_n(j)}e^{-iu\theta_n(j)} \right]
\end{eqnarray*} \end{small}Since $\theta_n(n+1-k)=-\theta_n(k)$, we get
\begin{small}$$\mbox{Im } \left[ \sum_{j=1}^{n} z_j(I_{n-1}) e^{iu\theta_n(j)} \right] = \mbox{Im }\left[ \sum_{j=1}^{\lfloor \frac{n+1}{2} \rfloor} \left( z_j(I_{n-1})-\overline{z_{n+1-j}(I_{n-1}} \right) e^{iu\theta_n(j)} \right]$$\end{small}and 
\begin{small}\begin{eqnarray*}
& &\mbox{Im } \left[ (i\delta)^{n-1}\sum_{j=1}^{n} z_j(I_{n-1}) U_{n-1}(\theta_n(j))e^{i(n-1)\theta_n(j)}e^{-iu\theta_n(j)} \right]  \\
& & \hspace*{7mm} = - \mbox{Im } \left[ (-i\delta)^{n-1}\sum_{j=1}^{\lfloor \frac{n+1}{2} \rfloor} \left(\overline{z_j(I_{n-1})}-z_{n+1-j}(I_{n-1}) \right) U_{n-1}(\theta_n(j))e^{-i(n-1)\theta_n(j)}e^{iu\theta_n(j)} \right] .
\end{eqnarray*} \end{small}Therefore, for all $1\leq j \leq \left\lfloor \frac{n+1}{2} \right\rfloor$, we get
\begin{small}$$ \overline{z_j(I_{n-1})}-z_{n+1-j}(I_{n-1})+(i\delta)^{n-1}U_{n-1}(\theta_n(j))e^{i(n-1)\theta_n(j)}\left( z_j(I_{n-1})-\overline{z_{n+1-j}(I_{n-1})} \right)=0 . $$\end{small}
\begin{flushright}
$\Box$
\end{flushright}

\begin{corollary} \label{cor:condpar} For all  $1\leq j \leq \lfloor \frac{n+1}{2} \rfloor$, there exists a real number $r_j(I_{n-1})$ (which may be negative) such that
$$\overline{z_j(I_{n-1})}-z_{n+1-j}(I_{n-1}) = r_j(I_{n-1})e^{i\left[\frac{n-1}{2}\left(\theta_n(j)+\frac{\pi \delta}{2} \right)+\frac{j\pi}{2}\right]}.$$
For all  $1\leq j \leq \lfloor \frac{n}{2} \rfloor$, there exists a real number $r_j(I_{n-2})$ (which may be negative) such that
$$\overline{z_j(I_{n-2})}-z_{n-j}(I_{n-2}) = r_j(I_{n-2})e^{i\left[\frac{n-2}{2}\left(\theta_{n-1}(j)+\frac{\pi \delta}{2} \right)+\frac{j\pi}{2}\right]} .$$ 
\end{corollary}
Proof: We have 
\begin{small}$$U_{n-1}(\theta_n(j))=U_{n-1}\left( \cos \left( \frac{j \pi}{n+1}\right)\right) =(-1)^{j+1} . $$\end{small}Therefore, if $\overline{z_j(I_{n-1})}-z_{n+1-j}(I_{n-1})= r_j e^{i\theta_j}$, the previous lemma gives
$$ r_j \left( e^{i\theta_j}+(i\delta)^{n-1}(-1)^{j-1}e^{i(n-1)\theta_j}e^{-i\theta_j} \right)=0. $$
Then $r_j=0$ or $e^{i\theta_j} = \pm e^{i\left[ \frac{n-1}{2}\left(\theta_n(j)+\frac{\pi \delta}{2} \right)+\frac{j\pi}{2} \right]}$.
\begin{flushright}
$\Box$
\end{flushright}
\paragraph*{Proof of proposition \ref{pro:solutionEDnoncontinue}}
From lemma \ref{lem:expressionparitedeg}, we have for $0\leq k \leq n-1$:
\begin{small}
$$ h|_{I_{n-(2k+1)}}(u) = \mbox{Im} \left[(i\delta)^k \sum_{j=1}^{\lfloor \frac{n+1}{2} \rfloor }  \left( z_j(I_{n-1})-\overline{z_{n+1-j}(I_{n-1})}\right) U_k(\theta_n(j))e^{i(u+k)\theta_n(j)}+z_n^\lambda (k)e^{i\lambda \left(u-\frac{n-2k-1}{2}\right)} \right]$$
\end{small}Then, corollary \ref{cor:condpar} gives:
\begin{small}\begin{eqnarray*} h|_{I_{n-(2k+1)}}(u)& =  & \mbox{Im} \left[\sum_{j=1}^{\lfloor \frac{n+1}{2} \rfloor } r_j(I_{n-1})U_k(\theta_n(j))e^{i\left[\left(u-\frac{n-2k-1}{2}  \right)\theta_n(j)-\frac{\pi}{2}\left(j+\delta \frac{n-2k-1}{2} \right)    \right]} \right. \\ & &  \hspace{79mm} \left. +r_n^\lambda(k)e^{i\left[\lambda\left( u-\frac{n-2k-1}{2}\right)+\theta_n^\lambda (k)  \right]} \right] 
\end{eqnarray*}\end{small}We obtain $h|_{I_{n-(2k+2)}}$ in the same way. To conclude, one can easily check that such a function satisfies the differential equation with temporal shifts (\ref{eqn:eddecalale}).
\begin{flushright}
$\Box$
\end{flushright}

\subsubsection*{Continuity conditions}
Since $\h_R$ satisfies the differential equation (\ref{eqn:eddecalale}), we have an explicit expression of $\h_R$ thanks to proposition \ref{pro:solutionEDnoncontinue}. However, due to the fact that $\h_R$ is a continuous function, there are some restrictions on coefficients $r_j$ which are associated to the explicit expression of $\h_R$. In this subsection, we prove that, except for at most a finite number of $R$ in $]\frac{n-1}{2},\frac{n}{2}[$, there is one and only one continuous function which satisfies the differential equation (\ref{eqn:eddecalale}) and we give the exact values of  the corresponding coefficients $r_j$. 

\begin{lemma} \label{lem:defzlambda} $\h_R$ is a continuous function if and only if coefficients $r_j(I_{n-1})$ and $r_j(I_{n-2})$ are satisfying the following linear system: for $0\leq k \leq n-1$, 
\begin{small}
\begin{eqnarray*}
& &\sum_{j=1}^{\lfloor \frac{n}{2}\rfloor}r_j(I_{n-2})U_k\left(\theta_{n-1}(j)\right)  \sin \left( \left[a_{n-1}-\frac{n-2}{2}\right]\theta_{n-1}(j)-\frac{\pi}{2}\left[ j+\delta \frac{n-2k-2}{2} \right]\right) \\
& & \hspace{5mm}- \sum_{j=1}^{\lfloor \frac{n+1}{2}\rfloor}r_j(I_{n-1})U_k\left(\theta_n(j)\right)  \sin \left( \left[a_{n-1}-\frac{n-1}{2}\right]\theta_{n}(j)-\frac{\pi}{2}\left[ j+\delta \frac{n-2k-1}{2} \right]\right) \\
& & \hspace{16mm} =  r_n^\lambda(k)\sin\left(\lambda\left[a_{n-1}-\frac{n-1}{2}\right]+\theta_n^\lambda (k) \right)-r_{n-1}^\lambda(k)\sin\left(\lambda\left[a_{n-1}-\frac{n-2}{2}\right]+\theta_{n-1}^\lambda (k) \right)\\
& &\hspace{108mm} = r_\lambda U_k\left( \lambda \right)\sin \left( \theta_\lambda-k\delta\frac{\pi}{2}\right)
\end{eqnarray*}
\end{small}where $r_\lambda$ and $\theta_\lambda$ refer to the modulus and the argument of the complex number:
$$z_\lambda =\frac{-2iwe^{-i\lambda a_{n-1}}}{U_n(\lambda)U_{n-1}(\lambda)}   \frac{1-(i\delta e^{i\lambda})^{n}U_n(\lambda)+(i\delta e^{i\lambda})^{n+1}U_{n-1}(\lambda)}{ -2i\delta e^{i\lambda}(\lambda+\delta \sin \lambda)} = \frac{-2iwe^{-i\lambda a_{n-1}}}{U_n(\lambda)U_{n-1}(\lambda)}\sum_{k=0}^{n-1}(i\delta e^{i\lambda})^k U_k(\lambda)  $$
\end{lemma}
Proof: We use the fact that $\h_R$ is a continuous function on $\inf I_j$ for $j=-(n-1)...n-1$. Precisely, $h_R$ is continuous on $a_{n-(2k+1)} >0 $ (ie $a_{n-1}, a_{n-3}...$) then
\begin{small}\begin{eqnarray*}  \mbox{for }0\leq k \leq \left\lfloor \frac{n-2}{2}\right\rfloor , \hspace{5mm} \lim_{u\rightarrow a_{n-(2k+1)}^{-}}\h_R|_{I_{n-2(k+1)}}(u) = \lim_{u\rightarrow a_{n-(2k+1)}^{+}}\h_R|_{I_{n-(2k+1)}}(u).
 \end{eqnarray*}\end{small}Moreover, $h_R$ is continuous on $a_{n-2(k+1)} <0$ (ie $a_{-n}, a_{-(n-2)},...$) then
 \begin{small}\begin{eqnarray*} \mbox{for }\left\lfloor \frac{n}{2}\right\rfloor \leq k \leq n-1, \hspace{5mm} \lim_{u\rightarrow a_{n-2(k+1)}^{-}}\h_R|_{I_{n-2(k+1)}}(u) = \lim_{u\rightarrow a_{n-2(k+1)}^{+}}\h_R|_{I_{n-(2k+1)}}(u).
\end{eqnarray*}\end{small}Therefore, for $0\leq k \leq n-1$, we have
\begin{small}
\begin{eqnarray*}
& &\sum_{j=1}^{\lfloor \frac{n}{2}\rfloor}r_j(I_{n-2})U_k\left(\theta_{n-1}(j)\right)  \sin \left( \left[a_{n-1}-\frac{n-2}{2}\right]\theta_{n-1}(j)-\frac{\pi}{2}\left[ j+\delta \frac{n-2k-2}{2} \right]\right) \\
& & \hspace{5mm}- \sum_{j=1}^{\lfloor \frac{n+1}{2}\rfloor}r_j(I_{n-1})U_k\left(\theta_n(j)\right)  \sin \left( \left[a_{n-1}-\frac{n-1}{2}\right]\theta_{n}(j)-\frac{\pi}{2}\left[ j+\delta \frac{n-2k-1}{2} \right]\right) \\
& & \hspace{16mm} =  r_n^\lambda(k)\sin\left(\lambda\left[a_{n-1}-\frac{n-1}{2}\right]+\theta_n^\lambda (k) \right)-r_{n-1}^\lambda(k)\sin\left(\lambda\left[a_{n-1}-\frac{n-2}{2}\right]+\theta_{n-1}^\lambda (k) \right).
\end{eqnarray*}
\end{small}Reciprocally, if $\h_R$ is continuous on  $\inf I_j$, since $\h_R$ is even and $-\inf I_j =\sup I_{-j}$, then $\h_R$ is a continuous function.
\begin{flushright}
$\Box$
\end{flushright}
We consider the real matrix $M_R$ which is associated to this linear system and $\det M_R$ refers to the determinant of $M_R$. Precisely, let $M_R=(m_{k,j})_{\tiny{\begin{array}{l} k=0,...,n-1 \\ j=1,...,n \end{array}} }$ defined by:
$$  \left\{\begin{array}{ll}  m_{k,j} =U_k\left(\theta_{n-1}(j)\right)  \sin \left( \left[a_{n-1}-\frac{n-2}{2}\right]\theta_{n-1}(j)-\frac{\pi}{2}\left[ j+\delta \frac{n-2k-2}{2} \right]\right) & \mbox{if } 1\leq j \leq \lfloor \frac{n}{2}\rfloor \\ 
 m_{k,j+\lfloor \frac{n}{2}\rfloor } =U_k\left(\theta_n(j)\right)  \sin \left( \left[a_{n-1}-\frac{n-1}{2}\right]\theta_{n}(j)-\frac{\pi}{2}\left[ j+\delta \frac{n-2k-1}{2} \right]\right) & \mbox{if } 1\leq j \leq \lfloor \frac{n+1}{2}\rfloor
\end{array} \right.$$
Therefore, the linear system of lemma \ref{lem:defzlambda} is
\begin{small}$$M_R \cdot \left( \begin{array}{c} r_1(I_{n-2}) \\ \vdots \\ r_{\lfloor \frac{n}{2}\rfloor}(I_{n-2}) \\ -r_1(I_{n-1}) \\ \vdots \\ -r_{\lfloor \frac{n+1}{2}\rfloor}(I_{n-1)}) 
\end{array}\right)  = r_{\lambda} \left( \begin{array}{c} U_0\left( \lambda \right)\sin \left( \theta_\lambda \right)\\ \vdots \\ U_k\left( \lambda \right)\sin \left( \theta_\lambda-k\delta\frac{\pi}{2}\right)   \\ \vdots \\ U_{n-1} \left( \lambda \right)\sin \left( \theta_\lambda-(n-1)\delta\frac{\pi}{2}\right)  \end{array}\right) . $$\end{small}
\begin{lemma}  \label{lem:systemelineairecontinuite}
Except for a finite number of values of $R$ in $]\frac{n-1}{2};\frac{n}{2}[$, the matrix $M_R$ is invertible.
\end{lemma}
Proof: We may decompose Chebychev polynomials as
\begin{small}\begin{eqnarray}
U_{2k}(X)=\sum_{j=0}^{k}a_{2k,j} X^{2j},\,  U_{2k+1}(X)=\sum_{j=0}^{k}a_{2k+1,j} X^{2j+1} \mbox{ with }a_{2k,k}=2^{2k} \mbox{ and }a_{2k+1,k}=2^{2k+1}. \label{eqn:exprpolyChe2}
\end{eqnarray}\end{small}We consider $N_R$ the real matrix which is defined by $N_R=(n_{k,j})$ such that, for $0\leq k\leq n-1$:
$$  \left\{\begin{array}{ll}  n_{k,j} =\theta_{n-1}(j)^k  \sin \left( \left[a_{n-1}-\frac{n-2}{2}\right]\theta_{n-1}(j)-\frac{\pi}{2}\left[ j+\delta \frac{n-2k-2}{2} \right]\right) & \mbox{if } 1\leq j \leq \lfloor \frac{n}{2}\rfloor \\ 
 n_{k,j+\lfloor \frac{n}{2}\rfloor } =\theta_n(j)^k  \sin \left( \left[a_{n-1}-\frac{n-1}{2}\right]\theta_{n}(j)-\frac{\pi}{2}\left[ j+\delta \frac{n-2k-1}{2} \right]\right) & \mbox{if } 1\leq j \leq \lfloor \frac{n+1}{2}\rfloor
\end{array} \right.$$
We consider rows of $M_R$ and $N_R$.  Precisely, let
 \begin{small}$$M_R = \left( \begin{array}{c} L_0 \\ \vdots \\ L_{n-1} \end{array} \right) \hspace{3mm} \mbox{and } \hspace{3mm} N_R = \left( \begin{array}{c} \tilde{L}_0 \\ \vdots \\ \tilde{L}_{n-1} \end{array} \right). $$\end{small}Thanks to relation (\ref{eqn:exprpolyChe2}), for $0\leq k \leq \lfloor \frac{n-1}{2}\rfloor$, we have
\begin{small}\begin{eqnarray*}
 L_{2k} = \sum_{j=0}^{k} a_{2k,j}(-1)^{k-i}\tilde{L}_{2j}, 
 \mbox{ and for }0\leq k \leq \lfloor \frac{n-2}{2}\rfloor, \hspace{1mm} L_{2k+1} =\sum_{j=0}^{k} a_{2k+1,j}(-1)^{k-i}\tilde{L}_{2j+1} .
\end{eqnarray*}\end{small}As a result, we have
$$ \det M_R = \left( \prod_{j=0}^{n-1}2^{k} \right) \det N_R =2^{\frac{n(n-1)}{2}}\det N_R .$$ 
In addition, we may write $N_R=(n_{k,j})$ with:
$$  \left\{\begin{array}{ll}  n_{k,j} =\frac{1}{2i}\left([i\delta\theta_{n-1}(j)]^k e^{i(R\theta_{n-1}(j)+\Phi_j)}- [-i\delta\theta_{n-1}(j)]^ke^{-i(R\theta_{n-1}(j)+\Phi_j)} \right) & \mbox{if } 1\leq j \leq \lfloor \frac{n}{2}\rfloor \\ 
 n_{k,j+\lfloor \frac{n}{2}\rfloor } = \frac{1}{2i}\left([i\delta\theta_{n}(j)]^k e^{i(R\theta_{n}(j)+\tilde{\Phi}_j)}- [-i\delta\theta_{n}(j)]^ke^{-i(R\theta_{n}(j)+\tilde{\Phi}_j)} \right) & \mbox{if } 1\leq j \leq \lfloor \frac{n+1}{2}\rfloor
\end{array} \right.$$ 
Due to the multilinearity of determinant, we have
$$\det N_R = \frac{(i\delta)^{\frac{n(n-1)}{2}}}{(2i)^n}\sum_{\epsilon=(\epsilon',\epsilon'')} V[\epsilon',\epsilon''](-1)^{\sigma(\epsilon)}e^{i\left[\sum_{j=1}^{\lfloor \frac{n}{2} \rfloor}\epsilon'(j)\left( R\theta_{n-1}(j) + \Phi_j\right) +\sum_{j=1}^{\lfloor \frac{n+1}{2} \rfloor}\epsilon''(j)\left( R\theta_{n}(j)  -\tilde{\Phi}\right)      \right]}$$
where the sum is running over functions $\epsilon : \{1,..,n\} \longmapsto \{-1,+1\}$ which may split in $\epsilon = (\epsilon',\epsilon'')$ where $\epsilon(j)=\epsilon'(j)$ if $1\leq j \leq \lfloor \frac{n}{2}\rfloor$ and $\epsilon(j+\lfloor \frac{n}{2}\rfloor)=\epsilon''(j)$ if $1\leq j \leq \lfloor \frac{n+1}{2}\rfloor$. In addition, $\sigma(\epsilon)= \left| \left\{1\leq j \leq n ;\, \epsilon(j)=-1\right\}\right|$.
Finally, $V[\epsilon',\epsilon'']$ refers to the determinant of the Vandermonde matrix associated to real numbers $(\epsilon'(1)\theta_{n-1}(1),..,\epsilon'(\lfloor \frac{n}{2}\rfloor)\theta_{n-1}(\lfloor \frac{n}{2}\rfloor), \epsilon''(1)\theta_n(1),..,\epsilon''(\lfloor \frac{n+1}{2}\rfloor)\theta_{n}(\lfloor \frac{n+1}{2}\rfloor))$.
Since  $V[\epsilon',\epsilon'']\neq 0$ for all $\epsilon$, we may conclude $R \longmapsto \det N_R$ vanishes at most a finite number of times on $]\frac{n-1}{2};\frac{n}{2}[$. The same result holds for $R\longmapsto \det M_R$ since $\det M_R  =2^{\frac{n(n-1)}{2}}\det N_R$. 
\begin{flushright}
$\Box$
\end{flushright}

\begin{corollary} \label{cor:expressionrj} For all $R$ in $]\frac{n-1}{2};\frac{n}{2}[$ such that $\det M_R \neq 0$, there exists one and only one continuous function satisfying the differential equation with temporal shifts $(\ref{eqn:eddecalale})$. Moreover:
\begin{itemize}
\item If $1\leq j \leq  \lfloor \frac{n}{2}\rfloor$, 
\begin{eqnarray} \label{eqn:expressionrjn}
r_j(I_{n-2}) = r_\lambda \sum_{k=0}^{n-1}U_k(\lambda) \sin \left(\theta_\lambda-k\delta \frac{\pi}{2}\right)\Delta_{k+1,j}.
\end{eqnarray} 
\item If $1\leq j \leq \lfloor \frac{n+1}{2}\rfloor$, 
\begin{eqnarray} \label{eqn:expressionrjnn}
r_j(I_{n-1}) = -r_\lambda \sum_{k=0}^{n-1}U_k(\lambda) \sin \left(\theta_\lambda-k\delta \frac{\pi}{2}\right)\Delta_{k+1,j+\lfloor \frac{n}{2}\rfloor}.
\end{eqnarray} 
\end{itemize}
with $\Delta_{k,j}=\frac{(-1)^{k+j}}{\det M_R}M_{k,j}$ where  $M_{k,j}$ denotes the minor of $M_R$ obtained by removing from $M_R$ its $k$-th row and $j$-th column.
\end{corollary}

\subsection{Exact value of the minimum} \label{ss:exactvaleur}
In this section, we finish proving theorem \ref{th:thm1}. Remember that some cases have been solved in corollary \ref{cor:remasimplificatrice}.
\paragraph*{}
Thanks to the previous section, we have an explicit expression of the optimal test function $\h_R$. Nevertheless, this explicit expression depends on the unknown parameter $\lambda$ which is related to $\tilde{\m}_R$ by the following relation:
$$\lambda^2 = 4\pi^2 \tilde{\m}_R$$
In order to conclude, we solve the equation $\tilde{\m}_R= \tilde{B}(\h_R)$ where $\tilde{\m}_R$ is the only unknown parameter.
\paragraph*{Two equations, one unknown parameter}
It is technically easier to express $\lambda$ instead of $\tilde{\m}_R$. Thanks to relation (\ref{eqn:defBdeh}), the relation $\tilde{\m}_R= \tilde{B}(\h_R)$ may be written:
\begin{eqnarray} \label{eqn:dernierequationaresoudrepouroptimisation} \lambda^2 = \frac{\int_{\R}\h_R'(u)^2du -\frac{\delta}{2}\int_{-1}^{1}\h_R'*\h_R'(u)du}{\int_{\R}\h_R(u)^2du+\frac{\delta}{2}\int_{-1}^{1}\h_R*\h_R(u)du + \varepsilon \left(\int_\R \h_R(u)du  \right)^2} 
\end{eqnarray}
Furthermore, the relation (\ref{eqn:equationfinalepremiereapparition}) may be written:
\begin{eqnarray} \label{eqn:equationfinalepremiereapparition2}
\frac{w}{\lambda}\cos \lambda R +\frac{\delta}{2}\int_{R-1}^{R}\h_R(x) dx +\varepsilon \int_{-R}^{R}\h_R(x) dx =0
\end{eqnarray}
As a result, $\lambda$ is satisfying two equations.
\paragraph*{}
\begin{lemma}
Relation (\ref{eqn:equationfinalepremiereapparition2}) implies relation (\ref{eqn:dernierequationaresoudrepouroptimisation}).
\end{lemma}
Proof: The differential equation with temporal shifts give
\begin{small}$$\int_{\R}\h_R'(u)^2du -\frac{\delta}{2}\int_{-1}^{1}\h_R'*\h_R'(u)du = \int_{\R} \h_R'(u) \left[\h_R'(u)+\frac{\delta}{2}\int_{u-1}^{u+1}\h_R'(t)dt \right]du = \int_\R \h_R'(u)\varphi'(u)du .
$$ \end{small}Thanks to the Volterra equation we have
\begin{small} 
$$ \int_{\R}\h_R(u)^2du +\frac{\delta}{2}\int_{-1}^{1}\h_R*\h_R(u)du = \int_\R \h_R(u)\varphi(u)du  +\frac{\delta}{2}\left(\int_\R \h_R(u)du\right) \left(\int_{R-1}^{R}\h_R(u)du\right) .$$
\end{small}Finally, using an integration by parts, we get
\begin{small}
$$ \int_\R \h_R'(u)\varphi'(u)du -\lambda^2 \int_\R \h_R(u)\varphi(u)du=-w\lambda \cos\lambda R \int_\R \h_R(u) . $$
\end{small}Therefore, we may write equation (\ref{eqn:dernierequationaresoudrepouroptimisation}) as
\begin{small}
$$ -\lambda^2 \left( \frac{w}{\lambda}\cos \lambda R +\frac{\delta}{2}\int_{R-1}^{R}\h_R (x) dx +\varepsilon \int_{-R}^{R}\h_R (x) dx  \right)\int_\R \h_R(u)=0. $$
\end{small}
\begin{flushright}
$\Box$
\end{flushright}
As a result, we use relation (\ref{eqn:equationfinalepremiereapparition2}) in order to determine $\lambda$.

\paragraph*{End of the proof of theorem \ref{th:thm1}} 

In this subsection, we finish  proving theorem \ref{th:thm1}. Precisely,
\begin{proposition} \label{pro:prodonnantlequationgenerale} If $G=SO^+$, $SO^-$ or $Sp$ with $n\geq 2$, then $\lambda_R :=\lambda $ is the smallest positive root of
\begin{small}
\begin{eqnarray} \label{eqn:expressionimpliciteminimum}
\frac{\delta}{\lambda}\cos\theta_\lambda - \sum_{k=0}^{n-1}U_k(\lambda)\sin\left(\theta_\lambda -k\delta \frac{\pi}{2} \right)\left[\frac{\delta \alpha_R(k)}{2}-1+\varepsilon \beta_R(k)  \right]+\frac{2\varepsilon}{\lambda}\sum_{k=0}^{n-1}U_k(\lambda)\cos\left(\theta_\lambda -k\delta \frac{\pi}{2} \right) = 0
\end{eqnarray}
\end{small}which is not a root of $U_nU_{n-1}$ and where $\alpha_R(k)$ and $\beta_R(k)$ are defined by:
\begin{small}
\begin{eqnarray*}  \alpha_R(k) & = &   2\sum_{j=1}^{\lfloor \frac{n}{2}\rfloor}\frac{\sin \left[\left(R-\frac{n}{2}\right)\theta_{n-1}(j)  \right]\sin\left[\frac{\pi}{2}\left(j+\delta \frac{n-2}{2}\right)\right]}{\theta_{n-1}(j)}\Delta_{k+1,j}  \\
& & \hspace{47mm}+ 2\sum_{j=1}^{\lfloor \frac{n+1}{2}\rfloor}\frac{\sin \left[\left(R-\frac{n-1}{2}\right)\theta_{n}(j)  \right]\sin\left[\frac{\pi}{2}\left(j+\delta \frac{n-1}{2}\right)\right]}{\theta_{n}(j)}\Delta_{k+1,j+\lfloor \frac{n}{2}\rfloor}  \end{eqnarray*}
\end{small}and
\begin{small}
\begin{eqnarray*}  \beta_R(k) & = &   2\sum_{j=1}^{\lfloor \frac{n}{2}\rfloor}\frac{\sin \left[\left(R-\frac{n}{2}\right)\theta_{n-1}(j)  \right]}{\theta_{n-1}(j)}\Delta_{k+1,j}\sum_{l=0}^{n-2}U_l(\theta_{n-1}(j)\sin\left[\frac{\pi}{2}\left(j+\delta \frac{n-2l-2}{2}\right)\right]  \\
& & \hspace{11mm} + 2\sum_{j=1}^{\lfloor \frac{n+1}{2}\rfloor}\frac{\sin \left[\left(R-\frac{n-1}{2}\right)\theta_{n}(j)  \right]}{\theta_{n}(j)}\Delta_{k+1,j+\lfloor \frac{n}{2}\rfloor}\sum_{l=0}^{n-1}U_l(\theta_n(j)) \sin\left[\frac{\pi}{2}\left(j+\delta \frac{n-2l-1}{2}\right)\right] . \end{eqnarray*}
\end{small}
\end{proposition}
Before proving this result, we need to prove some technical lemmas. We may write $\h_R$ as a sum of two non-continuous even functions. Let $\varphi_\lambda$ and $\psi$ which are defined on $\R \backslash\{a_{-n},...,a_n\}$ by:
\begin{itemize} \item
If $0\leq k \leq n-1$, let
\begin{eqnarray*}
\varphi_\lambda |_{I_{n-(2k+1)}}(u) =  r_n^\lambda(k) \sin \left(  \lambda \left[ u-\frac{n-2k-1}{2} \right]+\theta_n^\lambda(k) \right).
\end{eqnarray*}
\item If $0\leq k \leq n-2$, let
\begin{eqnarray*}
\varphi_\lambda |_{I_{n-2(k+1)}}(u) =  r_{n-1}^\lambda(k)\sin \left(  \lambda \left[ u-\frac{n-2k-2}{2} \right]+\theta_{n-1}^\lambda(k) \right).
\end{eqnarray*}
\item supp $\varphi_\lambda \subset [-R,R]$.
\item If $u \in \R \backslash\{a_{-n},...,a_n\}$, let
$$\psi (u) =  \h_R(u)-\varphi_\lambda(u) . $$
\end{itemize}
Even though $\varphi_\lambda$ and $\psi$ are not  continuous functions, they are smooths on each $I_k$. Furthermore, they satisfy
  $$  \h_R(u)=\psi (u) +\varphi_\lambda(u) \hspace{10mm} \mbox{and} \hspace{8mm}\varphi_\lambda '' = -\lambda^2 \varphi_\lambda .$$
\begin{lemma} \label{lem:integralespsi} We have
\begin{eqnarray*}
\int_{R-1}^{R} \psi (u) du & = &   r_\lambda \sum_{k=0}^{n-1}U_k(\lambda)\sin\left( \theta_\lambda-k\delta\frac{\pi}{2}\right) \alpha_R(k)\\
\mbox{and } \hspace{3mm} \int_{-R}^{R} \psi (u) du & = & r_\lambda \sum_{k=0}^{n-1}U_k(\lambda)\sin\left( \theta_\lambda-k\delta\frac{\pi}{2}\right) \beta_R(k) .
\end{eqnarray*}
\end{lemma}
Proof: We may write:
\begin{small}
\begin{eqnarray*}
\int_{R-1}^{R} \psi (u) du & = & \int_{I_{n-1}} \psi |_{I_{n-1}}(u)du + \int_{I_{n-2}} \psi |_{I_{n-2}}(u)du \\
& = & 2\sum_{j=1}^{\lfloor \frac{n}{2}\rfloor}r_j(I_{n-2})\frac{\sin \left[\left(R-\frac{n}{2}\right)\theta_{n-1}(j)  \right]}{\theta_{n-1}(j)}\sin\left[\frac{\pi}{2}\left(j+\delta \frac{n-2}{2}\right)\right] \\ 
& & \hspace{35mm} -2\sum_{j=1}^{\lfloor \frac{n+1}{2}\rfloor}r_j(I_{n-1})\frac{\sin \left[\left(R-\frac{n-1}{2}\right)\theta_{n}(j)  \right]}{\theta_{n}(j)}\sin\left[\frac{\pi}{2}\left(j+\delta \frac{n-1}{2}\right)\right] 
\end{eqnarray*}
\end{small}The result comes easily from relations (\ref{eqn:expressionrjn}) and (\ref{eqn:expressionrjnn}) of corollary \ref{cor:expressionrj}. 
Similarly, since
\begin{small}$$\int_{-R}^{R} \psi (u) du = \sum_{k=0}^{n-1} \int_{I_{n-1}}\psi|_{I_{n-2k-1}}(u-k)du + \sum_{k=0}^{n-2} \int_{I_{n-2}}\psi|_{I_{n-2k-2}}(u-k)du$$\end{small}we obtain the second part of this lemma.
\begin{flushright}
$\Box$
\end{flushright}
\begin{lemma} For $0\leq k \leq n-2$, we have
$$ \lim_{u\rightarrow \sup I_{n-2k-2} }\varphi_\lambda'|_{I_{n-2k-2}}(u)- \lim_{u\rightarrow \inf I_{n-2k-1} }\varphi_\lambda'|_{I_{n-2k-1}}(u) = \lambda U_k(\lambda)\emph{Re }\left( (i\delta)^{-k} z_\lambda \right) $$
and for $1\leq k \leq n-1$, 
$$ \lim_{u\rightarrow \sup I_{n-2k-1} }\varphi_\lambda'|_{I_{n-2k-1}}(u)- \lim_{u\rightarrow \inf I_{n-2k} }\varphi_\lambda'|_{I_{n-2k}}(u)= \lambda U_{n-k-1}(\lambda)\emph{Re }\left( (i\delta)^{-(n-k-1)} z_\lambda \right) $$
where  the complex number $z_\lambda$ is defined in lemma \ref{lem:defzlambda}.
\end{lemma}
Proof: Since $\h_R$ satisfies the differential equation with temporal shifts (\ref{eqn:eddecalale}), we may conclude $\h_R'$ is continuous on $]-R,R[$. Therefore, for  $0\leq k \leq n-2$, we get
\begin{small}\begin{eqnarray*}
& & \lim_{u\rightarrow \sup I_{n-2k-2} }\varphi_\lambda'|_{I_{n-2k-2}}(u)- \lim_{u\rightarrow \inf I_{n-2k-1} }\varphi_\lambda'|_{I_{n-2k-1}}(u)\\
 & = & \lambda \mbox{Re }  \left[ e^{i\lambda \left( a_{n-1} -\frac{n-2}{2} \right)}z_{n-1}^\lambda(k) - e^{i\lambda \left( a_{n-1} -\frac{n-1}{2} \right)}z_{n}^\lambda(k)\right] 
\end{eqnarray*}\end{small}and, for $1\leq k \leq n-1$, 
\begin{small}$$\lim_{u\rightarrow \sup I_{n-2k-1} }\varphi_\lambda'|_{I_{n-2k-1}}(u)- \lim_{u\rightarrow \inf I_{n-2k} }\varphi_\lambda'|_{I_{n-2k}}(u) = \lambda \mbox{Re }  \left[ e^{i\lambda \left( a_{n} -\frac{n-1}{2} \right)}z_{n}^\lambda(k) - e^{i\lambda \left( a_{n} -\frac{n}{2} \right)}z_{n-1}^\lambda(k-1)\right] .
$$\end{small}The result comes easily from relations (\ref{eqn:defznlambdak}) and  (\ref{eqn:relationmultiplicativitedesfibo}). 
\begin{flushright}
$\Box$
\end{flushright}
\begin{lemma} \label{lem:integralesphilambda} We have 
\begin{eqnarray*}
\int_\R \varphi_\lambda(u)du & = & -\frac{2 r_\lambda}{\lambda}\sum_{k=0}^{n-1} U_k(\lambda)\cos \left( \theta_\lambda-k\delta \frac{\pi}{2}\right) \\
\mbox{and } \hspace{3mm}  \int_{R-1}^{R} \varphi_\lambda(u)du &  = & -\frac{2w}{\delta \lambda} \cos \lambda R -\frac{2}{\lambda} r_\lambda \cos \theta_\lambda + 2\delta \emph{Re }\left[  i z_\lambda \sum_{k=0}^{n-1}(-i\delta)^kU_k(\lambda) \right] .
\end{eqnarray*}
\end{lemma}
Proof: We may write:
\begin{small}
\begin{eqnarray*}
\int_\R \varphi_\lambda(u)du & = &-\frac{1}{\lambda^2}\int_\R \varphi_\lambda ''(u) du =-\frac{1}{\lambda^2}\sum_{k=-(n-1)}^{n-1}[\varphi_\lambda '|_{I_k}]_{\inf I_k}^{\sup I_k} \\
 & = & -\frac{1}{\lambda^2} \left[ 2\lim_{u\rightarrow \sup I_{n-1} }\varphi_\lambda'|_{I_{n-1}}(u) +\sum_{k=-(n-1)}^{n-2} \lim_{u\rightarrow \sup I_{k} }\varphi_\lambda'|_{I_{k}}(u)- \lim_{u\rightarrow \inf I_{k+1} }\varphi_\lambda'|_{I_{k+1}}(u) \right]
\end{eqnarray*}
\end{small}Since $z_n^\lambda (0)= (i\delta)^{-(n-1)}e^{-i\lambda \left(a_n- \frac{n-1}{2}\right)}U_{n-1}(\lambda)z_\lambda$, we get
\begin{small}$$\lim_{u\rightarrow \sup I_{n-1} }\varphi_\lambda'|_{I_{n-1}}(u) = -\lambda \mbox{Re }\left(z_\lambda (i\delta)^{-(n-1)}U_{n-1}(\lambda)  \right)$$\end{small}and the previous lemma gives
\begin{small}$$ \int_\R \varphi_\lambda(u)du =\frac{-2}{\lambda}\sum_{k=0}^{n-1}U_k(\lambda) \mbox{Re }\left((i\delta)^{-k}z_\lambda  \right) . $$\end{small}Similarly, we have:
\begin{small}
\begin{eqnarray*}
& & \int_{R-1}^{R} \varphi_\lambda(u)du \\
 & & =\frac{-1}{\lambda^2} \left[\lim_{u \rightarrow \sup I_{n-1}} \varphi'|_{I_{n-1}}(u) + \lim_{u \rightarrow \sup I_{n-2}} \varphi'|_{I_{n-2}}(u)-\lim_{u \rightarrow \inf I_{n-1}} \varphi'|_{I_{n-1}}(u)-\lim_{u \rightarrow \inf I_{n-2}} \varphi'|_{I_{n-2}}(u) \right] 
\end{eqnarray*}
\end{small}The previous lemma gives 
\begin{small}$$\lim_{u \rightarrow \sup I_{n-2}} \varphi'|_{I_{n-2}}(u)-\lim_{u \rightarrow \inf I_{n-1}} \varphi'|_{I_{n-1}}(u)=\lambda \mbox{Re }(z_\lambda) . $$\end{small}In addition, since $z_{n-1}^\lambda (0) = i\delta \frac{U_n(\lambda)}{U_{n-1}(\lambda)}e^{i\lambda/2}z_n^\lambda (0)-2w\delta\lambda \mbox{Re }(e^{i\lambda R})$, we get
\begin{small}$$\lim_{u\rightarrow \inf I_{n-2} }\varphi_\lambda'|_{I_{n-2}}(u) =  -2w\delta \lambda \cos (\lambda R)-\lambda \mbox{Re } \left[U_n(\lambda) z_\lambda (i\delta)^{-n}  \right] . $$\end{small}Therefore, we may write
\begin{small}
\begin{eqnarray*}
\int_{R-1}^{R} \varphi_\lambda(u)du = -\frac{2w}{\delta \lambda} \cos \lambda R -\frac{2}{\lambda} r_\lambda \cos \theta_\lambda + \frac{1}{\lambda} \mbox{Re }\left[ z_\lambda (1-(i\delta)^{-n}U_n(\lambda)+(i\delta)^{-(n+1)}U_{n-1}(\lambda) )\right] .
\end{eqnarray*}
\end{small}The result comes easily from relation (\ref{eqn:seriegeneratricetronquee}).
\begin{flushright}
$\Box$
\end{flushright}
To conclude, since $\h_R= \psi+\varphi_\lambda$, we may easily transform equation (\ref{eqn:equationfinalepremiereapparition2}) thanks to both lemmas \ref{lem:integralespsi} and \ref{lem:integralesphilambda}. Proposition \ref{pro:prodonnantlequationgenerale} follows immediately due to the fact that $r_\lambda \neq 0$ (otherwise $\h_R=0$).

\paragraph*{Particular case}
In this subsection, we are assuming $n=2$ and we give a simpler expression of equation (\ref{eqn:expressionimpliciteminimum}) than in  proposition \ref{pro:prodonnantlequationgenerale}.
\begin{corollary} If $G=SO^+$, $SO^-$ or $Sp$ and $\frac{1}{2}<R<1$ then $\lambda_R$ is the smallest positive root of
\begin{small}
\begin{eqnarray*}& &(\delta+2\varepsilon)\frac{1-4\lambda^2}{\lambda}\left(\sin\lambda(1-R)-2\delta\lambda \cos\lambda R\right)\\ & & -\left[(\delta+2\varepsilon)(1-R)-1+4\varepsilon \right]\left[\cos\lambda (1-R)-2\delta \lambda \sin\lambda R-2\lambda \tan \Theta_R (\sin\lambda(1-R)-2\delta\lambda \cos\lambda R)  \right]=0
\end{eqnarray*}
\end{small}which is not a root of $U_1U_2$ and where
$$\Theta_R = \frac{1}{2}\left(R-\frac{1}{2}\right)+\frac{\pi}{2}\left(1+\frac{\delta}{2} \right) . $$
\end{corollary}
Proof: The matrix $M_R$ which appears in corollary \ref{cor:expressionrj} can be written
$$ M_R = \left( \begin{array}{cc}  -1 & \sin\left[\frac{1}{2}\left(\frac{1}{2}-R\right)-\frac{\pi}{2}\left(1+\frac{\delta}{2}\right)\right] \\
0 & \sin\left[\frac{1}{2}\left(\frac{1}{2}-R\right)-\frac{\pi}{2}\left(1-\frac{\delta}{2}\right)\right] \end{array} \right) . $$
Therefore, 
$$\det M_R = -\delta \cos \Theta_R \neq 0 .$$
We get $\Delta_{1,1}=-1$, $\Delta_{1,2}=0$, $\Delta_{2,1}=-\delta \tan \Theta_R$  and $\Delta_{2,2}=\frac{\delta}{\cos \Theta_R}$.
In addition, if $k\in\left\{1,2\right\}$, we have:
\begin{eqnarray*}
\alpha_R(k) & = &  2(R-1)\Delta_{k+1,1}-2\left(\cos\Theta_R+\delta \sin \Theta_R \right) \Delta_{k+1,2} \\
\beta_R(k) & = & \alpha_R(k)+2\left(\cos\Theta_R+\delta \sin \Theta_R \right) \Delta_{k+1,2} 
\end{eqnarray*}
Therefore:
$$\left\{ \begin{array}{l} \alpha_R(0)=\beta_R(0)=2(1-R) \\ \alpha_R(1)=2\delta(1-R)-2(\delta+\tan \Theta_R)\\ \beta_R(1)=2\delta(1-R)-4(\delta+\tan \Theta_R) \end{array} \right. $$
As a result, since $\delta \varepsilon=-\varepsilon$, we get
\begin{small}
\begin{eqnarray*}
& &\frac{\delta}{\lambda}\cos\theta_\lambda - \sum_{k=0}^{n-1}U_k(\lambda)\sin\left(\theta_\lambda -k\delta \frac{\pi}{2} \right)\left[\frac{\delta \alpha_R(k)}{2}-1+\varepsilon \beta_R(k)  \right]+\frac{2\varepsilon}{\lambda}\sum_{k=0}^{n-1}U_k(\lambda)\cos\left(\theta_\lambda -k\delta \frac{\pi}{2} \right)\\
& = & \frac{\delta+2\varepsilon}{\lambda}\cos\theta_\lambda- \left(\frac{\delta \alpha_R(0)}{2}+\varepsilon\beta_R(0)-1+4\varepsilon \right)\sin \theta_\lambda+2\lambda \delta  \left( \frac{\delta \alpha_R(1)}{2}+\varepsilon\beta_R(1)-1\right)\cos\theta_\lambda \\
& = & (\delta+2\varepsilon)\frac{1-4\lambda^2}{\lambda}\cos\theta_\lambda -\left[(\delta+2\varepsilon)(1-R)-1+4\varepsilon \right]\left[ \sin \theta_\lambda-2\lambda \tan (\Theta_R) \cos\theta_\lambda \right].
\end{eqnarray*}
\end{small}In addition, we have
\begin{small}\begin{eqnarray*} z_\lambda & = & \frac{-2iwe^{-i\lambda(1-R)}}{U_1(\lambda)U_2(\lambda)}\left(1+i\delta U_1(\lambda)e^{i\lambda} \right)\\
 & = &  \frac{-2w}{U_1(\lambda)U_2(\lambda)}\left(\sin\lambda(1-R)-2\delta\lambda \cos\lambda R+i\left[\cos\lambda (1-R)-2\delta \lambda \sin\lambda R \right]  \right) .
 \end{eqnarray*}\end{small}The result comes easily from these relations.
\begin{flushright}
$\Box$
\end{flushright} 
  
\paragraph*{Comments on the ``$Sp$ hypothesis''}
In the symplectic case with $R>1/2$, our proof of theorem \ref{th:thm1} is submitted to the ``$Sp$ hypothesis''. However, we prove that this theorem is still true even though $\m_R^2$ is an odd positive integer. Therefore, throughout this paragraph, we are assuming $G=Sp$, $R>1/2$ and $m_R = N^2$ where $N$ is an odd integer.  
\begin{lemma} \label{lem:pourref} $\h_R$ satisfies $\int_{R-1}^{R}\h_R(u)du=0$ and, for all $0\leq u\leq R$, 
$$\h_R(u)=\varphi(u) +\frac{\delta}{2}\int_{u}^{R}\h_R(t+1)-\h_R(t-1) dt$$
with
\begin{small} 
$$\varphi(u)=\left(\frac{-4}{N\pi}\int_{0}^{R}\left[ \h_R'(t)-\frac{\delta}{2}\left(\h_R(t+1)-\h_R(t-1)\right) \right]\sin\left(\frac{\pi N t}{2R} \right)dt\right) (\cos\left(\lambda u \right)-\cos(\lambda R) )\cdot \un_{[-R,R]}(u) . $$
\end{small}
\end{lemma}
Proof: In lemma \ref{lem:pour Sphypo}, we prove that  for all odd positive integer $n$,
\begin{small}$$(n^2-\m_R)\c_n =  \frac{ (-1)^\frac{n-1}{2}k_{\delta,\varepsilon}}{n}-(n^2-\m_R)\frac{\delta}{n\pi}(S_n*\h_R)(1) .$$\end{small}First, with $n=N$, we may deduce $k_{\delta,\varepsilon}=\int_{R-1}^R\h_R(u)du=0$.  Second, for all $n\neq N$, we get
\begin{small}$$\c_n = -\frac{\delta}{n\pi}(S_n*\h_R)(1) .  $$\end{small}Since $\h_R(u)=2\sum_{n\geq 0}\c_n \cos\left(\frac{\pi n u}{2R} \right)\cdot \un_{[-R,R]}(u)$, for all $u$ in $[0,R]$, we write
\begin{small}$$\h_R(u)=2\left[\c_N+\frac{\delta}{N\pi}S_N*\h_R(1) \right]\cos\left(\frac{\pi N u}{2R} \right)-\frac{2\delta}{\pi}{\sum_{n\geq 0}}^* \frac{S_n*\h_R(1)}{n}\cos\left(\frac{\pi n u}{2R} \right).$$\end{small}The sum in the right member of this equality has been computed in the proof of lemma \ref{lem:equationvolterra} and we have
\begin{small}$$ \c_N+\frac{\delta}{N\pi}S_N*\h_R(1) = \frac{-2}{N\pi}\int_{0}^{R}\left[ \h_R'(t)-\frac{\delta}{2}\left(\h_R(t+1)-\h_R(t-1)\right) \right]\sin\left(\frac{\pi N t}{2R} \right)dt . $$\end{small}
\begin{flushright}
$\Box$
\end{flushright}
Therefore, changing $w$ with
\begin{small}$$w=\frac{4 \lambda }{N\pi}\int_{0}^{R}\left[ \h_R'(t)-\frac{\delta}{2}\left(\h_R(t+1)-\h_R(t-1)\right) \right]\sin\left(\frac{\pi N t}{2R} \right)dt,  $$\end{small}this Volterra equation with temporal shift has been solved in section \ref{ss:optimaltestfunction}. Thus, we get an explicit expression of $\h_R$. Now, several cases may occur. First, this explicit expression of $\h_R$ doesn't satisfy the compatibility equation $\int_{R-1}^{R}\h_R(u)du=0$, then the ``$Sp$ hypothesis'' is true. Second, if $h_R$ satisfies the compatibility equation, then on account of the fact that the argument of the complex number $z_\lambda$ is independent of $w$, $\lambda$ is a root of equation (\ref{eqn:expressionimpliciteminimum}). Since $\m_R$ is the smallest critical value of $B$, $\lambda$ is still the smallest root of equation (\ref{eqn:expressionimpliciteminimum}).

\nocite{M,RR2,I}

\end{document}